\documentclass[11pt]{article}

\usepackage{amsmath,amssymb,amsfonts,enumitem,amsthm,accents,pstricks,pst-node,pstricks-add,mathrsfs,xy}
\usepackage{pst-solides3d}

\xyoption{all}

\usepackage{color}
\usepackage{makeidx}
\usepackage{graphicx}
\usepackage{setspace}
\usepackage{extarrows}
\usepackage{hyperref}

\numberwithin{equation}{section}

\usepackage{tocloft}

\usepackage[Symbolsmallscale]{upgreek} 

\newcommand*{\defeq}{\mathrel{\vcenter{\baselineskip0.5ex \lineskiplimit0pt
                     \hbox{\scriptsize.}\hbox{\scriptsize.}}}%
                     =}
                     
\newcommand{\sslash}{\!\mathbin{/\mkern-6mu/\!}}

\newtheorem{theorem}{Theorem}[section]
\newtheorem{lemma}[theorem]{Lemma}
\newtheorem{conjlemma}[theorem]{Conjectural Lemma}
\newtheorem{corollary}[theorem]{Corollary}
\newtheorem{proposition}[theorem]{Proposition}
\newtheorem{conjecture}[theorem]{Conjecture}
\newtheorem{question}[theorem]{Question}

\theoremstyle{definition}
\newtheorem{definition}[theorem]{Definition}
\newtheorem{remark}[theorem]{Remark}
\newtheorem{example}[theorem]{Example}

\def\ov#1{\overline{#1}}
\def\tn#1{\textnormal{#1}}
\def\mf#1{\mathfrak{#1}}
\def\wt#1{\widetilde{#1}}

\def\vr{\varrho}
\def\ll{\left\langle}
\def\rr{\right\rangle}
\def\mc{\mathcal}
\def\lra{\longrightarrow}

\def\ve{\varepsilon}

\newcommand{\abs}[1]{\left\vert #1 \right\vert}
\newcommand{\lrp}[1]{\left( #1 \right)}
\newcommand{\lrc}[1]{\left\{ #1 \right\}}

\def\bEq#1{\begin{equation}\label{#1}}
\def\eEq{\end{equation}}

\def\bsEq{\begin{equation*}}
\def\esEq{\end{equation*}}

\def\bDf#1{\begin{definition}\label{#1}}
\def\eDf{\end{definition}}

\def\bTh#1{\begin{theorem}\label{#1}}
\def\eTh{\end{theorem}}

\def\bCn#1{\begin{conjecture}\label{#1}}
\def\eCn{\end{conjecture}}

\def\bLm#1{\begin{lemma}\label{#1}}
\def\eLm{\end{lemma}}

\def\bCLm#1{\begin{conjlemma}\label{#1}}
\def\eCLm{\end{conjlemma}}

\def\bRm#1{\begin{remark}\label{#1}}
\def\eRm{\end{remark}}

\def\bEx#1{\begin{example}\label{#1}}
\def\eEx{\end{example}}

\def\bPr#1{\begin{proposition}\label{#1}}
\def\ePr{\end{proposition}}

\def\bCr#1{\begin{corollary}\label{#1}}
\def\eCr{\end{corollary}}

\def\bFg#1{\begin{figure}\label{#1}}
\def\eFg{\end{figure}}

\def\bPf{\begin{proof}}
\def\ePf{\end{proof}}

\def\bIt{\begin{itemize}[leftmargin=*]}
\def\eIt{\end{itemize}}

\def\bEn{\begin{enumerate}[label=$(\arabic*)$,leftmargin=*]}
\def\eEn{\end{enumerate}}

\def\bEnalph{\begin{enumerate}[label=$(\alph*)$,leftmargin=*]}
\def\eEnalph{\end{enumerate}}

\def\SL{\tn{SL}}
\def\SU{\tn{SU}}
\def\PGL{\P\tn{GL}}

\def\Hom{\tn{Hom}}

\def\cA{\mc{A}}
\def\cB{\mc{B}}
\def\cC{\mc{C}}
\def\cD{\mc{D}}
\def\cF{\mc{F}}
\def\cI{\mc{I}}

\def\cM{\mc{M}}
\def\cO{\mc{O}}

\def\cP{\mc{P}}

\def\cN{\mc{N}}
\def\cR{\mc{R}}

\def\cL{\mc{L}}

\def\cC{\mc{C}}

\def\cT{\mc{T}}

\def\R{\mathbb R}
\def\C{\mathbb C}
\def\Z{\mathbb Z}
\def\Q{\mathbb Q}
\def\P{\mathbb P}
\def\H{\mathbb H}
\def\T{\mathbb T}
\def\N{\mathbb N}

\def\mfi{\mf{i}}

\def\mfs{\mf{s}}
\def\mft{\mf{t}}

\def\la{\lambda}
\def\La{\Lambda}

\def\De{\Delta}
\def\de{\delta}

\def\si{\sigma}
\def\Si{\Sigma}
\def\al{\alpha}
\def\be{\beta}

\def\eset{\emptyset}

\def\G{\ifmmode{\cal G}\else{${\cal G}$}\fi}        
\def\spin{\ifmmode{\mathrm{spin}}\else{spin}\fi}
\def\Spin{\ifmmode{\mathrm{Spin}}\else{Spin}\fi}
\def\spinc{\ifmmode{\mathrm{spin}^c}\else{spin$^c$ }\fi}
\def\Spinc{\ifmmode{\mathrm{Spin}^c}\else{Spin$^c$ }\fi}

\def\sfrak{\ifmmode{\mathfrak s}\else{${\mathfrak s }$}\fi}  
\def\tfrak{\ifmmode{\mathfrak t}\else{${\mathfrak t }$}\fi} 
\def\sz{\ifmmode{{\tn\ss}}\else{\ss}\fi} 
\def\sbold{\ifmmode{\mbox{\boldmath$s$\unboldmath}}\else{\boldmath$s$\unboldmath}\fi}
\def\tbold{\ifmmode{\mbox{\boldmath$t$\unboldmath}}\else{\boldmath$t$\unboldmath}\fi} 

\usepackage{pict2e}

\makeatletter

\newcommand*\@KP@Large@frame[2]{%
    \setlength\unitlength{\fontdimen 22 #1\tw@}%
    \vrule \@width\z@ \@height 4\unitlength \@depth\tw@\unitlength
    \begin{picture}(6,2)(-3,-1)%
        \def\@KP@Radius     {3}%
        \def\@KP@Hole@radius{.5}
        \def\@KP@Diameter   {6}%
        #2%
    \end{picture}%
}
\newcommand*\@KP@Small@frame[2]{%
    \setlength\unitlength{\fontdimen 22 #1\tw@}%
    \vrule \@width\z@ \@height \thr@@\unitlength \@depth\@ne\unitlength
    \begin{picture}(4,2)(-2,-1)%
        \def\@KP@Radius     {2}%
        \def\@KP@Hole@radius{.5}
        \def\@KP@Diameter   {4}%
        #2%
    \end{picture}%
}

\newcommand*\@KP@Radius     {}
\newcommand*\@KP@Hole@radius{}
\newcommand*\@KP@Diameter   {}
\newcommand*\@KP@Shape@A{%
    \put(0,0){\circle{\@KP@Diameter}}%
}
\newcommand*\@KP@Shape@B{%
    \Line(-\@KP@Radius,\@KP@Radius )(\@KP@Radius,-\@KP@Radius)%
    \Line(-\@KP@Radius,-\@KP@Radius)(-\@KP@Hole@radius,-\@KP@Hole@radius)%
    \Line(\@KP@Radius ,\@KP@Radius )(\@KP@Hole@radius ,\@KP@Hole@radius )%
}
\newcommand*\@KP@Shape@E{%
    \Line(-\@KP@Radius,\@KP@Radius )(\@KP@Radius,-\@KP@Radius)%
     \Line(\@KP@Radius,\@KP@Radius )(-\@KP@Radius,-\@KP@Radius)%

}

\newcommand*\@KP@Shape@C{%
    \cbezier(-\@KP@Radius,\@KP@Radius )(0,0)(0,0)(\@KP@Radius,\@KP@Radius )%
    \cbezier(-\@KP@Radius,-\@KP@Radius)(0,0)(0,0)(\@KP@Radius,-\@KP@Radius)%
}
\newcommand*\@KP@Shape@D{%
    \cbezier(-\@KP@Radius,-\@KP@Radius)(0,0)(0,0)(-\@KP@Radius,\@KP@Radius)%
    \cbezier(\@KP@Radius ,-\@KP@Radius)(0,0)(0,0)(\@KP@Radius ,\@KP@Radius)%
}

\newcommand*\@KP@Atomic@mathpalette[1]{%
    \mathinner{
        \mathchoice{%
            \linethickness{.6\p@}
            \@KP@Large@frame \textfont {#1}%
        }{%
            \linethickness{.4\p@}
            \@KP@Small@frame \textfont {#1}%
        }{%
            \linethickness{.3\p@}
            \@KP@Small@frame \scriptfont {#1}%
        }{%
            \linethickness{.2\p@}
            \@KP@Small@frame \scriptscriptfont {#1}%
        }%
    }%
}

\newcommand*\KPA{\@KP@Atomic@mathpalette \@KP@Shape@A}
\newcommand*\KPB{\@KP@Atomic@mathpalette \@KP@Shape@B}
\newcommand*\KPC{\@KP@Atomic@mathpalette \@KP@Shape@C}
\newcommand*\KPD{\@KP@Atomic@mathpalette \@KP@Shape@D}
\newcommand*\KPE{\@KP@Atomic@mathpalette \@KP@Shape@E}

\makeatother

\makeindex
\begin{document}
\title{On compactifications of the $\SL(2,\C)$ character varieties\\ of punctured surfaces}
\author{Mohammad Farajzadeh-Tehrani\footnote{Partially supported by the NSF grant DMS-2003340}~~and Charles Frohman}
\date{\today}
\maketitle

\begin{abstract}
This paper addresses some conjectures and questions regarding the absolute and relative compactifications of the $\SL(2,\C)$-character variety of an $n$-punctured Riemann surface without boundary. We study a class of projective compactifications determined by ideal triangulations of the surface and prove explicit results concerning the boundary divisors of these compactifications. 
Notably, we establish that the boundary divisors are toric varieties and confirm a well-known conjecture asserting that the (dual) boundary complex of any (positive dimensional) relative character variety is a sphere.
In a different vein, we enhance and streamline Komyo's compactification method, which utilizes a projective compactification of $\SL(2,\C)$ to compactify the (relative) character varieties. Specifically, we construct a uniform relative compactification over the base space of $\C^n$ and determine its monodromy, addressing a question posed by Simpson.
\end{abstract}
\tableofcontents

\section{Introduction}\label{overview}

For $g,n\!\geq\! 0$, let $\Si_{g,n}$ denote an $n$-punctured Riemann surface of genus $g$ (without boundary).
Given an algebraic reductive Lie group $G$, the $G$-character variety of $\Si_{g,n}$ is the GIT quotient
\bEq{GIT-cR_e}
\cR_{g,n}(G)\defeq  \Hom\big(\pi_1(\Si_{g,n}),G\big)\sslash G,
\eEq
where $\Hom\big(\pi_1(\Si_{g,n}),G \big)$ is the affine variety of representations of the fundamental group $\pi_1(\Si_{g,n})$ into $G$, and $G$ acts by post-compositions with inner automorphisms of $G$; i.e., by
$$
\rho \lra \vr^{-1} \rho \vr \qquad \forall~\rho\in \Hom\big(\pi_1(\Si_{g,n}),G \big),~\vr\in G.
$$
In this paper, we are concerned with the case of  $n\!>\!0$ and $G=\SL(2,\C)$ (as well as the subgroups $\SL(2,\R),\SU(2)\subset \SL(2,\C)$ in Section~\ref{real}). For $G=\SL(2,\C)$, we will use the simplified notation $\cR_{g,n}\defeq\cR_{g,n}(\SL(2,\C))$. If $g=0$, we write $\cR_n$ in place of  $\cR_{g,n}$. We will follow the same notation convention throughout the paper.\\

Let $a_1, \ldots,a_n$ denote the peripheral loops around the punctures. For every $[\rho]\in \cR_{g,n}$, the trace functions $\tn{tr}(\rho(a_i))\in \C$ are well-defined and define an algebraic projection map
\bEq{pin_e}
\Pi_{g,n}\colon  \cR_{g,n} \lra~\C^n, \qquad [\rho] \lra  \mft=\big(t_i=\tn{tr}(\rho(a_i))\big)_{i=1}^n.
\eEq
In this paper\footnote{Often, relative character varieties are defined by fixing the conjugacy classes at punctures, and that is slightly different than the definition considered here when $t_i=\pm 2$ for some $i=1,\ldots,n$.}, the fibers $\cR_{g,n,\mft}=\Pi_{g,n}^{-1}(\mft)$ are called the \textbf{relative} $\SL(2,\C)$-character varieties of $\Si_{g,n}$, or the \textbf{fibers}\footnote{In the literature they are also called the {\bf slices} of $\cR_{g,n}$.} of $\cR_{g,n}$. The entire character variety $\cR_{g,n}$ as well as each fiber $\cR_{g,n,\mft}$ are quasi-projective affine varieties of complex dimensions $3(2g+n-2)$ and $2(3g+n-3)$, respectively. In this paper, we study divisorial compactifications of $\cR_{g,n}$ and $\cR_{g,n,\mft}$. 

\subsection{Geometric P=W conjecture and the main results}
It is well known that the character varieties $\cR_{g,n}(G)$ are closely related to several rich geometric and algebraic structures, including the Dolbeault (or Hitchin) moduli spaces of principal $G$-bundles \cite{H,S2}, the skein algebras of  surfaces \cite{BFK,PS,T}, and cluster varieties \cite{FG}. For instance, the isomorphism with Dolbeault moduli spaces implies that $\cR_{g}(G) = \cR_{g,0}(G)$ is hyperk\"ahler and admits a Lagrangian torus fibration.\\

A major open question in this area is the geometric P=W conjecture (stated below), which compares the topology of compactifications of $\cR_{g}(G)$ with those of the corresponding Hitchin moduli spaces $\cM_{g}(G)$. In $\cM_{g}(G)$, the points parametrize flat principal $G$-Higgs bundles. In the punctured case, the moduli spaces $\cM_{g,n}^{\mathrm{rel}}(G)$ are defined analogously, by imposing parabolic structures and additional conditions (denoted here simply by $\mathrm{rel}$) on the conjugacy classes of monodromies around the punctures.

There exists a canonical real analytic identification $\Psi\colon \cM_{g}(G) \lra \cR_{g}(G)$, known as the non-abelian Hodge correspondence (NAHC). This result stems from the foundational work of many mathematicians in the 1980s, notably Hitchin \cite{H}, Simpson \cite{S88}, Corlette \cite{C88}, and Donaldson \cite{D87}. A gentle review of the main ideas and the proof can be found in \cite[Theorem~8.7]{Thomas}. Relative versions of $\Psi$ between  $\cM_{g,n}^{\mathrm{rel}}(G)$ and slices of $\cR_{g,n}(G)$ can also be defined. Simpson \cite{S2} appears to have been the first to extend non-abelian Hodge theory to character varieties of punctured surfaces with arbitrary conjugacy classes. His work addresses the ``tame" case, which was subsequently generalized to the ``wild" case in \cite{BB}. Among more recent examples, \cite{ES} explicitly discusses non-Abelian Hodge correspondence (NAHC) in the punctured setting, and \cite{HHT} has results concerning the four-punctured sphere. \vskip.05in

The geometric P=W conjecture, as proposed by Katzarkov et. al. \cite[Conjecture~1.1]{KNPS}, and further refined in \cite{MMS},  aims to provide a deeper understanding of the asymptotic behavior of the transcendental map $\Psi$ with respect to suitable compactifications on both sides.

\begin{conjecture}(\cite[Conj~1.3.1]{MMS})\label{PWGeometric}
$\cR_{g,n}(G)$ admits a projective (relative dlt) log CY compactification 
\bEq{pinbar_e}
\Pi_{g,n}\colon  \ov\cR_{g,n}^{\tn{rel}} \lra \C^n
\eEq
such the the dual intersection complex of the boundary divisor exhibits a piecewise-linear (PL) homeomorphism to a sphere.
\end{conjecture}

\begin{remark}
Katzarkov et al.'s initial conjecture posits that the dual intersection complex of the boundary divisor is homotopy equivalent to a sphere.  This weaker version  has been recently proved for $G=\tn{GL}(N,\C)$ in the punctured case; c.f.  \cite{TS}.
\end{remark}

For $G\!=\!\tn{SL}(2,\C)$, the conjecture had been previously known to be true for $(g,n)\!=\!(0,4)$ and $(0,5)$. For $(g,n)\!=\!(0,4)$, it is well-known that $\Pi_4\colon\cR_4 \lra \C^4$ is the Fricke-Klein $4$-parameter (semi-universal) family of affine cubic surfaces. In this case, the natural relative compactification obtained by taking the completion of each fiber in $\C\P^3$ is a family of projective cubic surfaces for which the boundary divisor of each fiber is a triangle of projective lines in $\C\P^2$. The case of $(g,n)=(0,5)$ is substantially more complicated and was treated by Komyo and Simpson in \cite{K,S3}.\\

{\it The main results of this paper, encapsulated in Theorems~\ref{sk-comp_th}-\ref{Sphere_th}, confirm Conjecture~\ref{PWGeometric} for $G=\SL(2,\C)$ in a very strong sense for all $(g,n)$ with $n>0$. We find explicit compactifications for which the boundary divisors are  toric varieties and  provide a systematic approach for computing the associated polytopes/polytope complexes.
In general, the divisor added might not be a normal crossings divisor or $\Q$-Cartier, but it can be made so by performing some (toric) blowups.}

\begin{definition}\label{ideal-trig}
Given integers $g\geq 0$ and $n>0$ with $2g+n\geq 3$, an \textbf{ideal triangulation} of $\Si_{g,n}$ is a triangulation $\De$ of the closed surface $\Si_g=\ov\Si_{g,n}$ such that the set of vertices $V$ is the set of punctured points and the set of edges $E$ is a collection of  disjoint\footnote{i.e., they don't intersect in the interior.} arcs connecting the punctured points; see Section~\ref{Fil-from-int} for more details.
\end{definition}

\begin{definition}
We say a Weil divisor\footnote{In general, any subvariety.} $D\subset X$ is toric if each irreducible component of $D$ is a toric variety and the intersection of any sub-collection of irreducible components is a toric stratum of each component. The \textbf{moment polytope complex} of a toric divisor is the  complex obtained by gluing the moment polytopes of the irreducible components along their identified strata.
\end{definition}

\begin{theorem}(Absolute version)\label{sk-comp_th}
Given integers $g\geq 0$ and $n\!>\!0$, with $2g+n\geq 3$, an ideal triangulation of $\Si_{g,n}$ defines a normal compactification $\ov\cR_{g,n}$ of $\cR_{g,n}$ such that the boundary divisor $\cD_{g,n}\defeq\partial \ov\cR_{g,n}$ is an irreducible toric variety.
\end{theorem}
\begin{theorem}(Relative version)\label{Rel-sk-comp_th}
Given integers $g\geq 0$ and $n\!>\!0$, with $2g+n\geq 3$, an ideal triangulation of $\Si_{g,n}$ defines a normal relative compactification $\Pi_{g,n}\colon \ov\cR^{\tn{rel}}_{g,n}\lra \C^n$ such that, for each $\mft\in \C^n$, the boundary divisor $\cD_{g,n,\mft}\defeq\partial \ov\cR_{g,n,\mft}$ of the fiber $\ov\cR_{g,n,\mft}=\Pi_{g,n}^{-1}(\mft)\subset \ov\cR^{\tn{rel}}_{g,n}$ is toric and  independent of $\mft$. \end{theorem}

The two statements above differ in the following way. The projection map in $\Pi_{g,n}$ in \ref{pin_e} does not extend to the compatification $\ov\cR_{g,n}$ in Theorem~\ref{sk-comp_th}. The relative  compactification $\ov\cR^{\tn{rel}}_{g,n}$ is a blowup of $\ov\cR_{g,n}$ to which $\Pi_{g,n}$ extends. More specifically, let's define $\cD_{g,n,\mft}\subset \cD_{g,n}$ to be the boundary divisor of the closure $\ov\cR_{g,n,\mft}$ of $\cR_{g,n,\mft}$ inside $\ov\cR_{g,n}$. Then, we will show that $\cD_{g,n,\mft}$ is independent of~$\mft$ and we can simply refer to it as $\cD_{g,n}^{\tn{rel}}$. Further, we will also demonstrate that $\cD_{g,n}^{\tn{rel}}\subset \cD_{g,n}$ is a union of toric strata. 
The space $\ov\cR^{\tn{rel}}_{g,n}$ is a blowup of $\ov\cR_{g,n}$ along the base locus $\cD_{g,n}^{\tn{rel}}$ of the family $\{\ov\cR_{g,n,\mft}\}_{\mft\in \C^n}$. The fiber of the extended projection $\Pi_{g,n}$ over a point $\mft$ is the proper transform of~$\mathcal{R}_{g,n,\mft}$.  \\

In Section~\ref{Fil-from-int}, for every ideal triangulation $\De$, we identify a (saturated) finitely-generated monoid
$$
\La=\La\subset \Z_{\geq 0}^{E}\cong \Z_{\geq 0}^{3(2g+n-2)}
$$ 
whose associated real cone $\La_\R=\La \otimes \R\subset \R_{\geq 0}^{E}
$ 
is a non-degenerate rational polyhedral cone. The lattice points in $\Lambda$ are in one-to-one correspondence with curves on $\Si_{g,n}$. The intersection of $\La_\R$ with an affine hyperplane orthogonal to $\vec{1}=(1,\ldots,1)\in \R^{E}$ is the moment polytope $P_{\De}$ of the irreducible toric divisor $\cD_{g,n}$ in Theorem~\ref{sk-comp_th}. Let $P^{\tn{rel}}_{\De}\subset P_{\De}$ denote the union of strata (i.e., sub-polytopes) corresponding to $\cD_{g,n}^{\tn{rel}}\subset \cD_{g,n}$. More precisely, $P^{\tn{rel}}_{\De}$ is the union of those sub-polytopes in $P_{\De}$  that do not intersect the rays corresponding to the peripheral curves. The following confirms the main claim of Conjecture~\ref{PWGeometric}.
\begin{theorem}\label{Sphere_th}
For the relative compactification in\footnote{except for $(g,n)=(0,3)$ where $P^{\tn{rel}}_{\De}$ is empty} Theorem~\ref{Rel-sk-comp_th}, the moment polytope complex $P^{\tn{rel}}_{\De}$ is homeomorphic to a sphere.
\end{theorem}

\begin{remark}
In \cite{W}, Whang proved that the varieties $\cR_{g,n,\mft}$ are log Calabi-Yau. In an earlier version of this paper, we claimed that the result follows from Theorem~\ref{Sphere_th} for the following reason:
Let $K_{g,n}^{\tn{rel}}$ denote the canonical divisor of $\ov\cR_{g,n,\mft}$. By adjunction, the toric property in Theorem~\ref{Rel-sk-comp_th}, and Theorem~\ref{Sphere_th}, the restriction $\cO(K_{g,n}^{\tn{rel}}+{\cD}_{g,n}^{\tn{rel}})|_{{\cD}_{g,n}^{\tn{rel}}}$ is trivial. Since $\cD_{g,n}^{\tn{rel}}$ is ample, Grothendieck-Lefschetz theorem implies that  $K_{g,n}^{\tn{rel}}+{\cD}_{g,n}^{\tn{rel}}$ is trivial. However, as it was pointed out by the referee, using adjunction requires checking three conditions, namely (1) $\ov\cR_{g,n,\mft}$ is Cohen- Macaulay, (2) ${\cD}_{g,n}^{\tn{rel}}$ is $\Q$-Cartier, and (3) ${\cD}_{g,n}^{\tn{rel}}$ is ample; see  \cite[Proposition~5.73]{MK} and \cite[Section~4.1]{Kollar}. As some of the examples worked out in this paper show, condition (2) is not always satisfied. Thus, in order to fix the problem, one would either need to perform a blowup (which compromises ampleness) or work with the notion of “different” as studied in \cite[Section~4.1]{Kollar}. Both of these as well as establishing (1) require more work than we had anticipated. \qed
\end{remark}

\textbf{\large Idea of the proof of the main result.} By a Theorem of Bullock-Przytycki-Sikora \cite{Bu,PS} and Charles-March\'{e} \cite{CM}, the ring of regular functions $\C[\cR_{g,n}]$ of $\cR_{g,n}$ is canonically isomorphic to the specialization at $-1$ of the skein algebra of $\Si_{g,n}$.  Furthermore, an ideal triangulation of $\Si_{g,n}$ allows us to identify the generators of  $\C[\cR_{g,n}]$, as a vector space, with the lattice points of $\Lambda$ and define a suitable filtration on that. Theorems~\ref{sk-comp_th}-\ref{Rel-sk-comp_th} then follow from algebraic results of Frohman et al. \cite{AF,FK} on the correspondence between the product structure on the skein algebra and the additive structure on $\Lambda$. Proof of Theorem~\ref{Sphere_th} starts with a particular kind of triangulation (for which the sphere property can be proved directly) and shows that the sphere property is preserved through ``mutations" of~$\De$. In general, our approach is closely related to Manon's approach in \cite{Man}. In \cite{Man}, Manon uses skein relations to prove, for instance, that the ring of regular functions of $\cR_g$ is  presented by homogeneous skein relations.  Then, \cite[Thm.~1.6]{Man} provides some qualitative description of the resulting boundary divisor. 

\subsection{Compactifications obtained from $\ov\SL(2,\C)$}\label{part2}

A different approach to constructing a compactification of $\cR_{g,n}$ or $\cR_{g,n,\mft}$ is by replacing $G=\SL(2,\C)$ with a projective variety $\ov{G}=\ov\SL(2,\C)$ in (\ref{GIT-cR_e}). This is the method considered in \cite{K} (for $g=0$) and \cite{BLR}. In the second half of the paper, we exploit and streamline Komyo's compactification method in \cite{K}, and address a question in \cite{S} by constructing uniform relative compactifications over the base space $\C^n$ that sit in a flat fiber bundle with explicit monodromy maps.\\

More explicitly, in \cite{K}, a natural compactification $\ov\SL(2,\C)$ of $\tn{SL}(2,\C)$ (a quadratic threefold) is used to construct a compactification of $\cR_n$ or a single fiber of that. In Section~\ref{Kcomp-sec}, we revisit his construction and answer the following question raised in \cite[Sec.~1.3]{S} (attributed to Deligne).

\begin{question} What can be said as $\mft$ varies? e.g., can the compactification be done uniformly in $\mf{t}$ to obtain a relative  compactification (\ref{pinbar_e})
of the family (\ref{pin_e}) and what can be said about the monodromy?
\end{question}

To address this question, we will use GIT quotients 
$$
X_{\mf{a}}\defeq (\P^1)^{m} \sslash_{\cL_{\mf{a}}} \SL(2,\C)
$$ 
of $(\P^1)^{m}$ with respect to the diagonal M\"obius action, linearized by an ample line bundle $\cL_{\mf{a}}$, labelled by a tuple of positive integers $\mf{a}=\Z_+^m$. In Theorem~\ref{main_thm} and~\ref{real-main_thm}, we will need the case where $m=2(n-1)$ and $\mf{a}$ is \textbf{symmetric} in the sense that 
\bEq{atob_e}
\mf{a}=(b_1,b_1,\ldots,b_{n-1},b_{n-1})
\eEq
for some $(b_1,\ldots,b_{n-1})\in \Z_+^{n-1}$; see Definition~\ref{Syma_dfn}. For instance, there is a choice of $\mf{a}$ for which $X_{\mf{a}}$ is an explicit toric variety; see Lemma~\ref{toric_lmm}. \\

Every such $X_{\mf{a}}$ admits mutually commuting (holomorphic) involutions 
$$
\si_i\colon X_{\mf{a}}\lra X_{\mf{a}}, \qquad i=1,\ldots,n-1.
$$
Also, every $X_{\mf{a}}$ comes with a positive union $D_{\mf{a}}\!=\!\bigcup_{i=1}^{n-1} D_{\mf{a},i}$ of divisors  such that $D_{\mf{a},i}\subset \tn{Fix}(\si_i)$ and each of the involutions above preserves every $D_{\mf{a},i}$.\\

Let $\cA\cong \C^*$ denote the affine\footnote{We can consider the completion of this curve to a smooth conic curve in $\P^2$ by adding two points at the infinity. The family $\mf{X}^\circ_n$ constructed below and Theorem~\ref{main_thm} extend to this compactification. This is only useful if we study completions of the family (\ref{pin_e}) in the base direction.} curve $(t^2-s^2=4)\subset \C^2$, and 
\bEq{cB_e}
\cB= \cA_1\times \cdots\times \cA_{n-1}\times \C,
\eEq
where each $\cA_i$ is a copy of $\cA$ with coordinates $(t_i,s_i)$.
We define actions of $\Z_2$ on $\cB$, denoted similarly by $\si_i$, in the following way. For $i\!=\!1,\ldots,n-1$, let $\si_i$ act on $\cB$ by $(t,s)\!\lra\! (t,-s)$ on $\cA_i$ and fixing the other components. Let $\Gamma=\ll \si_1,\ldots,\si_{n-1}\rr\cong \Z_2^{n-1}$ denote the group generated by the involutions $\si_i$. Note that the projection
$$
\pi_{\cB}\colon \cB \lra \C^n,\quad \big((t_1,s_1),\ldots,(t_{n-1},s_{n-1}),t_n\big)\lra \mft=(t_1,\ldots,t_n)
$$ 
descends to an isomorphism $\cB/\Gamma \stackrel{\cong}{\lra} \C^n$. Fix a choice of symmetric $\mf{a}\in \Z_+^{2(n-1)}$ and let 
\bEq{Bcirc}
\aligned
&\cB^\circ=\pi_{\cB}^{-1}\big( \C^n_\circ\big), \quad \C^n_\circ\defeq (\C-\{\pm 2\})^{n-1}\times \C,\\
&\mf{X}^\circ_n\defeq (X_{\mf{a}}\times \cB^\circ)/\Gamma \quad \tn{and}\quad \mf{D}^\circ_{n}=(D_{\mf{a}}\times \cB^\circ)/\Gamma
\endaligned
\eEq
where $\Gamma$ acts on $X_{\mf{a}}$, $\cB$, and the divisor $D_{\mf{a}}$ in the way explained above. The projection map
$$
\mf{X}_n^\circ\lra \C^n_\circ, \qquad \Big[x \times \big((t_1,s_1),\ldots,(t_{n-1},s_{n-1}),t_n\big)\Big]\lra (t_1,\ldots,t_n)
$$
is well-defined and it is a flat $X_{\mf{a}}$-fiber bundle such that the monodromy map around both of the (omitted) complex hyperplanes $ (t_i=\pm 2)\subset \C^n
$ is  $\si_i$. Clearly, $(\mf{X}_n^\circ, \mf{D}^\circ_n)$ depend on $\mf{a}$; however, we exclude the subscript $\mf{a}$ to avoid a cluttered notation.

\begin{theorem}\label{main_thm}
For $n> 3$ and every choice of $(X_{\mf{a}},D_{\mf{a}})$ as above, there exists a relative compactification (\ref{pinbar_e}), also depending on the choice of $(b_1,\ldots,b_{n-1})$ in (\ref{atob_e}), such that the restriction 
$$
 \ov\cR^{\tn{rel},\circ}_n\defeq \Pi_n^{-1}(\C^n_\circ) \lra \C^n_\circ
$$  
can be identified with a divisor in $\mf{X}^\circ_n$ linearly equivalent to $\mf{D}^\circ_n$, identifying the given projection maps to $\C^n_\circ$.  \vskip.1in

Phrasing it differently, there exists a $\Gamma$-equivariant morphism 
$
\Psi\colon \cB^\circ \lra |D_{\mf{a}}|
$
into the linear system of $D_{\mf{a}}\subset X_{\mf{a}}$ such that the divisor 
$$
\Psi\big((t_1,s_1),\ldots,(t_{n-1},s_{n-1}),t_n\big)\subset X_{\mf{a}}
$$ is a compactification $\ov\cR_{n,\mft}$ of $\cR_{n,\mft}$ with boundary
$$
\cD_{n,\mf{t}}=\ov\cR_{n,\mft}-\cR_{n,\mft}=\ov\cR_{n,\mft} \cap D_{\mf{a}}.
$$
\end{theorem}

Note that every such $X_{\mf{a}}$ is rational. It is an interesting question whether $\ov\cR^{\tn{rel},\circ}_n \lra \C^n_\circ$ is a semi-universal family of hypersurfaces in $X_{\mf{a}}$ in the class of~$D_{\mf{a}}$ (Similarly to the Fricke-Klein semi-universal family of affine cubic surfaces).

\begin{remark}\label{mfM_rmk}
As we will explain in Section~\ref{RC}, the entire relative compactification $\ov\cR^{\tn{rel}}_n$ also sits as a divisor in some easy to describe space $\mf{M}_n$. The completed quotient space $\mf{X}_n\defeq (X_{\mf{a}}\times \cB)/\Gamma$ and $\mf{M}_n$ are birational and related by some explicit blowups and blowdowns. However, when some $t_i=\pm 2$, the fibers of the closure of $\ov\cR^{\tn{rel},\circ}_n$ in $\mf{X}_n$ include $\P^2$-components corresponding entirely to the identity matrix. This causes us to miss actual matrices corresponding to $\pm 2$, rendering the resulting compactification meaningless.
\end{remark}

The main step of the proof of Theorem~\ref{main_thm} is Proposition~\ref{Main_Prop} that provides a uniform parameterization of the slices of $\SL(2,\C)$ and exchanges the conjugation action of $\SL(2,\C)$ on each slice with the diagonal action of $\SL(2,\C)$ on $\P^1\times \P^1$ via m\"obius transformations.\\

Finally, in Section~\ref{real}, we prove the following theorem for the character varieties $\cR_{n}(\SL(2,\R))^{\tn{elliptic}}$ and $\cR_{n}(\SU(2))$ that sit in the ``real" locus $\cR_n^\R$ of $\cR_n$. Here, $\cR_{n}(\SL(2,\R))^{\tn{elliptic}}$, defined in (\ref{elliptic-ChVar}), is the subspace corresponding to matrices $A$ with $\abs{\tn{tr}(A)}<2$, and $\cR_{n}(\SU(2))$ is already compact.

\begin{theorem}\label{real-main_thm}
For $(X_{\mf{a}},D_{\mf{a}})$ as in Theorem~\ref{main_thm}, there are anti-holomorphic involutions 
$$
\eta,\tau\colon X_{\mf{a}}\lra X_{\mf{a}}
$$
such that, each component $D_{\mf{a},i}$ of $D_{\mf{a}}$ is preserved by $\tau$ and $\eta$ and $\tn{Fix}(\eta)\subset X_{\mf{a}}-D_{\mf{a}}$. Furthermore, for every $\mft\in (-2,2)^n$, with respect to the complex codimension one inclusions $(\ov\cR_{n,\mft},\cD_{n,\mft})\subset (X_{\mf{a}},D_{\mf{a}})$ given by $\Psi$ in Theorem~\ref{main_thm}, $\ov\cR_{n,\mft}$ is preserved by $\tau$ and $\eta$ and 
\bIt
\item $\tn{Fix}(\eta)\cap \ov\cR_{n,\mft}=\tn{Fix}(\eta)\cap \cR_{n,\mft}=\cR_{n,\mft}(\SU(2))\subset\cR^\R_{n,\mft},$  
\item $\tn{Fix}(\tau)\cap \ov\cR_{n,\mft} =\ov\cR_{n,\mft}(\SL(2,\R))^{\tn{elliptic}}\subset\cR^\R_{n,\mft}$ is a compactification of $\cR_{n,\mft}(\SL(2,\R))^{\tn{elliptic}}$ obtained by attaching the unbounded components of $\cR_{n,\mft}(\SL(2,\R))^{\tn{elliptic}}$ along the limiting real hypersurfaces 
$$
\ov\cR_{n,\mft}(\SL(2,\R))^{\tn{elliptic}}\cap D_{\mf{a}}.
$$
\eIt
\end{theorem}
\vskip.2in

\textbf{Acknowledgment}. We wish to thank Michael Thaddeus, J\'anos Koll\'ar, Jason Starr, Travis Mandel, Greg Muller, and Pinaki Mondal for answering our questions. We also wish to thank Travis Mandel and Mirko Mauri for their comments on the first version of the paper. We are grateful to the referee for a careful reading of the paper and pointing out some incorrect statements.

\section{Compactifications obtained from triangulations}\label{SkeinComp_sec}

In Section~\ref{filtration}, we review the results of \cite{Mon} on using filtrations defined by degree-like functions to compactify affine varieties.
In Section~\ref{skein-algebra_sec}, we review the definition and properties of the skein algebra of a surface.
In Section~\ref{Fil-from-int}, following \cite{AF}, we introduce a nice class of degree-like functions on skein algebra and study the resulting compactifications. 

\subsection{Basic facts about compactification via filtration}\label{filtration}
It is a well-understood fact in algebraic geometry that there is a correspondence between filtrations (satisfying some nice conditions) on the ring of regular functions of an affine variety $Z$ and (projective) compactifications or partial compactifications $\ov{Z}$ of $Z$. In this section, following the presentation of \cite{Mon}, we provide a quick review of the general definitions and statements that we will use in Section~\ref{Fil-from-int}. \\

Let $Z$ be an affine variety (over $\C$). A ($\Z$-)\textbf{filtration} $\mc{F}$ on the ring  of regular functions $\C[Z]$ of $Z$ is an increasing sequence of subspaces 
$$
\cdots F_{-1} \subset F_0\subset F_1 \subset \cdots \subset \C[Z]=\bigcup_{d\in \Z} F_d
$$
such that $F_d F_e \subset F_{d+e}$ for all $d,e\geq 0$ and $1\in F_0/F_{-1}$. 
To each filtration $\cF$, we can associate two graded algebras: 
\bEq{GAlg_e}
\C[Z][u]_{\cF}\defeq \bigoplus_{d\geq 0} F_d u^d \quad \tn{and}\quad \C[Z]_{\cF}\defeq \bigoplus_{d\geq 0} F_d/F_{d-1}.
\eEq
where $u$ is an indeterminate over $\C[Z]$. Filtrations correspond to {\it degree-like} functions on $\C[Z]$ in the following sense.

\bDf{dl-function}
A  \textbf{degree-like} function on $\C[Z]$ is an integer-valued function 
$\de\colon \C[Z]\!\to\! \Z\!\cup\!  -\infty$ satisfying 
\bEn
\item $\de(\C^*)=0$ and $\de(0)=-\infty$;
\item $\de(f_1+f_2)\leq \max\{\de(f_1),\de(f_2)\}$ for all $f_1,f_2\in \C[Z]$, with $<$ in the preceding equation implying $\de(f_1)=\de(f_2)$;
\item $\de(f_1f_2)\leq \de(f_1)+\de(f_2)$ for all $f_1,f_2\in \C[Z]$.
\eEn
We say $\de$ is a \textbf{semi-degree} if  $\de(f_1f_2)= \de(f_1)+\de(f_2)$ for all $f_1,f_2\in \C[Z]$. We say $\de$ is a \textbf{sub-degree} if $\de=\max\{\de_1,\ldots,\de_m\}$ such that $\de_1,\ldots,\de_m$ are semi-degrees. 
\eDf

The filtration induced by a degree-like function $\de$ is simply
$$
F_d\defeq \{f\in \C[Z]\colon \de(f)\leq d\} \qquad \forall~d \in \Z.
$$
Conversely, every filtration defines a degree-like function $\de(f)\defeq \inf \{d\in \Z\colon f\in F_d\}.$
When we are working with degree-like functions, we may alternatively denote $\C[Z][u]_{\cF}$ and $\C[Z]_{\cF}$ by $\C[Z][u]_{\de}$ and $\C[Z]_{\de}$, respectively. The motivation behind the notions of semi and sub-degrees is that the associated filtrations define compactifications that have a better behavior at infinity. The precise statements are explained below. By \cite[Thm.~1.3]{Mon}, every sub-degree $\de$ has a unique minimal presentation $\de=\max\{\de_1,\ldots,\de_m\}$ in which $\de_i$s are semi-degrees.

\begin{definition}
We call a filtration $\cF$  (equivalently, the associated degree-like function $\de$) \textbf{finitely-generated} if the graded algebra $\C[Z][u]_{\cF}$ is finitely generated. We call $\cF$ \textbf{projective} if it is finitely generated and $F_0=\C$.
\end{definition}

In general, by \cite[Prop.~2.8]{Mon}, for every filtration $\cF$,  
$$
\ov{Z}^\cF \defeq \tn{Proj}(\C[Z][u]_{\cF})
$$ 
is a reduced scheme containing $Z$ as a Zariski dense open subset. The complement $\partial  \ov{Z}^\cF=\ov{Z}^\cF -Z$ of $Z$ in $\ov{Z}^\cF$ is the zero set of the ideal $\ll u\rr$ generated by $u$ and is isomorphic (as a scheme) to $\tn{Proj}(\C[Z]_{\cF})$.  Extra assumptions on $\cF$ result in a better behavior of $\ov{Z}^\cF$ along $\partial  \ov{Z}^\cF$.\\

If $\cF$ is finitely generated, then $\ov{Z}^\cF$ is a closed sub-variety of a {\it generalized weighted projective space}\footnote{Following \cite[Definition~2.7]{Mon}, a generalized weighted projective space is a weighted projective space possibly with negative weights.} and $\partial\ov{Z}^\cF$ is the support of an effective ample divisor. If in addition $\cF$ is a projective, then $\ov{Z}^\cF$ is a projective variety. We will need the first condition for constructing a relative compactification as in (\ref{pinbar_e}) and the second one for constructing a normal (overall) compactification of $\cR_n$.  \\

Filtrations defined by semi and sub-degrees define compactifications that have a better behavior at infinity.
{\it If $\de=\max\{\de_1,\ldots,\de_m\}$ is a sub-degree, then $\ov{Z}^\de$ is non-singular in codimension one at infinity and the number of irreducible components of $\partial\ov{Z}^\de$ is $m-1$ or $m$ depending on whether one of $\de_i$ is the zero degree-like function or not. In fact, if $Z$ is normal, then so is $\ov{Z}^\de$.} Note that if $\de$ is a semi-degree, then $\partial\ov{Z}^\de$ is irreducible.

\begin{example}
Every $n$-dimensional toric variety is the compactification $\ov{\T}^\de$ of the $n$-dimensional complex torus $\T=(\C^*)^n$ with respect to a sub-degree on $\C[\T]=\C[x_1^\pm, \ldots, x_n^\pm]$ determined by the boundary equations of the associated polytope; c.f. \cite[Sec.~1.2]{Mon}. 
\end{example}

\bRm{sub-condition}
If $\C[Z][u]_{\de}$ is Noetherian, then $\de$ is a sub-degree if and only if $\de(f^k)=k\de(f)$ for all $f\in \C[Z]$ and $k\geq 0$. 
If $\de$ is projective and $\partial\ov{Z}^\de$ is the support of an effective ample divisor but $\de$ is not a sub-degree, then by \cite[Thm~5.10]{Mon}, the degree-like function $\wt\de$ defined by an integer multiple of $\lim_{k\lra \infty} \frac{\de(f^k)}{k\de(f)}$
is a sub-degree such that $\ov{Z}^{\wt{\de}}$ is the normalization at infinity of $\ov{Z}^{\de}$.
\eRm 

Remark~\ref{sub-condition} will be useful in situations like the following scenario to prove that the compactification is normal without directly finding the semi-degrees that define a sub-degree filtration.  Suppose $Y\subset Z$ is an affine sub-variety and $\de$ is a projective semi-degree on $\C[Z]$ with the associated filtration $\cF$. The image of $\cF$ under the surjective projection map $\C[Z]\lra \C[Y]$ gives us filtration $\cF^Y$ on $\C[Y]$ that corresponds to the closure $\ov{Y}$ of $Y$ in the projective variety $\ov{X}^\cF$ (Note that since $\de$ is a semi-degree, $\partial\ov{X}^\cF$ is an irreducible divisor). The induced filtration $\cF^Y$ does not necessarily come from a sub-degree on $\C[Y]$ as $\ov{Y}$ can be singular in co-dimension one. For instance, $\ov{Y}$ can be a nodal curve in $\P^2=\ov{\C^2}$ that has a node on the boundary divisor $\partial \;\ov{\C^2}=\P^1$. Remark~\ref{sub-condition} provides a way of checking whether $\ov{Y}$ is normal, or if it is not, to find a sub-degree on $\C[Y]$ that defines the normalization of $\ov{Y}$. 

\subsection{Skein algebra}\label{skein-algebra_sec}

For $g,n\geq 0$, the \textbf{skein algebra} $\tn{Sk}_{A}(\Si_{g,n})$ of $\Si_{g,n}$ is the $\Z[A^\pm]$-algebra generated by isotopy classes of framed links in $\Si_{g,n}\times \R$ modulo the relations 
$$
\aligned
    \left<L\cup\KPA\right> &= (-A^{2}-A^{-2})\langle L\rangle
        &&\text{disjoint union with trivial loop}  \\
    \left<\KPB\right> &= A\left<\KPC\right> + A^{-1}\left<\KPD\right>
        &&\text{skein relation,} 
        \endaligned
$$ 
with the product structure defined by stacking two links on the top of each other. A \textbf{curve} on $\Si_{g,n}$ is an embedded circle. A \textbf{multi-curve}  on $\Si_{g,n}$ is the union of finitely many disjoint curves (with vertical framing in $\Si_{g,n}\times \R$) such that none of them bounds a disc (is trivial). By \cite[Thm.~IX.7.1]{P1}, isotopy classes of multi-curves (together with the empty multi-curve  $\cong$ 1) form a basis of $\tn{Sk}_{A}(\Si_{g,n})$ as $\Z[A^\pm]$-module. On the other hand, by a Theorem of Bullock-Przytycki-Sikora \cite{Bu,PS} and Charles-March\'{e} \cite{CM}, the specialization $\tn{Sk}_{g,n}\defeq \tn{Sk}_{A=-1}(\Si_{g,n})$ of the skein algebra at $A = -1$ is canonically isomorphic, via the (negative of) trace function along loops, to the ring of regular functions $\C[\cR_{g,n}]$ of $\cR_{g,n}$. The curves $a_1,\ldots,a_n$ around the punctures, called the \textbf{peripheral} curves, are in the center of $\tn{Sk}_{A}(\Si_{g,n})$. Therefore, we can view $\tn{Sk}_{A}(\Si_{g,n})$ as an $\Z[A^{\pm}][a_1,\ldots, a_n]$-module. This inclusion is a quantization of the projection $\Pi_{g,n}$ in (\ref{pin_e}). 

 \begin{remark}\label{spin_rmk}
 The classical limit $\tn{Sk}_{A=1}(\Si_{g,n})$ is isomorphic to $\tn{Sk}_{A=-1}(\Si_{g,n})$, just not canonically. The isomorphism depends on a choice of a spin structure; see \cite{B}. The former behaves better when generalizing the skein algebra to include arcs and producing cluster algebras \cite{FST}.  We will use the non-canonical isomorphism $\Phi\colon \tn{Sk}_{A=-1}(\Si_{g,n})\lra \tn{Sk}_{A=1}(\Si_{g,n})$ in the proof of Theorem~\ref{sk-comp_th}.
 \end{remark}
 
 In the classical limit, isotopy classes of multi-curves (together with the empty curve $\cong 1$) form a vector space  basis for $\tn{Sk}_{g,n}$ over $\C$. The product structure (denoted by $\ast$) is by taking union of two curves (in transverse position), and using the skein relations as reduction rules  
\bEq{skein-relation}
    \left<L\cup\KPA\right> =-2\langle L\rangle,\qquad 
    \left<\KPE\right> + \left<\KPC\right> + \left<\KPD\right> =0
    \eEq
    to expand the product as a linear combination of multicurves. 
 As before, if $g=0$, we will simplify the notation $\tn{Sk}_{0,n}$ to $\tn{Sk}_{n}$. \\
 
 For $\mft=(t_1,\ldots,t_n)\in \C^n$, the inclusion $\cR_{g,n,\mft} \subset \cR_{g,n}$ corresponds to a projection map from $\tn{Sk}_{g,n}$ to the ring of regular functions $\tn{Sk}_{g,n,\mft}\defeq \C[\cR_{g,n,\mft}]$ of $\cR_{g,n,\mft}$. As a vector space, there is a right-inverse for this projection and $\tn{Sk}_{g,n,\mft}$ can be identified with the subspace $\tn{Sk}^{\tn{rel}}_{g,n}\subset \tn{Sk}_{g,n}$ generated by multi-curves that do not include any peripheral curve. Under this linear identification, however, the product structure $\ast_{\mft}$ induced on $\tn{Sk}^{\tn{rel}}_{g,n}$ will (be different from $\ast$ and) depends on $\mft$ - the induced product structure $\ast_{\mft}$ is obtained by replacing the $i$-th peripheral curve $a_i$ with  the constant $-t_i$, whenever such curves arise in the resolution of crossing points.  

\subsection{Filtration by intersection number (Proof of Theorem~\ref{sk-comp_th})}\label{Fil-from-int}

In this section, following the works of Frohman et al. \cite{AF,FK}, we introduce a class of filtrations on $\tn{Sk}_{g,n}=\C[\cR_{g,n}]$ that are defined in terms of the geometric intersection number of multi-curves with certain collection of arcs, and show that the resulting compactifications have interesting properties.  \\

An \textbf{arc} on the closure $\Si_g=\ov{\Si}_{g,n}$ is the image of an embedding $[0,1]\lra \Si$ such that $0$ and $1$ are mapped to the puncture points and $(0,1)$ is mapped to the interior $\Si_{g,n}$ of $\Si$. Two arcs are called {\bf disjoint} if they don't intersect in the interior.  \\

Given two curves $c_1$ and $c_2$, the {\bf geometric intersection number} of $c_1$ and $c_2$ is the quantity
$$
i(c_1,c_2)= \substack{\min \\ c'_1\sim c_1,~ c'_2\sim c_2 }~ \# (c'_1\cap c'_2)\in \Z_{\geq 0},
$$ 
where the minimum runs over the pairs $(c'_1,c'_2)$ of transversally intersecting curves in the compactly-supported ambient isotopy classes of $(c_1,c_2)$. For instance, $i(c,c)=0$ and $i(c,a)=0$ for every curve $c$ and any peripheral curve $a$. The intersection number between a curve and an arc is defined similarly, as they can only intersect at the interior points of the arc. For a more general variation of this definition and its application in constructing cluster algebras see \cite[Dfn.~8.4]{FST}.\\
 
For a pair $(c_1,c_2)$ of transverse curves/arcs, a {\bf bigon} is an embedded disc $B\subset \Si$ such the interior is disjoint from $c_1\cup c_2$ and the boundary $\partial B$ decomposes as a union of closed arcs $a \cup b$, where $a\subset c_1$ and $b\subset c_2$.
By Bigon Criterion (c.f. \cite{FLP}), a transverse pair of curves/arcs $(c_1,c_2)$ realizes the geometric intersection number $i(c_1,c_2)$ if and only if there are no bigons. If there is a bigon, there is always an innermost bigon, whose interior is disjoint from $c_1\cup c_2$ (so one can remove it and continue inductively to remove all the bigons and achieve the minimum intersection number).

\begin{remark}
Given two finite sets of curves/path $S_1=\{a_1, \cdots, a_k\} $ and $S_2=\{b_1, \cdots, b_\ell\}$, there exist (isotopically) equivalent collections $S'_1=\{a'_1, \cdots , a'_k \}$ and $S'_2=\{b'_1, \cdots , b'_\ell\}$ such that $i(a_r,b_s)=\# a_s'\cap b'_s$ for all $r,s$. Therefore, the definition of geometric intersection naturally extends to finite collections of curves/arcs satisfying the identity $i(S_1,S_2)= \sum_{r,s} i(a_r,b_s)$. Furthermore, intersection number is symmetric; i.e., $i(S_1,S_2)=i(S_2,S_1)$.
\end{remark}

\bDf{dl-intersection}
Give a finite collection of curves and arcs $S$, define $
\de_S\colon \tn{Sk}_{g,n}\lra \Z_{\geq 0} \cup -\infty
$ by sending any  finite non-zero linear combination $f= \sum_{s} \la_s \mf{m}_s \in \tn{Sk}_{g,n}$ of multi-curves $\mf{m}_s$  to 
\bEq{max-def}
\de_S(f)\defeq \max_{s\colon \la_s\neq 0}\{ i(\mf{m}_s,S)\}.
\eEq
We define $\de_S$ on $\C^*$ (non-zero multiples of a contractible curve) to be $0$ and $\de_S(0)= -\infty$.  \eDf

\begin{lemma}
For any collection $S$ of curves and arcs, $\de_S$ is a degree-like function.
\end{lemma}

\bPf
The first property in Definition~\ref{dl-function} holds by definition and the second one follows from the max function on the righthand side of (\ref{max-def}). If $\mf{m}$ and $\mf{m}'$ are transverse multi-curves, by definition, we have $\mf{m}\ast \mf{m}'=\sum_{s} \la_s \mf{m}_s$ where each multi-curve $\mf{m}_s$ is obtained by forming a simultaneous resolution $\wt{\mf{m}}_s$ of the crossing points (known as a \textbf{state}) and then substituting the contractible curves with $-2$ to get $\mf{m}_s$. We can assume that these resolutions are happening in sufficiently small balls so that none of the states $\wt{\mf{m}}_s$ includes any of the intersection points with $S$. Therefore, 
$$
\de_S(\mf{m}_s)\leq \de_{S}(\mf{m})+\de_S(\mf{m}')\qquad \forall~s~\tn{s.t}~\la_s\neq 0,
$$ 
with the equality happening if and only if there are no bigons in $\wt{\mf{m}}_s\cap S$. For arbitrary elements 
$$
f_1=\sum_s \la_{1,s}\, \mf{m}_{s} \quad \tn{and}\quad f_2=\sum_t \la_{2,t}\, \mf{m}_{t},
$$
the product $f_1\ast f_2$ decomposes as
$$
f_1\ast f_2= \sum_{s,t} \la_{1,s}\la_{2,t} \,\mf{m}_{s}\ast \mf{m}_{t}=\sum_{s,t,r} \la_{1,s}\la_{2,t} \la_{s,t,r}\,\mf{m}_{s,t,r},
$$
where it is possible to have $\mf{m}_{s,t,r}=\mf{m}_{s',t',r}$ for some $(s',t')\neq(s,t)$. This could result in cancellations on the righthand side, if the coefficients of certain identical multi-curves in the right summation add up to zero. Nevertheless,
\bEq{product-expansion}
\aligned
&\de_S(f_1\ast f_2)\leq \tn{max}_{s,t,r}~i(\mf{m}_{s,t,r},S)=i(\mf{m}_{s_0,t_0,r_0},S)\qquad \tn{for some}~(s_0,t_0,r_0),\\
\tn{and}\quad&i(\mf{m}_{s_0,t_0,r_0},S) \leq ~i(\mf{m}_{s_0},S)+i(\mf{m}_{t_0},S)\leq \de_S(f_1)+\de_S(f_2).
\endaligned
\eEq
Therefore, we conclude that the third property in Definition~\ref{dl-function} holds as well.
\ePf

\begin{remark}
The sum of two degree-like functions $\de_1$ and $\de_2$ is \underline{not} necessarily a degree-like function, because the second property in Definition~\ref{dl-function} may fail to hold for  $\de_1+\de_2$. In particular, for $S=S_1\cup S_2$, there is no obvious relation between $\de_S$ and $\de_{S_1}+\de_{S_2}$.\end{remark}

Recall from Definition~\ref{ideal-trig} that an ideal triangulation of $\Si_{g,n}$  is a maximal collection $\De$ of non-isotopic disjoint arcs between the $n$ puncture points of $\Si_g=\ov{\Si}_{g,n}$. An ideal triangulation may include so-called \textbf{folded} triangles. A folded triangle is a triangle in which two sides are supported on the same arc; see~Figure~\ref{folded-T_fig}. Some ideal triangulations of $\Si_{1,1}$ and $\Si_{0,4}$ are illustrated in Figures~\ref{gn11-fig} and  \ref{n4EX}. Every triangle has three \textbf{corners}, determined by a choice of two of the three sides (the folded arc of a folded triangle contributes two sides to the unfolding of the triangle). Every ideal triangulation of $\Si_{g,n}$ is made of $n$ vertices, $3(n+2g-2)$ arcs, $2(n+2g-2)$ triangles, and thus $6(n+2g-2)$ corners. We denote the corresponding sets by $V(\De)$, $E(\De)$, $T(\De)$, and $\Theta(\De)$, respectively. When $\De$ is fixed in any discussion, we will simply write $V$, $E$, $T$, and $\Theta$, respectively. We follow the same convention for the other objects associated to $\De$.\\

\begin{figure}
\begin{pspicture}(-15,-.8)(0,.8)
\psset{unit=1cm}

\pscircle*(-7,1){0.1}
\pscircle*(-7,0){0.1}\rput(-6.8,0.5){\small $c$}
\psline(-7,0)(-7,1)
\pscircle(-7,0){1}

\end{pspicture}
\caption{A triangle folded on the arc $c$}
\label{folded-T_fig}
\end{figure}

Let 
$$
[N]\defeq\{1,\ldots,N\}\qquad \forall~N\geq 0\quad \tn{with}~~[0]=\eset.
$$
\vskip.1in

\textbf{Proof of Theorem~\ref{sk-comp_th}.}
As we explain below, Theorem~\ref{sk-comp_th} is essentially a consequence of \cite[Prop.~3.2, Thm~3.4]{AF}. {\it We show that the degree-like function $\de_\De\colon \tn{Sk}_{g,n}\lra \Z_{\geq 0} \cup -\infty$ associated to any ideal triangulation is a projective semi-degree and identify the graded algebra $\C[\cR_{g,n}]_{\de_\De}=(\tn{Sk}_{g,n})_{\de_\De}$ with the semigroup algebra of a specific saturated affine semi-group.}\\ 

The argument below mostly highlights and explains the main arguments in the proofs of \cite[Prop.~3.2, Thm~3.4]{AF}  and adds some extra details needed in the current context. We skip some of the fine details and refer to \cite{AF}. \\

Let $E=\{e_1,\ldots,e_N\}$, $N=3(n+2g-2)$. Labeling the edges as $e_1,\ldots,e_N$ fixes an order which we will use later in the proof. Also, let $\Theta=\{\theta_1,\ldots,\theta_{2N}\}$ denote the set of corners (ordered in an arbitrary way, if needed).

\begin{definition}
An {\bf admissible coloring} of $\De$ is a non-negative integer-vector $v=(v_{e})_{e\in E}\in \Z_{\geq 0}^{E}$ such that for every triangle $\{e,e',e''\}\in T$ we have
\bEq{admissible_e}
2\mid v_{e}+v_{e'}+v_{e''}\qquad \tn{and}\qquad v_{e}\leq v_{e'}+v_{e''},~v_{e'}\leq v_{e}+v_{e''},~v_{e''}\leq v_{e}+v_{e'}.
\eEq
Define $\La=\Lambda(\De) \subset \Z_{\geq 0}^{E}$ to be the subset of admissible colorings.
\end{definition}

\begin{remark}
Note that if $\{e,e',e''\}$ is a folded triangle, then two of the edges are the same, say $e'=e$, and  (\ref{admissible_e}) reads 
$$
2\!\mid\! v_{e''}\qquad \tn{and}\qquad 2v_{e}\geq  v_{e''}.
$$
\end{remark}
\begin{remark}
Whenever needed, an ordering on $E$ gives an isomorphism $\Z_{\geq 0}^{E}\cong \Z_{\geq 0}^N$ and we will write $v=(v_i)_{i\in [N]}$ instead of $v=(v_{e_i})_{i\in [N]}$.
\end{remark}

It is clear from (\ref{admissible_e}) that $(\Lambda,+)$ is monoid with the identity element $\vec{0}=(0_e)_{e\in  E}$. However, due to the first condition in (\ref{admissible_e}), it is not clear whether $\Lambda$ is isomorphic to the lattice points of a finitely generated integral cone\footnote{A strongly convex rational polyhedral cone.} in some $\R_{\geq 0}^k$ or not. The following change of coordinate linear map in the proof of \cite[Prop.~3.2]{AF} provides such an embedding
$$
\iota \colon \Lambda\lra \Z_{\geq 0}^{\Theta}\cong \Z_{\geq 0}^{2N}.
$$
If $\{e,e',e''\}$ is a triangle in $\De$ and $\theta$ is the corner between $e$ and $e'$, given $v=(v_e)_{e\in E} \in \Lambda$, define the $\theta$-component $u_\theta$ of the vector $u=\iota(v)$  to be 
\bEq{Corner_e}
u_\theta=\frac{1}{2}(v_e+v_{e'}-v_{e''}).
\eEq
By (\ref{admissible_e}), $u=(u_\theta)_{\theta\in \Theta}\in \Z_{\geq 0}^{\Theta}$. It is easy to see that $\iota$ is injective. Furthermore, $\iota(\Lambda)$ is the set of non-negative vectors $u$ for which the sum of the two corner numbers on one side of any edge is equal to the sum of the corner numbers on the other side of that edge. For each edge $e$, the latter is a linear equation of the form 
$$
L_e(u)=(u_{\theta_1}+u_{\theta_2})-(u_{\theta_3}+u_{\theta_4})=0,
$$ 
where $\theta_i$ are the corners adjacent\footnote{For the double edge $e$ of a folded triangle, two of these angles are the same and $L_e$ takes the form $L_e(u)=u_{\theta_1}-u_{\theta_3}$.}  to $e$.
In other words, $\iota(\La)$ is defined by a set of $N$ linear equations in $\Z_{\geq 0}^{2N}$, or $\iota(\La)$ is the intersection of a half-dimensional integrally defined subspace in $\R^{2N}$ with $\Z_{\geq 0}^{2N}$.\\

Let $\Lambda_\R=\Lambda\otimes \R\subset \R_{\geq 0}^E$. The last paragraph also shows the codimension-$1$ boundary faces of $\Lambda_\R$ are defined by equations $u_\theta=0$ for some $\theta\in \Theta$. In other words, we have the following.    

\begin{corollary}\label{interior-point}
The interior of $\La$ consists of those admissible colorings $v$ for which none of the triangle inequalities in (\ref{admissible_e}) is satisfied.
In particular, the admissible coloring $(2,2,\ldots,2)\in \La$ is an interior point of $\La$.
\end{corollary}

 \noindent
Geometric intersection numbers with any collection $S$ of curves on $\Si_{g,n}$   defines an admissible coloring 
$$
v(S)= (i(e,S))_{e\in E} \in \Lambda.
$$ 
By \cite[Prop.~3.2]{AF}, the map 
\bEq{MC=Color}
\mf{m}\lra v(\mf{m})
\eEq
that sends any multi-curve to the corresponding admissible coloring \underline{identifies} the set of multi-curves on $\Si_{g,n}$ with $\Lambda$. We will denote the inverse of~(\ref{MC=Color}) by $v\lra \mf{m}(v)$, for all $v\in \La$. It is clear that 
\bEq{deDe=inner-product}
\de_{\De}(\mf{m})=v(\mf{m})\cdot \vec{1}_E \in \Z_{\geq 0},
\eEq
where $\vec{1}_E=(1_e)_{e\in E}\in \Z^E$ ($v\cdot \vec{1}_E=\sum_{e} v_e$ for all $v\in \Z^E$).
Also, note that $2 \,v\cdot \vec{1}_E= \iota(v)\cdot \vec{1}_\Theta$; i.e.,
$$
i(\mf{m}(v),\De),\quad v\cdot \vec{1}_E,\quad \tn{and}\quad \frac{1}{2} \iota(v)\cdot \vec{1}_\Theta
$$ 
define the same filtration on $\tn{Sk}_{g,n}$.  The multi-curve $\mf{m}(v)$ associated to $v\in \La$ is constructed in the following way. For each $e\in E$ with $v_e>0$, place $v_e$ parallel line segments orthogonal to $e$, each of them intersecting $e$ at one point; see Figure~\ref{LS_fig}. The admissibility condition means there exists a unique way to connect those line segments to one another in each triangle $t$ to make a multi-curve $c$ whose geometric intersection numbers with the edges of the ideal triangles realizes the coloring. The number of arcs formed at each angle $\theta$ will be $u_\theta$.
\begin{corollary}\label{theta-to-peripheral_rmk}
If $u_\theta\neq 0$ for all $\theta$ adjacent to a puncture $p_i$, then the multi-curve $\mf{m}(v)$ contains the peripheral curve $a_i$ as one of its components. 
\end{corollary}
\begin{figure}[t]
\begin{pspicture}(-10,-1)(0,0)
\psset{unit=.7cm}


\psline[linecolor=red](-7.1,-1.3)(-7.1,-.7)
\psline[linecolor=red](-6.9,-1.3)(-6.9,-.7)
\psline[linecolor=red](-7.3,-1.3)(-7.3,-.7)
\psline[linecolor=red](-6.7,-1.3)(-6.7,-.7)

\psline[linecolor=red](-6.8,-.2)(-6.2,.2)
\psline[linecolor=red](-6.7,-.35)(-6.1,.05)

\psline[linecolor=red](-7.8,.2)(-7.2,-.2)
\psline[linecolor=red](-7.9,.05)(-7.3,-.35)

\psline(-7,1)(-8,-1)
\psline(-7,1)(-6,-1)
\psline(-8,-1)(-6,-1)

\psline[doubleline=true,doublesep=1.5pt]{->}(-4,0)(-2.5,0)

\psline[linecolor=red](-0.1,-1.3)(-0.1,-.26)
\psline[linecolor=red](0.1,-1.3)(0.1,-.19)
\psline[linecolor=red](-0.3,-1.3)(-0.3,-.34)
\psline[linecolor=red](.3,-1.3)(.3,-.34)

\psline[linecolor=red](.1,-.2)(.733,.2)
\psline[linecolor=red](.3,-.35)(.9,.05)

\psline[linecolor=red](-0.8,.2)(-0.1,-.263)
\psline[linecolor=red](-0.9,.05)(-0.3,-.35)

\psline(0,1)(-1,-1)
\psline(0,1)(1,-1)
\psline(-1,-1)(1,-1)

\end{pspicture}
\caption{Line segments on the edges of a triangle $\{e,e',e''\}$ for $v\in\La$ satisfying $(v_e,v_{e'},v_{e''})=(4,2,2)$. }
\label{LS_fig}
\end{figure}

By the identification~(\ref{MC=Color}), the set $\cC(\De)$ of the generators of the monoid $\La$ corresponds to certain set of ``indecomposable curves" on $\Si_{g,n}$ (Note that every multi-curve is decomposable to curves). Theorem~\ref{indecomposable_th} in Section~\ref{indecomposable_sec} explicitly characterizes the generating set $\cC(\De)$ as ``barbell trees". Often, but not always, $\cC(\De)$ contains the set of peripheral curves $\{a_1,\ldots,a_n\}$. Nevertheless, $v(a_1),\ldots,v(a_n)$ are linearly independent and, by Corollary~\ref{interior-point}, since
$$
v(a_1)+\cdots+v(a_n)=2\times \vec{1}_E,
$$
$v(a_1)+\cdots+v(a_n)$ defines an interior point of $\La_\R$. By Theorem~\ref{indecomposable_th}, the set of generators of the real and rational cones $\La_\R$ and $\La_\Q=\La\otimes \Q$ belong to a sub-set $\cC_\Q(\De)\subset \cC(\De)$ of barbell trees called ``simple barbel trees".

\begin{remark}\label{Toric_rmk}
In the toric geometry terminology of \cite{CLS}, $\Lambda$ is a  {\bf saturated affine semi-group}. The projective toric variety associated to $\Lambda$ is the projective variety corresponding to the projective scheme
$$
\tn{Proj}(\C[\Lambda]),
$$ 
where the semigroup algebra $\C[\Lambda]$ is the vector space over $\C$ with $\Lambda$ as a basis. The multiplication on $\C[\Lambda]$ is induced by the semigroup structure and $\C[\Lambda]$ is $\Z$-graded\footnote{For the sake of defining Proj.} by $v\cdot \vec{1}_E$. More precisely, if we think of $\Z_{\geq 0}^{E}$ as the character lattice of a torus so that $v\in \Z_{\geq 0}^{E}$ gives the character $x^v$, then
$$
\C[\Lambda]=\lrc{\sum_{v\in \La}\la_v x^v\mid \la_v \in \C~\tn{and}~\la_v= 0~\tn{for all but finitely many}~v}
$$
is an $\Z^E$-graded algebra with multiplication induced by $x^v\cdot x^{v'}=x^{v+v'}$.
The map $\mf{m}\lra x^{v(\mf{m})}$ is not an algebra isomorphism between $\tn{Sk}_{g,n}$ and $\C[\Lambda]$. As the following argument shows, we can think of $\tn{Sk}_{g,n}$ as an algebra generated (as a vector space) by $v\in \Lambda$ such that product structure $\ast$ is of the form 
\bEq{expansion-ast}
x^v\ast x^{v'}=\la_{v,v'} x^{v+v'} + \tn{``lower order terms"}. 
\eEq
Furthermore, we can use the algebra isomorphism in Remark~\ref{spin_rmk} to make $\la_{v,v'} =1$. As a corollary of this argument, with notation as in (\ref{GAlg_e}), we obtain a non-canonical isomorphism 
$$
\C[\cR_{g,n}]_{\de}=(\tn{Sk}_{g,n})_{\de}\cong \C[\Lambda].
$$
We conclude that the boundary divisor $\cD_{g,n}$ of the compactification
$$
\ov\cR_{g,n}=\tn{Proj}\big(\C[\cR_{g,n}][u]_{\de}\big)
$$ 
is an irreducible (projective) toric variety whose moment polytope $P_\De$ is the intersection of $\La_\R$ with the hyperplane $v\cdot \vec{1}_E\equiv \tn{constant}>0$ in $\R_{\geq 0}^N$. The vertices of $P_\De$ correspond to extremal rays in $\La_\R$ and thus to indecomposable curves in $\cC_\Q(\De)$. The faces of $P_\De$ are defined by equations $u_\theta=0$, for all $\theta\in \Theta$.
\end{remark}

In order to state the precise form of (\ref{expansion-ast}) and prove that $\de_\De$ is a semi-degree, we will utilize a strong partial order and a total order defined below. The second one requires an ordering of $E$.\\

For every two vectors $v,v'\in \Z^E$, we say $v$ is {\bf strongly smaller} than $v'$, and write $v\prec v'$, iff $0\neq (v'-v)\in\Z_{\geq 0}^E$. On the other hand, an ordering on $E$ defines a lexicographic (total) ordering on $\Z_{\geq 0}^E\cong \Z_{\geq 0}^N$. We will denote this total order by $<$. It is clear that $v\preceq v'$ iff $v\leq v'$ with respect to all lexicographic orders on $\Z_{\geq 0}^E$. \\

By the identification (\ref{MC=Color}), for two multi-curves $\mf{m}$ and $\mf{m}'$ on $\Si_{g,n}$, we write $\mf{m}\prec \mf{m}'$ or $\mf{m}< \mf{m}'$ iff the corresponding admissible coloring vectors $v=v(\mf{m})$ and $v'=v(\mf{m}')$ satisfy $v\prec v'$ or $v<v'$, respectively. We extend the total order $<$ on multi-curves to linear combinations in $\tn{Sk}_{g,n}$ in the following way. For every $f= \sum_s \la_s \mf{m}_s \in \tn{Sk}_{g,n}$, among those terms $\mf{m}_s$ in\footnote{i.e., $\la_s\neq 0$.} $f$ for which $\de_\De(\mf{m}_s)$ is maximal\footnote{This extra condition is not required in \cite{AF}.} (i.e. $\de_\De(\mf{m}_s)=\de_\De(f)$), we define the the {\bf lead term} $\tn{ld}(f)$ to be the unique largest multi-curve $\mf{m}_s$ with respect to $<$. We write $f<f'$ iff either $\de_\De(f)<\de_\De(f')$ or $\de_\De(f)=\de_\De(f')$ and $\tn{ld}(f)<\tn{ld}(f')$.  \\

\begin{theorem}(\cite[Thm~3.4]{AF})\label{Thm3.4}
For every two multi-curves $\mf{m}=\mf{m}(v)$ and $\mf{m'}=\mf{m}(v')$, we have 
$$
\mf{m}\ast \mf{m}'= (-1)^{i(\mf{m},\mf{m}')} \mf{m}(v+v')+\sum_{s} \la_s \mf{m}_s\qquad \tn{s.t.}\quad \mf{m}_s\prec \mf{m}(v+v')\quad \forall~s.
$$
The first term $\mf{m}(v+v')$ on right is called the \textbf{geometric sum} of $\mf{m}$ and  $\mf{m}'$ in $\tn{Sk}_{g,n}$ and will be denoted by $\mf{m}\,\#\, \mf{m}'$. 
\end{theorem}
By~(\ref{deDe=inner-product}), $\mf{m}\,\#\, \mf{m}'$ is the unique term realizing $\de_\De(\mf{m}\ast \mf{m}')$. The explicit  description of the geometric sum $\mf{m}\,\#\, \mf{m}'$ is as follows. First, we can assume that the crossing points of $\mf{m}$ and  $\mf{m}'$ happen in the interior of the triangles of $\De$. The edges of $\De$ cut each of $\mf{m}$ and  $\mf{m}'$ into a sequence of line segments. If $\al$ and $\beta$ are two line segments of $\mf{m}$ and  $\mf{m}'$ inside a triangle $\{e,e',e''\}$ crossing at one point $q$, one of the edges, say $e$, makes a triangle with $\al$ and $\beta$; see Figure~\ref{littleT}. Among the two possible smoothings of $q$, there is one and only one  (the middle picture in Figure~\ref{littleT}) that includes a bigon bounded on one side by $e$. Resolving this bigon decreases the $e$-factor of the sum vector $v+v'$. Choosing the Right diagram at each crossing, \cite[Thm~3.4]{AF} shows that we obtain a multi-curve that corresponds to $v+v'$ and has no trivial components. Since we are working with the classical limit $A=-1$, by the second relation in (\ref{skein-relation}),  the coefficient of $\mf{m}\# \mf{m}'$ is $(-1)^{i(\mf{m},\mf{m}')}$ (The coefficient will be $1$ if use the classical limit $A=1$).\\

\begin{figure}
\begin{pspicture}(-11,-1.3)(0,1.2)
\psset{unit=1cm}

\pscircle*(-7,1){0.1}
\pscircle*(-8,-1){0.1}
\pscircle*(-6,-1){0.1}

\psline(-7,1)(-8,-1)
\psline[linecolor=red](-7.25,.5)(-6.5,-1) 
\rput(-6.6,-.4){\small $\al$}

\psline(-7,1)(-6,-1)
\psline[linecolor=red](-6.75,.5)(-7.5,-1)
\rput(-7.35,-.4){\small $\beta$}

\psline(-8,-1)(-6,-1)\rput(-7,-1.2){\small $e$}\rput(-7,0.3){\small $q$}

\psline[doubleline=true,doublesep=1.5pt]{->}(-5.5,0)(-4,0)
\rput(-4.5,0.5){\small two resolutions}

\pscircle*(-1,1){0.1}
\pscircle*(-2,-1){0.1}
\pscircle*(0,-1){0.1}
\psline(-1,1)(-2,-1)

\psline(-1,1)(0,-1)
\psline(-2,-1)(0,-1)\rput(-1,-1.2){\small $e$}

 \psline[linearc=.2,linecolor=red](-0.5,-1)(-1,0)(-1.5,-1)
 \psline[linearc=.2,linecolor=red](-1.25,.5)(-1,0)(-0.75,.5)

\pscircle*(2,1){0.1}
\pscircle*(1,-1){0.1}
\pscircle*(3,-1){0.1}
\psline(2,1)(1,-1)
\psline(2,1)(3,-1)
\psline(1,-1)(3,-1)\rput(2,-1.2){\small $e$}
 \psline[linearc=.4,linecolor=red](2.5,-1)(2,0)(2.25,.5) \psline[linearc=.4,linecolor=red](1.75,.5)(2,0)(1.5,-1)

\end{pspicture}
\caption{Left- two line segments crossing in a triangle. Middle- the resolution that creates a bigon. Right- the resolution that defines the geometric sum. }
\label{littleT}
\end{figure}

We are now ready to prove that $\de_\De$ is a projective semi-degree and finish the proof of Theorem~\ref{sk-comp_th}. For arbitrary elements 
$$
f_1=\sum_{s} \la_{1,s} \mf{m}_{s} \quad \tn{and}\quad f_2=\sum_t\la_{2,t} \mf{m}_{t},
$$
the product $f_1\ast f_2$ decomposes as
$$
f_1\ast f_2= \sum_{s,t} \la_{1,s}\la_{2,t} \,\mf{m}_{s}\ast \mf{m}_{t}=\sum_{s,t,r} \la_{1,s}\la_{2,t} \la_{s,t,r}\,\mf{m}_{s,t,r}.
$$
For $i=1,2$, we may assume $f_i$ is $\de_\De$-homogenous, i.e.
$$
\aligned
&\de_{\De}(\mf{m}_{s})=\de_{\De}(f_1)\qquad  \forall\;s~~\tn{s.t}~~\la_{1,s}\neq 0,\\
&\de_{\De}(\mf{m}_{t})=\de_{\De}(f_2)\qquad  \forall\;t~~\tn{s.t}~~\la_{1,t}\neq 0,
\endaligned
$$
because, by the Theorem~\ref{Thm3.4}, the multi-curves $\mf{m}_{s,t,r}$ arising from the products of  multi-curves with $\de_\De$ below maximum satisfy $\de_\De(\mf{m}_{s,t,r})<\de_\De(f)+\de_\De(f').$\\

With this assumption, let $\mf{m}^1$ denote the lead term of $f_1$, $\mf{m}^2$ denote the lead term of $f_2$ and $\mf{m}^{12}=\mf{m}^1\#\mf{m}^2$ denote the geometric sum of $\mf{m}^{1}$ and  $\mf{m}^{2}$. First, by the Theorem~\ref{Thm3.4} we have 
$$
\de_{\De}(\mf{m}^{12})=\de_{\De}(\mf{m}^{1})+\de_{\De}(\mf{m}^{2})= \de_{\De}(f_1)+\de_{\De}(f_2).
$$
On the other hand, since $\La$ is a well-ordered monoid, $\mf{m}^{12}$ is larger in order $>$ from any other geometric sum $\mf{m}_{1,s}\#\mf{m}_{1,t}$ other than itself. We conclude that $\de_\De(f_1\ast f_2)=\de_{\De}(\mf{m}^{12})=\de_{\De}(f_1)+\de_{\De}(f_2)$.\\

Finally, it is clear that the filtration $\cF$ corresponding to $\de_\De$ satisfies $F_0=\C$. That  $\C[Z][u]_{\cF}$ is finitely-generated follows from the fact that $\La$ is a finitely generated integral cone in the positive quadrant. We conclude that $\de_\De$ is a projective semi-degree.\\

By Theorem~\ref{Thm3.4}, the product map $\ov\ast$ on the associated graded algebra $(\tn{Sk}_{g,n})_{\de_\De}$ is induced by
\bEq{AssGradProd_e}
\mf{m}\,\ov\ast\, \mf{m}'= (-1)^{i(\mf{m},\mf{m}')} \mf{m}\#\mf{m}'
\eEq
on the (vector space) basis elements.
Under the identification~(\ref{MC=Color}), the latter matches with the product structure on $\C[\La]$ in Remark~\ref{Toric_rmk}, except for the extra sign factor $(-1)^{i(\mf{m},\mf{m}')}$. The sign factor will become $1$ if we use the classical limit $A=1$ instead of $A=-1$. Therefore, the isomorphism $\Phi$ in Remark~\ref{spin_rmk} induces an isomorphism $(\tn{Sk}_{g,n})_{\de_\De}\cong \C[\Lambda]$. We conclude that the boundary divisor 
$$
\cD_{g,n}=\tn{Proj}\big((\tn{Sk}_{g,n})_{\de_\De}\big)
$$
of the projective compactification $\ov\cR_{g,n}= \ov\cR_{g,n}^{\de_\De}$ determined by the semi-degree $\de_\De$ 
is isomorphic to the irreducible (projective) toric variety $\tn{Proj}(\C[\La])$ whose moment polytope $P_\De$ is the intersection of $\La_\R$ with the hyperplane $v\cdot \vec{1}_E\equiv \tn{constant}>0$ in $\R_{\geq 0}^N$.  This finishes the proof of Theorem~\ref{sk-comp_th}.\qed

\subsection{Relative compactification (Proof of Theorem~\ref{Rel-sk-comp_th})}\label{Pfof-Rel-sk-comp_th}

In this section, first, we study the closure $\ov\cR_{g,n,\mft}$ of each fiber $\cR_{g,n,\mft}\subset \cR_{g,n}$ in the compactification $\ov\cR_{g,n}$ defined by an ideal triangulation $\De$. We show that the boundary divisor
$$
\cD_{g,n,\mft}=\ov\cR_{g,n,\mft}\cap \cD_{g,n}
$$
is a fixed toric sub-variety $\cD_{g,n}^{\tn{rel}}\subset \cD_{g,n}$ independent of the choice of $\mft\in \C^n$. Intuitively, in order to separate the compactified fibers $\ov\cR_{g,n,\mft}$ and construct a relative compactification as in (\ref{pinbar_e}), we need to pass to a birational modification (e.g. blowup) of $\ov\cR_{g,n}$ along $\cD_{g,n}^{\tn{rel}}$. This corresponds to changing the filtration on $\tn{Sk}_{g,n}$. We will define an explicit modification $\wt\de_\De$ of the semi-degree $\de_\De$ which defines the desired relative compactification. The new degree like function $\wt\de_\De$ can not be a semi-degree because $\cD_{g,n}^{\tn{rel}}$ is not irreducible. We show (indirectly) that it is a sub-degree. This implies that the resulting relative compactification is also normal; meaning that for instance  $\ov\cR_{g,n,\mft}$ does not have nodal singularities along $\cD_{g,n}^{\tn{rel}}$.\\

For $\mft=(t_1,\ldots,t_n)\in \C^n$, the inclusion $\iota_{\mft}\colon \cR_{g,n,\mft} \lra \cR_{g,n}$ corresponds to a projection map 
$$
\iota_\mft^*\colon \tn{Sk}_{g,n}\lra \tn{Sk}_{g,n,\mft}
$$ 
to the ring of regular functions $\tn{Sk}_{g,n,\mft}\defeq \C[\cR_{g,n,\mft}]$ of $\cR_{g,n,\mft}$. If $\cF$ is the filtration on $\tn{Sk}_{g,n}$ corresponding to $\de_\De$, the filtration $\cF_\mft$ on $\C[\cR_{g,n,\mft}]$ that defines the closure $\ov\cR_{g,n,\mft}$ is the image of $\cF$ under $\iota^*_\mft$. Recall from the last paragraph of Section~\ref{skein-algebra_sec} that $\tn{Sk}_{g,n,\mft}$ can be identified with the subspace $\tn{Sk}^{\tn{rel}}_{g,n}\subset \tn{Sk}_{g,n}$ generated by multi-curves that do not include any peripheral curve, the degree-like function defining $\cF_\mft$ is simply the restriction of $\de_\De$ to $\tn{Sk}^{\tn{rel}}_{g,n}$ (which is independent of $\mft$), but the product structure $\ast_{\mft}$ induced on $\tn{Sk}^{\tn{rel}}_{g,n}$ will be different from $\ast$ and depends on $\mft$. \\

\textbf{Proof of Theorem~\ref{Rel-sk-comp_th}.} For each multi-curve $\mf{m}\in \tn{Sk}_{g,n}$, define $\la_{\mft}(\mf{m})~\ov{\mf{m}}\in \tn{Sk}^{\tn{rel}}_{g,n}$ so that $\ov{\mf{m}}$ is the multi-curve obtained by removing all the peripheral curves from $\mf{m}$ and $\la_{\mft}(\mf{m})\in \C$ is the number obtained by exchanging $a_i$ for $t_i$. It is clear that 
$$
\de_\De(\mf{m})\geq \de_\De(\ov{\mf{m}})
$$
with the equality happening if and only if $\mf{m}\in \tn{Sk}^{\tn{rel}}_{g,n}$ (in which case $\la_\mft(\mf{m})=1$).
On the other hand, recall from Theorem~\ref{Thm3.4} that 
$$
\mf{m}\ast \mf{m}'= (-1)^{i(\mf{m},\mf{m}')} \mf{m}\#\mf{m}'+\sum_{s} \la_s \mf{m}_s\qquad \tn{s.t.}\quad \mf{m}_s\prec \mf{m}\#\mf{m}'\quad \forall~s.
$$
For $\mf{m},\mf{m}'\in \tn{Sk}^{\tn{rel}}_{g,n}$, it is possible that $\mf{m}\#\mf{m}'$ contains a peripheral curve. Therefore,
$$
\mf{m}\ast_{\mft} \mf{m}'= (-1)^{i(\mf{m},\mf{m}')} \la_{\mft}(\mf{m}\#\mf{m}')\, \ov{\mf{m}\#\mf{m'}}+\sum_{s}  \la_s \la_{\mft}(\mf{m}_s) \ov{\mf{m}}_s
$$ 
and 
$$
\de_\De(\mf{m}\ast_{\mft} \mf{m})\leq \de_\De(\mf{m})+\de_\De(\mf{m}')
$$ 
with the equality happening if and only if $\mf{m}\#\mf{m'}\in \tn{Sk}^{\tn{rel}}_{g,n}$. The divisor $\cD_{g,n,\mft}$ is defined by the associated graded algebra of $\tn{Sk}^{\tn{rel}}_{g,n}$. By the argument above, the product structure $\ov\ast_\mft$ induced on the graded algebra $[\tn{Sk}^{\tn{rel}}_{g,n}]_{\de_\De}$ is determined by 
$$
\mf{m}\ov\ast_{\mft} \mf{m}'=\begin{cases}
(-1)^{i(\mf{m},\mf{m}')} \mf{m}\#\mf{m}' &\tn{if}~~\mf{m}\#\mf{m'}\in \tn{Sk}^{\tn{rel}}_{g,n}\\
0 & \tn{otherwise}
\end{cases}\quad \forall~\mf{m},\mf{m}'\in \tn{Sk}_{g,n}^{\tn{rel}}
$$
on the basis elements. In particular, $\ov\ast_{\mft}$ and therefore $\cD_{g,n,\mft}$ is independent of $\mf{t}$. The ideal defining $\cD_{g,n,\mft}\subset \cD_{g,n}$ is generated by the characters corresponding to peripheral curves. Therefore, $\cD_{g,n,\mft}\subset \cD_{g,n}$ is a toric subvariety.\\

Recall that the projection
$$
\Pi_{g,n}\colon \cR_{g,n}\lra \C^n
$$
is dual to considering $\tn{Sk}_{g,n}$ as an $\C[a_1,\ldots,a_n]$-module, where $\C[a_1,\ldots,a_n]$ is the sub-algebra generated by peripheral curves. Thus, relative compactifications of $\cR_{g,n}$ correspond to filtrations $\wt\cF$ with $\wt{F}_0=\C[a_1,\ldots,a_n]$. Therefore, to define a degree-like function $\wt\de_\De$/filtration $\wt\cF$ that produces a relative compactification of $\cR_{g,n}$ with fibers the same as $\ov\cR_{g,n,\mft}$ above, we need a degree-like function $\wt\de_\De$ satisfying 
$$
\wt\de_\De(\mf{m})=\de_\De(\ov{\mf{m}}) \qquad \forall~\mf{m}\in \tn{Sk}^{\tn{rel}}_{g,n}, \quad \wt\de_\De(a_i)=0\quad i=1,\ldots,n.
$$

For every arc/edge $e\in \De$, let $\wt{e}$ be the unique (up to isotopy) curve which is the boundary of a sufficiently small thickening of $e$. Let $\wt\De$ be the collection of these curves. It follows immediately from Bigon Criterion that 
$\wt\de_\De\defeq \frac{1}{2}\de_{\wt{\De}}$ satisfies 
$$
\wt\de_\De(c)=
\begin{cases}
\de_\De(c) &\tn{if}~~c\neq a_1,\ldots,a_n\\
0 & \tn{otherwise},
\end{cases}
$$
for every curve $c$. It is immediate from this definition that $\wt{F}_0=\C[a_1,\ldots,a_n]$ and the compactification corresponding to the degree-like function  $\wt\de_\De$ is a relative compactification whose fibers are~$\ov\cR_{g,n,\mf{t}}$.\\

It is just left to show that $\wt\de_\De$ is a sub-degree. By Remark~\ref{sub-condition}, we will simply show that 
$$
\wt\de_\De(f^k)=k\wt\de_\De(f)\qquad \forall~f\in \tn{Sk}^{\tn{rel}}_{g,n}
$$
to make the desired conclusion. Suppose $f= \sum_{s} \la_s \mf{m}_s \in \tn{Sk}^{\tn{rel}}_{g,n}$ and $\mf{m}_0=\tn{ld}(f)$ is the leading term. Then, the leading term of $f^k$ is $\mf{m}_0^k\in \tn{Sk}^{\tn{rel}}_{g,n}$.  Therefore, with respect to the product $\ast_\mft$, too, the term $\mf{m}_0^k$ appears in the expansion of $f\ast_\mft \cdots \ast_\mft f$. This finishes the proof of Theorem~\ref{Rel-sk-comp_th}.

\begin{remark}
It would be interesting to find geometrically defined semi-degrees $\de_1,\ldots,\de_m$ so that $\wt\de_\De=\max(\de_1,\ldots,\de_m)$. 
\end{remark}

\subsection{Indecomposable curves and barbell trees}\label{indecomposable_sec}
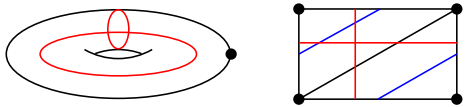
In this section, we give a complete characterization of the set $\cC=\cC(\De)$ of ``indecomposable curves"  that generate the monoid $\Lambda$ associated to an ideal triangulation $\De$ of $\Si_{g,n}$. We also describe the generating set of the cone $\La_\R$. For instance, Figure~\ref{gn11-fig} illustrates an ideal triangulation of $\Si_{1,1}$ and the corresponding set of generators of $\La.$

\begin{figure}[h]
\begin{pspicture}(4,.6)(18,2.4)
\psset{unit=.3cm}

\psellipse[linewidth=.07](35,5)(5,2)
\psellipse[linecolor=red,linewidth=.07](35,5)(3.5,1)

\psarc[linewidth=.07](35,7.8){3}{240}{300}
\psarc[linewidth=.07](35,2.2){3}{70}{110}
\psellipse[linecolor=red,linewidth=.07](35,6.1)(0.5,.9)

\pscircle*(40,5){.25}

\psline[linewidth=.07](43,3)(43,7)(50,7)(50,3)(43,3)
\psline[linewidth=.07](43,3)(50,7)
\psline[linewidth=.07,linecolor=blue](46.5,3)(50,5)
\psline[linewidth=.07,linecolor=blue](46.6,7)(43,5)

\pscircle*(43,3){.25}
\pscircle*(43,7){.25}
\pscircle*(50,7){.25}
\pscircle*(50,3){.25}
\psline[linecolor=red,linewidth=.07](45.5,3)(45.5,7)
\psline[linecolor=red,linewidth=.07](43,5.5)(50,5.5)

\end{pspicture}
\caption{Red and blue curves generate $\La\subset \N^3$. Black lines in the rectangle are the edges of the triangulation.}
\label{gn11-fig}
\end{figure}
Similarly, Figures~\ref{gn11-fig} Left and Right illustrate two ideal triangulations of $\Si_{0,4}$ and the corresponding set of generators of $\La.$ The one on the left has $7$ and the one on the right has $8$ generators. Readers not interested in the combinatorial properties of $\La$ may skip this section.

\begin{figure}[h]
\begin{pspicture}(-10,-2)(0,1.5)
\psset{unit=.7cm}

\pscircle*(-7,1){0.1} \pscircle[linecolor=blue](-7,1){.25}
\pscircle*(-8,-1){0.1}\pscircle[linecolor=blue](-8,-1){.25}
\pscircle*(-6,-1){0.1}\pscircle[linecolor=blue](-6,-1){.25}
\pscircle[linecolor=blue](-7,-.3){2.1}
\psellipse[rot=-30,linecolor=red](-7.5,0)(.5,1.8)
\psellipse[rot=30,linecolor=red](-6.5,0)(.5,1.8)
\psellipse[rot=90,linecolor=red](-7,-1)(.5,1.8)
\psline(-7,1)(-7,2.5)
\psline(-8,-1)(-9.5,-2)
\psline(-6,-1)(-4.5,-2)

\psline(-7,1)(-8,-1)
\psline(-7,1)(-6,-1)
\psline(-8,-1)(-6,-1)


\pscircle*(1,1){0.1} \pscircle[linecolor=red](1,1){.25}
\pscircle*(0,-1){0.1}\pscircle[linecolor=red](0,-1){.25}
\pscircle*(2,-1){0.1}\pscircle[linecolor=red](2,-1){.25}
\pscircle[linecolor=red](1,-.3){2.1}
\psellipse[rot=-30,linecolor=red](0.5,0)(.5,1.8)
\psellipse[rot=30,linecolor=red](1.5,0)(.5,1.8)
\psellipse[rot=90,linecolor=red](1,-1)(.5,1.8)
\psline(1,1)(1,2.5)
\psline(0,-1)(-1.5,-2)
\psline(2,-1)(3.5,-2)

\psline(1,1)(0,-1)
\psline(1,1)(2,-1)
\psline(1,1)(1,-2.8)

\psarc[linecolor=blue](0,-1){.35}{180}{360}
\psarc[linecolor=blue](2,-1){.35}{180}{360}
\psline[linecolor=blue](-0.35,-1)(-.35,.9)
\psline[linecolor=blue](1.65,-1)(1.65,1.3)
\psline[linecolor=blue](2.35,-1)(2.35,.9)
\psline[linecolor=blue](.35,-1)(.35,1.3)
\psarc[linecolor=blue](1,.9){1.35}{0}{180}
\psarc[linecolor=blue](1,1.3){.65}{0}{180}

\end{pspicture}
\caption{In each diagram red and blue curves generate $\La\subset \N^6$. Black lines are the edges of a triangulation, with the half-open edges meeting at the fourth punctured point located at $\infty$.}
\label{n4EX}
\end{figure}

\begin{figure}
\begin{pspicture}(2,-1.5)(5,0)
\psset{unit=.8cm}

 \pscircle*(7,-1){0.1}
\pscircle*(9,-1){0.1}
\psline(7,-1)(9,-1)
\psarc(8,-2){1.4}{45}{135}
\psarc(8,0){1.4}{225}{315}

 \pscircle*(12,-.85){0.1}
 \pscircle*(13,0){0.1}
  \pscircle*(11,0){0.1}
 \pscircle*(12,-2){0.1}
\psline(12,-2)(12,-.85)(13,0)(12,-2)
\psline(13,0)(11,0)(12,-.85)
\psline(11,0)(12,-2)

 \pscircle*(18,0){0.1}
  \pscircle*(16,0){0.1}
 \pscircle*(18,-2){0.1}
  \pscircle*(16,-2){0.1}
\psline(18,0)(16,0)(16,-2)(18,-2)(18,0)
\psarc(17,-1){1.4}{-45}{45}
\psarc(17,-1){1.4}{135}{225}

\end{pspicture}
\caption{Left-dual graph of Figure~\ref{gn11-fig}; Middle-dual graph of Figure~\ref{n4EX}.Left; and Right-dual graph of Figure~\ref{n4EX}.Right}
\label{DG_fig}
\end{figure}
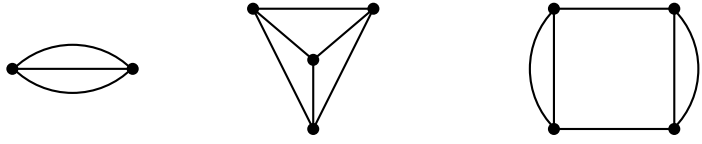

With notation as before, suppose $\De$ is an ideal triangulation of $\Si_{g,n}$ with the set of triangles $T$ and edges $E$.  The dual graph $\Gamma$ of $\De$ is a trivalent graph which has a vertex for each triangle $t$ in $T$ and a dual edge $e^\perp$ for each edge $e$ in $E$ connecting the triangles on the two sides of $e$. Here, $e^\perp$ corresponds to a line segment transversely intersecting $e$ in one point. Thus, the vertices of $\Gamma$ are indexed by $T$ and the edges are still indexed by $E$. Double edges of folded triangles generate loops\footnote{i.e., edges ending on the same vertex at the two ends.} in $\Gamma$. Figure~\ref{DG_fig} illustrates the dual graphs of the triangulations in Figures~\ref{gn11-fig} and  \ref{n4EX}. By the one-to-one correspondence between the edges of $\De$ and $\Gamma$, every admissible coloring $v\in \La$ can also be thought of as a coloring of the edges of $\Gamma$ satisfying the admissibility conditions
\bEq{admissible_e2}
2\mid v_{e_1^\perp}+v_{e_2^\perp}+v_{e_3^\perp}\qquad \tn{and}\qquad v_{e_i^\perp}\leq v_{e_j^\perp}+v_{e_k^\perp}\quad \forall~i=1,2,3,~\{j,k\}=\{1,2,3\}-\{i\},
\eEq  
at each vertex $t$ of $\Gamma$. We will continue to write $v_e$ instead of $v_{e^\perp}$, for simplicity.\\

An admissible coloring $v\in \La$ is called {\bf indecomposable} whenever $v=v'+v''$, with $v',v''\in \La$, implies $v'=0$ or $v''=0$.
Recall from (\ref{MC=Color}) that there is one-to-one correspondence between the set of multi-curves and admissible colorings $v\in  \La$. Clearly, the generators of $\La$ are the indecomposable colorings and correspond to curves (i.e., multi-curves with only one component). The curves corresponding to indecomposable colorings will be called {\bf indecomposable curves}. We will characterize the set $\cC(\De)$ of indecomposable admissible colorings/curves.  \\

For each $v\in \La$, define the \textbf{support of} $v$ to be the positively-weighted/colored subgraph $\Gamma_v$ of $\Gamma$ consisting of the edges $e^\perp$ with $v_e>0$ and weighted/colored with $v_e\in \Z_+$. Each edge of $\Gamma$ has two ends, even if both ends of the edge lie at the same vertex. The ends of an  edge inherits a coloring from the color of the edge.   For each vertex $t\in \Gamma$, we say $v$ or $\Gamma_v$ has type $(a,b,c)$ at $t$ if the ends of the edges at $t$ in $\Gamma$ are assigned the colors $(a,b,c)$.  Note that by (\ref{admissible_e2}), $2\!\mid\! a+b+c$; i.e., the number of odd colors at each $t$ is even.

\bDf{barbell-tree}
We say a weighted/colored subgraph $G$ of $\Gamma$ is a \textbf{barbell tree} if $G$ is connected, the triples of colors at each vertex $t$ are one of the followings
\bEq{barbell-types}
(2,2,2),\quad (2,2,0), \quad (2,1,1),\quad (0,1,1), \quad (0,0,0),
\eEq
and the result of collapsing all the edges labeled $1$ is a tree. 
\eDf
\vskip.1in

Note that every barbell tree $G$ is the support graph of some admissible coloring $v\in \La$ because (\ref{admissible_e2}) is satisfied at any vertex $t$ of $\Gamma$. 
The terminology ``barbell tree" follows from the following geometric description of $G$.  Since $\Gamma$ is trivalent and each vertex contains none or two edges colored with $1$, the set of edges of weight $1$ form a disjoint union of embedded circles (we call it {\bf bells}) in $\Gamma$. The rest of the edges in $G$ are colored with $2$ and each connected component of them is a tree that connects some of the bells. None of the edges colored by $2$ is a loop and every such edge is connected to another edge on both ends. If we collapse all the bells to vertices, we get a graph $\ov{G}$ which we require to be a tree. We say a bell is a leaf-bell if the corresponding vertex in $\ov{G}$ is a leaf. In other words, either that bell is the entire $G$ or it is connected to a single edge of weight $2$ at a vertex of type $(2,1,1)$.  Every barbell tree that includes a edge of weight $2$ has at least two leaf-bells.

\bDf{simpleBB_dfn}
We say $G$ is a  \textbf{simple barbell tree} if is either a single bell (a cycle of $1$-colored edges) or it is a barbell, i.e., two bells attached by a chain of $2$-colored edges.
\eDf

 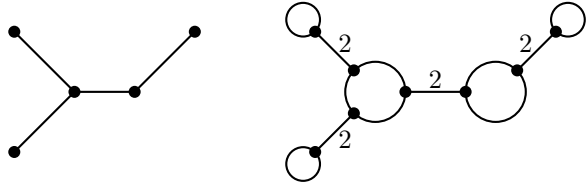
\begin{figure}
\begin{pspicture}(5,-2.0)(15,0)
\psset{unit=.8cm}

  \pscircle*(11,0){0.1}
 \pscircle*(11,-2){0.1}
  \pscircle*(12,-1){0.1}
   \pscircle*(13,-1){0.1}
    \pscircle*(14,0){0.1}
\psline(11,0)(12,-1)(13,-1)(14,0)
\psline(11,-2)(12,-1)

  \pscircle*(16,0){0.1}  \pscircle(15.8,0.2){0.28}
 \pscircle*(16,-2){0.1} \pscircle(15.8,-2.2){0.28}
  \pscircle(17,-1){0.5} 
  \pscircle*(17.5,-1){0.1}\pscircle*(16.643,-.643){0.1}\pscircle*(16.643,-1.357){0.1}
   \pscircle(19,-1){0.5}
     \pscircle*(18.5,-1){0.1}\pscircle*(19.357,-.643){0.1}
    \pscircle*(20,0){0.1}\pscircle(20.2,0.2){0.28}
\psline(16,0)(16.643,-.643)\rput(16.5,-.2){\small $2$}
\psline(16,-2)(16.643,-1.357)\rput(16.5,-1.8){\small $2$}
\psline(17.5,-1)(18.5,-1)\rput(18,-.8){\small $2$}
\psline(19.357,-.643)(20,0)\rput(19.5,-.2){\small $2$}

\end{pspicture}
\caption{A barbell tree (Right) and its underlying tree (Left).}
\label{Bmodel}
\end{figure}

\begin{example} 
Figure~\ref{Bmodel}-Right illustrates a barbell tree and Figure~\ref{Bmodel}-Left is the tree obtained by collapsing the bells into vertices. This example is not simple.
Figures~\ref{gn11-fig-EBT}, \ref{n4EX-EBT}, and~\ref{n4EX-2-EBT} show all the barbell trees embedded in the dual graphs of Figures~\ref{DG_fig}-Left, Middle, and Right, respectively. They are all simple.
 \begin{figure}
\begin{pspicture}(4,-1.5)(5,0)
\psset{unit=1cm}

 \pscircle*(7,-1){0.1}
\pscircle*(9,-1){0.1}
\psline[linecolor=red,linewidth=1.5pt](7,-1)(9,-1)
\psarc[linecolor=red,linewidth=1.5pt](8,-2){1.4}{45}{135}
\psarc(8,0){1.4}{225}{315}

 \pscircle*(11,-1){0.1}
\pscircle*(13,-1){0.1}
\psline(11,-1)(13,-1)
\psarc[linecolor=red,linewidth=1.5pt](12,-2){1.4}{45}{135}
\psarc[linecolor=red,linewidth=1.5pt](12,0){1.4}{225}{315}

 \pscircle*(15,-1){0.1}
\pscircle*(17,-1){0.1}
\psline[linecolor=red,linewidth=1.5pt](15,-1)(17,-1)
\psarc(16,-2){1.4}{45}{135}
\psarc[linecolor=red,linewidth=1.5pt](16,0){1.4}{225}{315}

\end{pspicture}
\caption{Red sub-graphs are the barbell trees of the dual graph in Figure~\ref{DG_fig}-Left.}
\label{gn11-fig-EBT}
\end{figure}
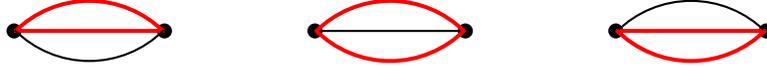

 \begin{figure}
\begin{pspicture}(6,-1.5)(12,0)
\psset{unit=.7 cm}

 \pscircle*(12,-.85){0.1}
 \pscircle*(13,0){0.1}
  \pscircle*(11,0){0.1}
 \pscircle*(12,-2){0.1}

\psline(12,-2)(12,-.85)(13,0)(12,-2)
\psline(13,0)(11,0)(12,-.85)
\psline(11,0)(12,-2)

\pscircle*(14.5,-.85){0.1}
 \pscircle*(15.5,0){0.1}
  \pscircle*(13.5,0){0.1}
 \pscircle*(14.5,-2){0.1}
\psline(14.5,-2)(14.5,-.85)(15.5,0)(14.5,-2)
\psline(15.5,0)(13.5,0)(14.5,-.85)
\psline(13.5,0)(14.5,-2)

\pscircle*(17,-.85){0.1}
 \pscircle*(18,0){0.1}
  \pscircle*(16,0){0.1}
 \pscircle*(17,-2){0.1}
\psline(17,-2)(17,-.85)(18,0)(17,-2)
\psline(18,0)(16,0)(17,-.85)
\psline(16,0)(17,-2)

\pscircle*(19.5,-.85){0.1}
 \pscircle*(20.5,0){0.1}
  \pscircle*(18.5,0){0.1}
 \pscircle*(19.5,-2){0.1}
\psline(19.5,-2)(19.5,-.85)(20.5,0)(19.5,-2)
\psline(20.5,0)(18.5,0)(19.5,-.85)
\psline(18.5,0)(19.5,-2)

\pscircle*(22,-.85){0.1}
 \pscircle*(23,0){0.1}
  \pscircle*(21,0){0.1}
 \pscircle*(22,-2){0.1}
\psline(22,-2)(22,-.85)(23,0)(22,-2)
\psline(23,0)(21,0)(22,-.85)
\psline(21,0)(22,-2)

\pscircle*(24.5,-.85){0.1}
 \pscircle*(25.5,0){0.1}
  \pscircle*(23.5,0){0.1}
 \pscircle*(24.5,-2){0.1}
\psline(24.5,-2)(24.5,-.85)(25.5,0)(24.5,-2)
\psline(25.5,0)(23.5,0)(24.5,-.85)
\psline(23.5,0)(24.5,-2)

\pscircle*(27,-.85){0.1}
 \pscircle*(28,0){0.1}
  \pscircle*(26,0){0.1}
 \pscircle*(27,-2){0.1}
\psline(27,-2)(27,-.85)(28,0)(27,-2)
\psline(28,0)(26,0)(27,-.85)
\psline(26,0)(27,-2)

\psline[linecolor=red,linewidth=1.5pt](11,0)(12,-.85)(13,0)(11,0)
\psline[linecolor=red,linewidth=1.5pt](14.5,-.85)(15.5,0)(14.5,-2)(14.5,-.85)
\psline[linecolor=red,linewidth=1.5pt](17,-2)(16,0)(17,-.85)(17,-2)
\psline[linecolor=red,linewidth=1.5pt](20.5,0)(18.5,0)(19.5,-2)(20.5,0)
\psline[linecolor=red,linewidth=1.5pt](23,0)(21,0)(22,-.85)(22,-2)(23,0)
\psline[linecolor=red,linewidth=1.5pt](25.5,0)(23.5,0)(24.5,-2)(24.5,-.85)(25.5,0)%
\psline[linecolor=red,linewidth=1.5pt](28,0)(27,-.85)(26,0)(27,-2)(28,0)

\end{pspicture}
\caption{Red sub-graphs are the barbell trees of the dual graph in Figure~\ref{DG_fig}-Middle.}
\label{n4EX-EBT}
\end{figure}

 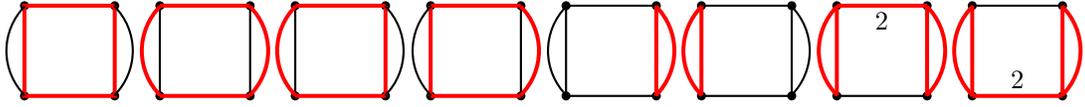
\begin{figure}[t]
\begin{pspicture}(-1.5,-1.5)(15,0)
\psset{unit=.6cm}

 \pscircle*(2,0){0.1}
  \pscircle*(0,0){0.1}
 \pscircle*(2,-2){0.1}
  \pscircle*(0,-2){0.1}
\psline[linecolor=red,linewidth=1.5pt](2,0)(0,0)(0,-2)(2,-2)(2,0)
\psarc(1,-1){1.4}{-45}{45}
\psarc(1,-1){1.4}{135}{225}

 \pscircle*(5,0){0.1}
  \pscircle*(3,0){0.1}
 \pscircle*(5,-2){0.1}
  \pscircle*(3,-2){0.1}
\psline[linecolor=red,linewidth=1.5pt](5,0)(3,0)
\psline(3,0)(3,-2)
\psline[linecolor=red,linewidth=1.5pt](3,-2)(5,-2)
\psline(5,-2)(5,0)
\psarc[linecolor=red,linewidth=1.5pt](4,-1){1.4}{-45}{45}
\psarc[linecolor=red,linewidth=1.5pt](4,-1){1.4}{135}{225}

 \pscircle*(8,0){0.1}
  \pscircle*(6,0){0.1}
 \pscircle*(8,-2){0.1}
  \pscircle*(6,-2){0.1}
\psline[linecolor=red,linewidth=1.5pt](8,0)(6,0)
\psline(6,0)(6,-2)
\psline[linecolor=red,linewidth=1.5pt](6,-2)(8,-2)
\psline[linecolor=red,linewidth=1.5pt](8,-2)(8,0)
\psarc(7,-1){1.4}{-45}{45}
\psarc[linecolor=red,linewidth=1.5pt](7,-1){1.4}{135}{225}

 \pscircle*(11,0){0.1}
  \pscircle*(9,0){0.1}
 \pscircle*(11,-2){0.1}
  \pscircle*(9,-2){0.1}
\psline[linecolor=red,linewidth=1.5pt](11,0)(9,0)
\psline[linecolor=red,linewidth=1.5pt](9,0)(9,-2)
\psline[linecolor=red,linewidth=1.5pt](9,-2)(11,-2)
\psline(11,-2)(11,0)
\psarc[linecolor=red,linewidth=1.5pt](10,-1){1.4}{-45}{45}
\psarc(10,-1){1.4}{135}{225}

 \pscircle*(14,0){0.1}
  \pscircle*(12,0){0.1}
 \pscircle*(14,-2){0.1}
  \pscircle*(12,-2){0.1}
\psline(14,0)(12,0)
\psline(12,0)(12,-2)
\psline(12,-2)(14,-2)
\psline[linecolor=red,linewidth=1.5pt](14,-2)(14,0)
\psarc[linecolor=red,linewidth=1.5pt](13,-1){1.4}{-45}{45}
\psarc(13,-1){1.4}{135}{225}

 \pscircle*(17,0){0.1}
  \pscircle*(15,0){0.1}
 \pscircle*(17,-2){0.1}
  \pscircle*(15,-2){0.1}
\psline(17,0)(15,0)
\psline[linecolor=red,linewidth=1.5pt](15,0)(15,-2)
\psline(15,-2)(17,-2)
\psline(17,-2)(17,0)
\psarc(16,-1){1.4}{-45}{45}
\psarc[linecolor=red,linewidth=1.5pt](16,-1){1.4}{135}{225}

 \pscircle*(20,0){0.1}
  \pscircle*(18,0){0.1}
 \pscircle*(20,-2){0.1}
  \pscircle*(18,-2){0.1}
\psline[linecolor=red,linewidth=1.5pt](20,0)(18,0)\rput(19,-.35){\small 2}
\psline[linecolor=red,linewidth=1.5pt](18,0)(18,-2)
\psline(18,-2)(20,-2)
\psline[linecolor=red,linewidth=1.5pt](20,-2)(20,0)
\psarc[linecolor=red,linewidth=1.5pt](19,-1){1.4}{-45}{45}
\psarc[linecolor=red,linewidth=1.5pt](19,-1){1.4}{135}{225}

 \pscircle*(23,0){0.1}
  \pscircle*(21,0){0.1}
 \pscircle*(23,-2){0.1}
  \pscircle*(21,-2){0.1}
\psline(23,0)(21,0)
\psline[linecolor=red,linewidth=1.5pt](21,0)(21,-2)
\psline[linecolor=red,linewidth=1.5pt](21,-2)(23,-2)\rput(22,-1.65){\small 2}
\psline[linecolor=red,linewidth=1.5pt](23,-2)(23,0)
\psarc[linecolor=red,linewidth=1.5pt](22,-1){1.4}{-45}{45}
\psarc[linecolor=red,linewidth=1.5pt](22,-1){1.4}{135}{225}

\end{pspicture}
\caption{Red sub-graphs are the barbell trees of the dual graph in Figure~\ref{DG_fig}-Right.}
\label{n4EX-2-EBT}
\end{figure}

\end{example}

\begin{theorem}\label{indecomposable_th} A curve $c$ or equivalently an admissible coloring $v\in \La$ is indecomposable if and only if the support $\Gamma_v$ of $v$ in $\Gamma$ is a barbell tree. In other words, there is a one-to-one correspondence between the set $\cC(\De)$ of generators of $\La$ and barbell trees embedded in~$\Gamma$. Similarly, there is a one-to-one correspondence between the set\footnote{The same statements hold with $\R$ in place of $\Q$.} $\cC_\Q(\De)$ of generators of $\La_\Q$ and simple barbell trees embedded in~$\Gamma$.\end{theorem}

\begin{proof}
Define a \textbf{chain} of $2$-colored edges in a barbell three to be a connected sequence of $2$-colored edges such that the vertex between every two consecutive edges has type $(2,2,0)$. Every maximal chain in a barbell tree has two ending vertices of the following possible types: $(2,2,2)$, $(2,1,1)$. \\

First, we show that if the support $\Gamma_v$ of $v\in \La$ is a barbell tree, then it is indecomposable. Suppose $v\in \La$ is decomposable and $v=v'+v''$ is a non-trivial decomposition of $v$. If $v_e=2$, then either $v'_e=v''_e=1$ or one of them, say $v'_e$ is $2$ and the other is $0$. Note that removing $e^\perp$ makes $\Gamma_v$ disconnected with two non-trivial connected components (if $v_e=2$, then $e^\perp$ is not a loop).\\

If $v'_e=2$,  by (\ref{barbell-types}) and the triangle inequality at vertices of $\Gamma_{v''}$ and $\Gamma_{v'}$ (i.e. the fact that the colorings $(2,0,0)$ and $(2,1,0)$ are not admissible), $\Gamma_{v'}$ includes the maximal chain of $2$-colored edges containing $e^\perp$. 
Therefore, for each maximal chain of $2$-colored edges 
\begin{itemize}
\item either the entire chain belongs to one of $\Gamma_{v'}$ or $\Gamma_{v''}$, 
\item or the weight $2$ is split evenly between $v$ and $v'$ on all of them.
\end{itemize}

First, we show that the second scenario is impossible. In the second case, the complement of this chain in both $\Gamma_{v'}$ or $\Gamma_{v''}$ will be disconnected. However, this is impossible because any curve in $\mf{m}(v')$ will define a loop in $\Gamma_{v'}$: for a path starting on this chain of $1$-colored edges, once it leaves this chain on one side, it can't comeback to it from the other side; that is a contradiction. Therefore, each maximal chain of $2$-colored edges belongs entirely to $\Gamma_{v'}$ or $\Gamma_{v''}$. If $e^\perp$ is the last edge of such a chain on $\Gamma_{v'}$ ending at the vertex $t$, then $t$ can be of one of the following types:
\bIt 
\item $(2,1,1)$  if $t$ is the connecting point with a bell;
\item $(2,2,2)$ if $t$ is the meeting point of three chains. 
\eIt
In the second case, all three chains meeting at $t$ must belong to $\Gamma_{v'}$; otherwise, the triangle inequality in one of $\Gamma_{v'}$ or $\Gamma_{v''}$ the will be violated at $t$.
In the first case, by the triangle inequality again, the edges with weight $1$ must belong to $\Gamma_{v'}$ as well. We can inductively make the same conclusion by moving along the $1$-weighted edges of the bell until we reach another vertex of type $(2,1,1)$ where another chain of $2$-weighted edges is connected to the bell. If the new $2$-colored edge $e^\perp$ found in the last triple belongs to $\Gamma_{v''}$, we get a contradiction because the same argument applied to $e^\perp$ shows that the $1$-colored edges should belong to $\Gamma_{v''}$. Therefore, unless one of $\Gamma_{v'}$ or $\Gamma_{v''}$ is empty, we always arrive at a vertex where the triangle inequality can not be satisfied!  This finishes the argument that every barbell tree defines an indecomposable admissible coloring.    \\

Conversely, suppose $v\in \La$ is indecomposable, we show that $\Gamma_v$ is a barbell tree.  Divide the edges $E$ into three types:
\bIt
\item $E_{\tn{odd}}=\{e\in E\colon v_e~\tn{is odd}\};$
\item $E_{\tn{even}}=\{e\in E\colon 2\mid v_e~\tn{and}~v_e\neq 0\};$
\item $E_0=\{e\in E\colon v_e=0\}$.
\eIt
The partition of $E$ into these three sets gives a decomposition of every vector $u\in \Z^E$ as
$$
u=(\al,\beta,\gamma)\in \Z^{E_{\tn{odd}}}\oplus \Z^{E_{\tn{even}}}\oplus \Z^{E_0}
$$ 
Define 
$$
v'=(\vec{1},\vec{2},\vec{0})\in \Z^{E_{\tn{odd}}}\oplus \Z^{E_{\tn{even}}}\oplus \Z^{E_0}.
$$
It is easy to see that both $v'$ and $v-v'$ are admissible colorings; see the sentence before Definition~\ref{barbell-tree}. Therefore, $v-v'$ must be zero. We conclude that every vertex in $v$ has one of the types listed in (\ref{barbell-types}). If the result of collapsing all the edges labeled $1$ is not a tree, one can subtract a cycle of edges from $v$, which defines an admissible coloring, and show that $v$ is decomposable. This finishes the proof of the first statement of Theorem~\ref{indecomposable_th}.\\

If $\Gamma_{v}$ is a barbell tree which is not simple, then either $\Gamma_v$ has a vertex $t$ of type $(2,2,2)$, or it has a bell with $2$ or more $2$-colored edges attached to it (this is a non-leaf bell). In the first case, the complement of $t$ is a disjoint union of $3$ connected components $G_1,G_2,G_3$. It is easy to check that each $G\setminus G_i$ is a barbell tree with one fewer vertex of type $(2,2,2)$. Let $v_i$ be the admissible coloring corresponding to $G\setminus G_i$. It is clear from the construction that 
$$
2v=v_1+v_2+v_3;
$$
i.e., the vector $v$ is a positive $\Q$-linear combination of $v_1,v_2,v_3$. \\

Similarly, in the second case, let $B$ be a bell in $\Gamma_v$ with $k\geq 2$ $2$-colored edges attached to it; see Figure~\ref{Bell_fig}. 
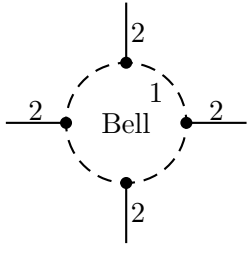
\begin{figure}
\begin{pspicture}(6,-2.3)(15,0)
\psset{unit=.8cm}


  \pscircle[linestyle=dashed](17,-1){1} \rput(17,-1){Bell}\rput(17.5,-.5){1}
 \pscircle*(18,-1){0.1}\psline(18,-1)(19,-1)\rput(18.5,-.8){2}
 \pscircle*(16,-1){0.1}\psline(16,-1)(15,-1)\rput(15.5,-.8){2}
  \pscircle*(17,0){0.1}\psline(17,0)(17,1)\rput(17.2,.5){2}
   \pscircle*(17,-2){0.1}\psline(17,-2)(17,-3)\rput(17.2,-2.5){2}
   
\end{pspicture}
\caption{A bell with four $2$-colored edges attached to it.}
\label{Bell_fig}
\end{figure}
Removing the edges in $B$, we get $k$ connected components, each of which contains one of these $2$-colored edges. Going around the circle, we name these connected components $G_1,\ldots,G_k$. For each $i=1,\ldots,k$, with $k+1\equiv 1$, define $\Gamma_i$ to be the following colored graph. Let $e^\perp_1,\ldots,e^\perp_s$ be the edges in $B$ that goes from the attaching point of $G_i$ to the attaching point of $G_{i+1}$. Define $\Gamma_i$ to be the union of $G_i$, $G_{i+1}$, and $e^\perp_1,\ldots,e^\perp_s$ with color $2$; see Figure~\ref{tentacles_fig}. As in the previous case, it is easy to check that each $\Gamma_i$ is a barbell tree with one fewer bell. Let $v_i$ be the admissible coloring corresponding to $\Gamma_i$.
It is also clear from the construction that 
$$
2v=v_1+\cdots+v_k;
$$
i.e., the vector $v$ is positive $\Q$-linear combination of $v_1,\ldots,v_k$.\\

In each case, continuing inductively on each $\Gamma_{v_i}$, we obtain an equation expressing $v$ in terms of simple barbell trees with rational coefficients. Therefore, simple barbell trees generate $\La_\Q$.\\

To prove the other direction of the second statement in Theorem~\ref{indecomposable_th}, suppose $\Gamma_v$ is a simple barbell tree and $v=\sum_{i=1}^k \la_i v_i$, $\la_i\in \Q_{>0}$ for $i=1,\ldots,k$, where $\Gamma_{v_i}$ are simple barbell trees. We show that $v=v_i$ for all $i=1,\ldots,k$. By definition, each $\ov\Gamma_i$ is a chain (possibly empty). Therefore, each $\ov\Gamma_{v_i}$ is a subchain of $\ov\Gamma_v$. On the other hand, every barbell tree which is not a single bell has at least two leaf-bells. Therefore, for every $i=1,\ldots,k$, either $\Gamma_{v_i}=\Gamma_v$ or $\Gamma_{v_i}$ is one of the bells of $\Gamma_v$. By throwing away those $v_i$ such that $v=v_i$ and adjusting the coefficients $\la_i$ accordingly we can assume that each $\Gamma_{v_i}$ is one of the bells of $\Gamma_v$. We conclude that $\Gamma_v$ must be a single bell itself and $v=v_i$ for all $i=1,\ldots,k$.

 \begin{figure}
\begin{pspicture}(10,-1)(15,1.5)
\psset{unit=.8cm}

  \psarc(17,-1){1}{0}{90} \rput(17.5,-.5){2}
 \pscircle*(18,-1){0.1}\psline(18,-1)(19,-1)\rput(18.5,-.8){2}
  \pscircle*(17,0){0.1}\psline(17,0)(17,1)\rput(17.2,.5){2}
  
    \psarc(22,-1){1}{90}{180} \rput(21.5,-.5){2}
 \pscircle*(21,-1){0.1}\psline(21,-1)(20,-1)\rput(20.5,-.8){2}
  \pscircle*(22,0){0.1}\psline(22,0)(22,1)\rput(22.2,.5){2}
 
   \psarc(25,1){1}{180}{270} \rput(24.5,.5){2}

 \pscircle*(24,1){0.1}\psline(24,1)(23,1)\rput(23.5,1.2){2}
   \pscircle*(25,0){0.1}\psline(25,0)(25,-1)\rput(25.2,-.5){2} 
   
   \psarc(26,1){1}{270}{360} \rput(26.5,.5){2}

 \pscircle*(27,1){0.1}\psline(27,1)(28,1)\rput(27.5,1.2){2}
   \pscircle*(26,0){0.1}\psline(26,0)(26,-1)\rput(26.2,-.5){2}

\end{pspicture}
\caption{The graphs $\Gamma_i$ associated to Figure~\ref{Bell_fig}.}
\label{tentacles_fig}
\end{figure}
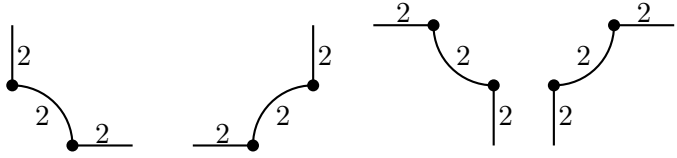

\end{proof}

\subsection{Flower triangulations (Proof of Theorem~\ref{Sphere_th})}\label{Bu_sec}
\begin{figure}
\begin{pspicture}(-4.5,-1.8)(5,1)
\psset{unit=.7cm}

\pscircle*(1.5,0){0.1} \rput(2,0){\small $p_6$}\psline(1.5,0)(0,0)
\pscircle*(.75,1.3){0.1} \rput(1,1.7){\small $p_1$}
\psline(.75,1.3)(0,0)\rput(.35,1){\tiny $\al_1$}\rput(.35,1.9){\tiny $\beta_1$}
\psline(.75,1.3)(0,0)\rput(.35,3.3){\tiny $\gamma$-curvs}
\pscircle*(-.75,1.3){0.1}\rput(-1,1.7){\small $p_2$}\psline(-.75,1.3)(0,0)
\pscircle*(-1.5,0){0.1}\rput(-2,0){\small $p_3$}\psline(-1.5,0)(0,0)
\pscircle*(-.75,-1.3){0.1}\rput(-1,-1.7){\small $p_4$}\psline(-.75,-1.3)(0,0)
\pscircle*(.75,-1.3){0.1}\rput(1,-1.7){\small $p_5$}\psline(.75,-1.3)(0,0)
\pscircle*(0,0){0.1} 

  \psplot[polarplot,algebraic,linewidth=.7pt,plotpoints=2000]{0}{TwoPi}{2.8*cos(3*x)*cos(3*x)}
 \psplot[polarplot=true,linearc=.2]{30}{150}{x 30 sub 1.5 mul sin 5 mul 3 min}
 \psplot[polarplot=true,linearc=.2]{150}{270}{x 150 sub 1.5 mul sin 5 mul 3 min}
 \psplot[polarplot=true,linearc=.2]{270}{390}{x 270 sub 1.5 mul sin 5 mul 3 min}

 \pscircle*(7,-1){0.1}\pscircle(7,-1.4){0.4}
\pscircle*(8,-1){0.1}\pscircle(8,-1.4){0.4}
\psline(7,-1)(7.5,0)(8,-1)

\pscircle*(9,-1){0.1}\pscircle(9,-1.4){0.4}
\pscircle*(10,-1){0.1}\pscircle(10,-1.4){0.4}
\psline(9,-1)(9.5,0)(10,-1)

\pscircle*(11,-1){0.1}\pscircle(11,-1.4){0.4}
\pscircle*(12,-1){0.1}\pscircle(12,-1.4){0.4}
\psline(11,-1)(11.5,0)(12,-1)

\psline(9.5,2)(7.5,0) \pscircle*(9.5,2){0.1} \pscircle*(7.5,0){0.1}
\psline(9.5,2)(9.5,0) \pscircle*(9.5,0){0.1}
\psline(9.5,2)(11.5,0) \pscircle*(11.5,0){0.1}

\end{pspicture}
\caption{Left- Flower triangulation of $(g,n)=(0,7)$. The un-labeled center puncture is $p_7$. Right- Dual graph of the triangulation on left.}
\label{Flower_fig}
\end{figure}
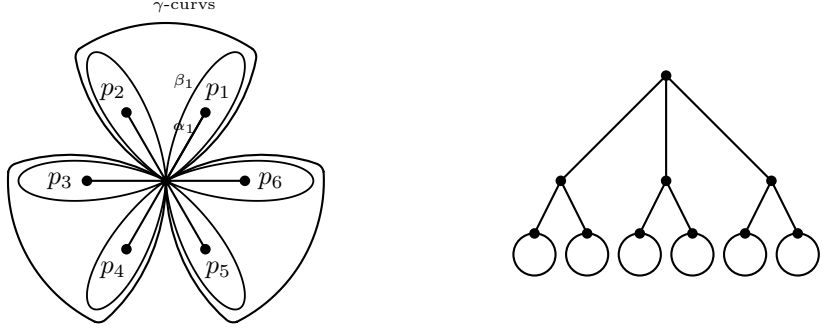

In this section, we discuss the example of flower triangulations shown in Figure~\ref{Flower_fig} and described below. 
This example has interesting properties making it useful for calculations. For a flower triangulation, we show that indecomposable curves/barbell trees in Theorem~\ref{indecomposable_th} are in one-to-one correspondence with ``complexity-zero" algebra generators of Bullock in \cite[Thm.~1]{Bu}. We also use this example and mutations to prove Theorem~\ref{Sphere_th}. \\

A \textbf{flower} triangulation is triangulation $\De$ that includes the following collection of arcs/edges. If $g=0$, we require $n\geq 4$. Let $\al_1,\ldots,\al_{n-1}$ be a set of mutually disjoint arcs connecting the $n$-th puncture point $p_n$ to $p_i$, for $i=1,\ldots,n-1$. To this collection we add a set of mutually disjoint arcs $\beta_1,\ldots,\beta_{n-1}$, where $\beta_i$ starts and ends at $p_n$ and is the boundary of a sufficiently small thickening $B_i$ of $\al_i$; see Figure~\ref{Flower_fig}-Left for genus $0$ and $n=7$. In other words, the only puncture included inside the disk $B_i$ is $p_i$ and the arc $\al_i$ lies inside $B_i$ (touching the boundary at $p_n$). For each $i=1,\ldots,n-1$, $\{\al_i,\beta_i,\alpha_i\}$ is a folded triangle (these will be the only folded triangles in $\De$). If $g=0$, the unfolding of the complement of $\bigcup_{i\in [n-1]} B_i$ is a polygon with $n-1$ edges $\beta_1,\ldots,\beta_{n-1}$ such that all the vertices are copies of $p_n$. 
In general, however, the complement of $\bigcup_{i\in [n-1]} B_i$ is a genus $g$ surface with $n-1$ boundary components and corners (where all the corners are copies of $p_n$).
We triangulate this complement by adding $6g+n-4$ more arcs $\gamma_1,\ldots,\gamma_{6g+n-4}$ that also start and end at (copies of) $p_n$. The dual graph $\Gamma$ of a genus zero flower is a trivalent  graph made of a tree $\cT$ with $n-1$ leaves and $n-1$ loops $\ell_1,\ldots,\ell_{n-1}$ attached to the leaves. The entire $\Gamma$ is a barbell tree whose bells are $\ell_1,\ldots,\ell_{n-1}$. These bells/loops correspond to the folded triangles $\{\al_i,\beta_i,\beta_i\}$, $i=1,\ldots,n-1$. The dual graph $\Gamma$ of Figure~\ref{Flower_fig}-Left is shown in Figure~\ref{Flower_fig}-Right. Since $\cT$ is a tree, barbell trees in $\Gamma$ are in one-to-one correspondence with  non-trivial subsets of $\{1,\ldots,n-1\}$. Given $\eset \neq B \subset \{1,\ldots,n-1\}$, the barbell tree corresponding to $B$ is the union of $\{\ell_i\}_{i\in B}$ and the minimal subtree in $\cT$ (colored with $2$) that contains the attaching leaves of these loops. Furthermore, it is straightforward to check that the set of generating curves $\cC(\De)=\{c_B\}_{\eset \neq B \subset \{1,\ldots,n-1\}}$ corresponding to these $2^{n-1}-1$ barbell trees matches the generating set of Bullock in \cite[Thm.~1]{Bu}. The latter corresponds to the handle decomposition given by the $\al$-curves.

\begin{remark}
In this example, simple barbell trees correspond to the significantly shorter list of ${n \choose 2}$ curves  $\{c_1,\ldots,c_{n-1}\}\cup \{c_{ij}\}_{i\neq j}$. They generate $\La_\Q$.
\end{remark}

 In the higher genus case, the argument generalizes and gives the 
 $2^{\tn{rank}\,H_1(\Si_{g,n},\Z)-1}-1$ generators of \cite[Thm.~1]{Bu} but the  tree $\cT$ will be replaced with a graph $G$ with $H_1(G)\cong \Z^{2g}$ (making it hard to draw). \\

{\bf Proof of Theorem~\ref{Sphere_th}.}
First, we explicitly show that Theorem~\ref{Sphere_th} holds for flower triangulations. Then, we use mutations and the identity \cite[(9.5)]{FST} to prove it for any arbitrary ideal triangulation.\\

Given a flower triangulation, we decompose every vector $v\in \Z^E$ as 
$$
v=(v^\al,v^\beta,v^\gamma), \qquad v^\al,v^\beta\in \Z^{n-1},\quad v^\gamma\in \Z^{6g+n-4},
$$
such that $v^\al$, $v^\be$, and $v^\gamma$ record the geometric intersection numbers with  $\al$, $\beta$, and $\gamma$-curves, respectively. For each $i=1,\ldots,n-1$, let $\vec{1}_i$ denote the $i$-th standard basis vector in $\Z^{n-1}$. 
If $a_i$ is the peripheral curve centered at $p_i$, it is clear that the corresponding admissible coloring vector $v(a_i)\in \Lambda$ is 
$$
\begin{cases}
v(a_i)=(\vec{1_i},\vec{0},\vec{0})& \tn{if}~~i=1,\ldots,n-1,\\
v(a_n)=(\vec{1}_{[n-1]},2\times \vec{1}_{[n-1]},2\times \vec{1}_{[6g+n-4]})& \tn{if}~~i=n.
\end{cases}
$$
Also, if $\mf{m}=\mf{m}(v)$ is a multi-curve that contains none of $a_1,\ldots,a_{n-1}$, then $v^\be=2v^\al$. Let $\La'$ denote the sub-monoid generated by $v(a_1),\ldots,v(a_{n-1})$, and $\La''$ denote the sub-monoid generated by the rest of the generators in $\cC(\De)$. By the calculations above, $v\in \La''$ if and only if $v^\beta=2v^\al$. Furthermore, the vector subspaces 
$$
\La'_\R=\La'\otimes \R,~\La''_\R=\La''\otimes \R\subset \R^{3(2g+n-2)}.
$$ 
are complimentary. We conclude that $\La\cong \La'\oplus \La''$.  Let $P''$ denote the polytope obtained by intersecting  $\La''_\R$ with the hyperplane $v\cdot \vec{1}=\tn{constant}\!>\!0$. Since, $P_\De^{\tn{rel}}$ is the union of those sub-polytopes of $P_\De$ that do not contain the points corresponding to $v(a_1),\ldots,v(a_n)$,  $P_\De^{\tn{rel}}$ is the union of those sub-polytopes of $P''$ that do not contain the point corresponding to $v(a_n)$. Since $v(a_1)+\ldots+v(a_n)=2\times \vec{1}_{E}\in \La$, we conclude that $v(a_n)$ is an interior vector of the cone $\La''_\R$. Therefore, $P_\De^{\tn{rel}}$ is simply the boundary of the  polytope  $P''$. We conclude that $P_\De^{\tn{rel}}$ is a sphere. \\

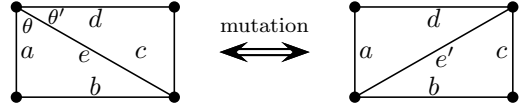
\begin{figure}
\begin{pspicture}(4,.3)(18,2.4)
\psset{unit=.3cm}

\psline[linewidth=.07](43,3)(43,7)(50,7)(50,3)(43,3)
\psline[linewidth=.07](43,3)(50,7)
\pscircle*(43,3){.25}
\pscircle*(43,7){.25}
\pscircle*(50,7){.25}
\pscircle*(50,3){.25}
\rput(43.5,5){\small $a$}
\rput(49.5,5){\small $c$}
\rput(46.5,3.5){\small $b$}
\rput(46.5,6.5){\small $d$}
\rput(47,4.8){\small $e'$}

\psline[linewidth=.07](28,3)(28,7)(35,7)(35,3)(28,3)
\psline[linewidth=.07](35,3)(28,7)
\pscircle*(28,3){.25}
\pscircle*(28,7){.25}
\pscircle*(35,7){.25}
\pscircle*(35,3){.25}
\psline[doubleline=true,doublesep=1.5pt]{<->}(37,5)(41,5)\rput(39,6.2){\scriptsize mutation}

\rput(28.5,5){\small $a$}\rput(28.5,6.2){\scriptsize $\theta$}\rput(29.8,6.6){\scriptsize $\theta'$}
\rput(33.5,5){\small $c$}
\rput(31.5,3.5){\small $b$}
\rput(31.5,6.5){\small $d$}
\rput(31,4.8){\small $e$}

\end{pspicture}
\caption{Mutation operation that changes one diagonal to another in a square.}
\label{Mutation-2}
\end{figure}

Every two triangulations can be connected by a sequence of mutations. A mutation is the operation of changing one diagonal $e$ to the other one $e'$ in a square; see Figure~\ref{Mutation-2}. With notation as in Figure~\ref{Mutation-2}, by \cite[(9.5)]{FST}, the colors $v_e$ and $v_{e'}$ of any admissible coloring vector $v$ are related by 
\bEq{mutation_eq}
v_{e'}=\tn{max}(v_a+v_c,v_b+v_d)-v_{e}.
\eEq

\noindent
Since we already know that Theorem~\ref{Sphere_th} holds for one triangulation, its is sufficient to prove the following lemma.
\begin{lemma}
Suppose $\De$ and $\De'$ are related by a mutation. Then $P_\De^{\tn{rel}}$ and $P_{\De'}^{\tn{rel}}$ are homeomorphic.
\end{lemma}

\begin{proof}
With notation as in  Figure~\ref{Mutation-2}, let $\wt{E}=E(\De)\cup \{e'\}= E(\De')\cup \{e\}$. By  (\ref{mutation_eq}), the image $\wt\La_\R$ of $\La_\R=\La_\R(\De)$ in $\R^{\wt{E}}$ is the graph of the piece-wise linear function (\ref{mutation_eq}).
Let 
$$
\La_\R^+=\{v\in \La_\R\colon v_{a}+v_{c}\geq v_{b}+v_{d
}\}\quad \tn{and}\quad \La_\R^-=\{v\in \La_\R\colon v_{a}+v_{c}\leq v_{b}+v_{d
}\}.$$
Then, $\wt\La_\R$ is a union of two cones $\wt\La_\R^\pm$ isomorphic to $\La_\R^{\pm}$ that intersect along the wall 
$$
\wt\La_\R^0\cong \La_\R^0=\{v\in \La_\R\colon v_{a}+v_{c}= v_{b}+v_{d
}\};
$$
see Figure~\ref{Break_fig}-Left.\newline
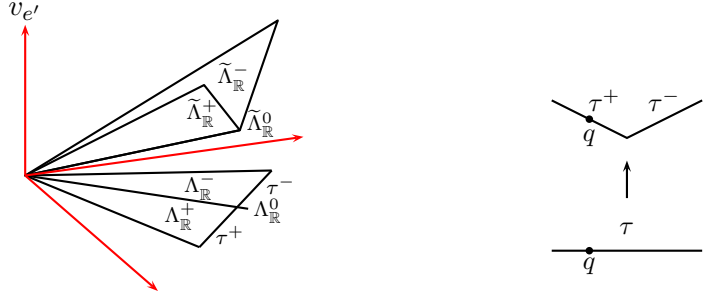
\begin{figure}
\begin{pspicture}(-3,-1.6)(3,2.5) 
\psline(7,-1)(9,-1)\rput(8,-.7){\small $\tau$}
\pscircle*(7.5,-1){0.05}\rput(7.5,-1.2){\small $q$}
\pscircle*(7.5,.75){0.05}\rput(7.5,.5){\small $q$}
\psline(7,1)(8,.5)(9,1)\rput(7.7,1){\small $\tau^+$}\rput(8.5,1){\small $\tau^-$}
\psline{->}(8,-.3)(8,0.2)
\psset{viewpoint=25 -10 10 ,Decran=30}
\psSolid[object=new,
  action=draw,
  sommets= 
  0 0 0 
  3 1 0 
 1 3 0 
 2.1 2.1 0 
 3 1 2 
 1 3 2 
 2 2 1
 ,
  faces={%
    [0 1]
    [0 2]
    [0 3]
    [1 2]
    [0 6 4]
    [0 6 5]
 }]%
\axesIIID[axisnames={,,v_{e'}},linecolor=red](4,4,2)%
\psPoint(3,1.9,0){U}
\uput[u](U){\scriptsize $\mathrm{\La_\R^0}$}
\psPoint(2.8,.85,0){U}
\uput[u](U){\scriptsize $\mathrm{\La_\R^+}$}
\psPoint(2,1.5,0){U}
\uput[u](U){\scriptsize $\mathrm{\La_\R^-}$}
\psPoint(2,2.3,.7){U}
\uput[u](U){\scriptsize $\mathrm{\wt\La_\R^0}$}
\psPoint(2.4,1.3,1.1){U}
\uput[u](U){\scriptsize $\mathrm{\wt\La_\R^+}$}
\psPoint(2.4,1.7,1.5){U}
\uput[u](U){\scriptsize $\mathrm{\wt\La_\R^-}$}
\psPoint(3.4,1.2,0){U}
\uput[u](U){\scriptsize $\mathrm{\tau^+}$}
\psPoint(2.2,2.5,0){U}
\uput[u](U){\scriptsize $\mathrm{\tau^-}$}

\end{pspicture}
\caption{Left: A schematic description of $\La_{\R}^\pm$, $\wt\La_{\R}^\pm$, $\La^0_\R$, $\wt\La_\R^0$, $\tau$, and $\tau^\pm$. Right: A hypothetical bad breaking. }\label{Break_fig}
\end{figure}

Homeomorphically, the polytope $P_\De$ is any slice of the cone $\La_\R$. At the level of polytopes, the process above corresponds to folding/breaking the polytope $P_\De$ along the hyperplane defined by the equation $v_{a}+v_{c}= v_{b}+v_{d
}$ into a union of two polytopes $P_\De^\pm$ intersecting along a codimension-$1$ face $P_\De^0$ ($P_\De^\pm$ are slices of $\La_\R^\pm\cong \wt\La_\R^\pm$.)\\

Recall that $P_\De^{\tn{rel}}$ is the union of those sub-polytopes in $P_\De$ whose closures do not contain the intersection points of the rays in $\La_\R$ corresponding to the peripheral curves. We define $P_\De^{\pm,\tn{rel}}$ similarly. A priori, for the following hypothetical reason, it is not clear that  
\bEq{Breaking-polytope_e}
P_\De^{+,\tn{rel}}\cup P_\De^{-,\tn{rel}}\cong_{C^0} P_\De^{\tn{rel}}
\eEq 
(the lefthand side could be larger). Suppose $\tau$ is a face of $P_\De$. For the sake of visualizing, assume $\tau$ is a $1$-dimensional polytope (an interval), as shown in Figure~\ref{Break_fig}-Right. Let $q\in \tau$ be a point corresponding to a peripheral curve. Thus, (the interior of) $\tau$ does not belong to $P_\De^{\tn{rel}}$. Let $\tau^\pm$ be the two segments of  $\tau$ ($\tau=\tau^+\cup \tau^-$) in $P_\De^{\pm,\tn{rel}}$ obtained by breaking/folding $\tau$ along the hyperplane $P_\De^0$. If $q\in \tau^+$ and $q\notin \tau^-$ as shown in Figure~\ref{Break_fig}-Right, then $\tau^+$ does not belong to $P_\De^{+,\tn{rel}}$ while $\tau^-$ does belong to $P_\De^{-,\tn{rel}}$. Therefore, (\ref{Breaking-polytope_e}) does not hold. The following argument shows that this does not happen and  (\ref{Breaking-polytope_e})  holds. \\

In terms of the corner variables $u_{\theta}$ and $u_{\theta'}$ defined in (\ref{Corner_e}) and shown in Figure~\ref{Mutation-2}-Left, the equality $v_{a}+v_{c}= v_{b}+v_{d
}$ is equivalent to $u_{\theta}=u_{\theta'}$. It is easy to see from the latter description that the admissible coloring vectors of peripheral curves belong to $\La_\R^0$.  Since the admissible coloring vectors of peripheral curves belong to $\La_\R^0$, we conclude that $P_\De^{\tn{rel}}$ is homeomorphic to $P_\De^{+,\tn{rel}}\cup_{P_\De^{0,\tn{rel}}}P_\De^{-,\tn{rel}}$; in other words, if a face $\tau$ in $P_\De$ contains a point corresponding to peripheral curves, both pieces $\tau^\pm$ of $\tau$ will contain that point (and their interiors must be excluded to get $P_\De^{\pm,\tn{rel}}$).  The same argument can be applied to $\De'$ (i.e., by writing $v_e$ as a function of $v_{e'}$ and $v_{a},\ldots,v_{d}$) to show that 
$$
P_\De^{\tn{rel}}\cong_{C^0} P_\De^{+,\tn{rel}}\cup P_\De^{-,\tn{rel}}=P_{\De'}^{+,\tn{rel}}\cup P_{\De'}^{-,\tn{rel}}\cong_{C^0}P_{\De'}^{\tn{rel}}.
$$
We conclude that $P_\De^{\tn{rel}}$ and $P_{\De'}^{\tn{rel}}$ are homeomorphic.
\end{proof}
\qed

\begin{remark}\label{Flower-g_rmk}
The calculations about show that flower triangulations are useful for calculating $\cD_{g,n}^\tn{rel}$ and $P_\De^\tn{rel}$, as it allows us to restrict to the smaller cone $\La''$ generated by $2^{2g+n-1}-(n+1)$  vectors.
Furthermore, since the $\beta$ component of the vectors in $\La''$ is twice the $\al$-component, we can drop the $\beta$-component and naturally identify $\La''$ with a sub-monoid of $\Z^{2n-5}$. We call these the reduced vectors and denote them by $\ov{v}$. We have $\tn{deg}(\ov{v})=3\deg(v^\al)+\deg(v^\gamma)$. In the genus zero case, to find $P^{\tn{rel}}_\De$, we only need the basic barbell trees that correspond to subsets $B$ with $|B|=2$. In other words, $P^{\tn{rel}}_\De$ is the boundary of a slice of the cone generated by the reduced vectors 
$$
\{\ov{v}_{ij}=\ov{v}(c_{ij})\}\subset \R^{2n-5}.
$$

\end{remark}
 
\subsection{Example of 5-punctured sphere}\label{n5example}
 
In this section, we use a flower triangulation of the five-punctured sphere and the corresponding $6$ generators of $\La''_\R$ to describe  $\cD_{5}^{\tn{rel}}$. This example has also been previously studied by Komyo and Simpson in \cite{K,S3}. Our approach provides a substantially sharper and shorter result.\\

The full monoid $\La$ is generated by $15$ elements. The full set of algebra relations among these $15$ generators can be obtained by reducing the full set of quantum algebra relations found in \cite[Appendix]{CL}, but we only study $\cD_{5}^{\tn{rel}}$ and identify it as a union of six\footnote{This number depends on the choice of triangulation.} toric pieces glued along their boundary to form a moment polytope complex homeomorphic to $S^3$. \\

By Remark~\ref{Flower-g_rmk}, in order to find $P^{\tn{rel}}_\De$, we need to understand the geometry of the reduced cone $\La''_\R$ generated by the reduced vectors $\{\ov{v}_B\}_{|B|=2}\in \Z_{\geq 0}^5$.  \begin{figure}
\begin{pspicture}(-4.5,-1.2)(5,1)
\psset{unit=.7cm}

\pscircle*(0,1){0.1} \rput(0,1.4){\small $p_1$}\psline(0,1)(0,0)
\pscircle*(1,0){0.1}\rput(1.4,0){\small $p_2$}\psline(1,0)(0,0)
\pscircle*(0,-1){0.1}\rput(0,-1.4){\small $p_3$}\psline(0,-1)(0,0)
\pscircle*(-1,0){0.1}\rput(-1.4,0){\small $p_4$}\psline(-1,0)(0,0)
\pscircle*(0,0){0.1} 

\psplot[polarplot,algebraic,linewidth=.7pt,plotpoints=2000]{0}{TwoPi}{2*cos(2*x)*cos(2*x)}
\psplot[polarplot=true,linearc=.3]{-45}{135}{x 45 add 1 mul sin 5 mul 2.5 min}

%
 \pscircle*(7,0){0.1}\pscircle(7,-0.4){0.4}
\pscircle*(8,0){0.1}\pscircle(8,-0.4){0.4}
\psline(7,0)(7.5,1)(8,0)

\pscircle*(9,0){0.1}\pscircle(9,-0.4){0.4}
\pscircle*(10,0){0.1}\pscircle(10,-0.4){0.4}
\psline(9,0)(9.5,1)(10,0)
\psline(7.5,1)(9.5,1) \pscircle*(7.5,1){0.1} \pscircle*(9.5,1){0.1}
\end{pspicture}
\caption{Left- Flower triangulation of $(g,n)=(0,5)$. The un-labeled center puncture is $p_5$. Right- Dual graph of the triangulation on left.}
\label{Flower5_fig}
\end{figure}
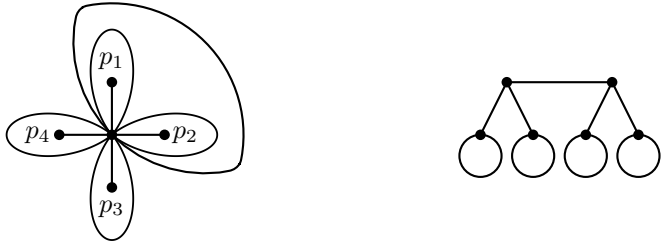
These reduced vectors are 
$$
\aligned
&\ov{v}_{12}=((1,1,0,0),0),~~\ov{v}_{34}=((0,0,1,1),0),\\
&\ov{v}_{13}=((1,0,1,0),2),~~\ov{v}_{24}=((0,1,0,1),2),\\
&\ov{v}_{14}=((1,0,0,1),2),~~\ov{v}_{23}=((0,1,1,0),2),\\
\endaligned
$$
satisfying one relation
$$
\ov{v}_{13}+\ov{v}_{24}=\ov{v}_{14}+\ov{v}_{23}.
$$

\begin{remark}\label{rmk225}
The integral generators of $\La''$ also include 
$$
\aligned
&\ov{v}_{123}=((1,1,1,0),2),~~\ov{v}_{124}=((1,1,0,1),2),\\
&\ov{v}_{134}=((1,0,1,1),2),~~\ov{v}_{234}=((0,1,1,1),2),\\
\endaligned
$$
satisfying  
\bEq{Redundantvs_e}
\aligned
& \ov{v}_{12}+\ov{v}_{23}+\ov{v}_{13}=2\ov{v}_{123},\\
& \ov{v}_{12}+\ov{v}_{24}+\ov{v}_{14}=2\ov{v}_{124},\\
& \ov{v}_{13}+\ov{v}_{34}+\ov{v}_{14}=2\ov{v}_{134},\\
& \ov{v}_{23}+\ov{v}_{34}+\ov{v}_{24}=2\ov{v}_{234}.\\
\endaligned
\eEq
\end{remark}
The original admissible coloring vectors $v_B$ are obtained by replacing $\ov{v}_B=(v^\al_B,v^\gamma_B)$ with $v_B=(v^\al_B,2v^\al_B,v^\gamma_B)$; e.g., $v_{13}=((1,0,1,0),(2,0,2,0),2)\in \La\subset \Z^9$.
The degrees of the corresponding curves/variables $c_B$ are determined by the sum of the coefficients of the full vectors $v_B$. Therefore, 
$$
\aligned
& \tn{deg}(c_{14})=\deg(c_{23})=\tn{deg}(c_{13})=\deg(c_{24})=8\\
&\tn{deg}(c_{12})=\deg(c_{34})=6,~~\tn{deg}(c_{ijk})=11
\endaligned
$$
\qed

\begin{remark}
Note that the map 
$$
\C\P[6,8,8,11]\supset (c_{12}\;c_{23}\;c_{13}= c_{123}^2)\lra [c_{12} \colon c_{23}\colon c_{13}]\in \C\P[6,8,8]
$$
is an isomorphism because 
$$
[c_{12}\colon c_{23} \colon c_{13} \colon c_{123}]=[(-1)^6c_{12}\colon (-1)^8c_{23} \colon(-1)^8 c_{13} \colon (-1)^{11} c_{123}]=[c_{12}\colon c_{23} \colon c_{13} \colon - c_{123}]
$$
in $\C\P[6,8,8,11]$. Similar phenomena happens for the other relations in~(\ref{Redundantvs_e}). This illustrates why there is no vertex corresponding to $\ov{v}_{ijk}$ in $P_\De^{\tn{rel}}$. \qed\end{remark}

We conclude that $P^{\tn{rel}}_\De$ consist of six $3$-dimensional polytopes whose vertices can be labeled as $q_{ij}$; $P^{\tn{rel}}_\De$ is the join of the square of vertices $q_{13},q_{14},q_{23},q_{24}$ (the vertices involved in the relation above Remark~\ref{rmk225}) and the line segment connecting $q_{12}$ to $q_{34}$. Geometrically, $q_{i_0j_0}$ is the point in the weighted projective space $\C\P[\deg(c_{ij})]$ where the $(i_0j_0)$-th variable $c_{i_0j_0}$ is $1$ and the rest are $0$. The six polytopes of $P^{\tn{rel}}_\De$ come in two types:
\bIt
\item $\ll q_{12},q_{34}, q_{13},q_{14}\rr \cong \ll q_{12},q_{34}, q_{14},q_{24}\rr\cong\ll q_{12},q_{34}, q_{24},q_{23}\rr\cong \ll q_{12},q_{34}, q_{23},q_{13}\rr\cong$ moment polytope of $\C\P[6\colon 6\colon 8 \colon 8]$;
\item $\ll q_{12},q_{13},q_{24},q_{14},q_{23}\rr\cong \ll q_{34},q_{13},q_{24},q_{14},q_{23}\rr \cong $ Pyramid over a square.
\eIt

\begin{figure}
\begin{pspicture}(-8,-2.5)(3,1) 
\psset{viewpoint=20 30 10 ,Decran=15}
\psSolid[
  object=new,
  action=draw,
  sommets= 
  0 0 0 
  0 6 0 
  6 0 0 
   6 6 0 
  3 3 4  
 3 3 -4  
 ,
  faces={%
    [2 3 4]
    [3 1 4]
    [0 1 4]
    [0 2 4]
     [2 3 5]
    [1 3 5]
    [1 0 5]
    [0 2 5]
 }]%
\psPoint(3,3,4){U}
\uput[u](U){\small $\mathrm{q_{12}}$}
\psPoint(3,3,-4){D}
\uput[d](D){\small$\mathrm{q_{34}}$}
\psPoint(0,6,0){R}
\uput[r](R){\small $\mathrm{q_{13}}$}
\psPoint(6,0,0){L}
\uput[l](L){\small $\mathrm{q_{24}}$}
\psPoint(0,0,0){B}
\uput[r](B){\small $\mathrm{q_{23}}$}
\psPoint(6,6,0){F}
\uput[l](F){\small $\mathrm{q_{14}}$}

\end{pspicture}
\caption{Union of the moment polytopes $\ll q_{12},q_{13},q_{24},q_{14},q_{23}\rr$ and  $\ll q_{34},q_{13},q_{24},q_{14},q_{23}\rr$.}
\label{MP}
\end{figure}
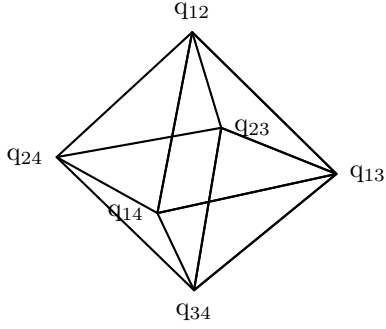

Figure~\ref{MP} shows the union of the last two pieces. The first four polytopes (each of which is a $3$-simplex) are attached to the outside of the union of these two along triangles, all together forming a topological manifold homeomorphic to $S^3$.

\section{Compactifications obtained from $\ov\SL(2,\C)$}\label{Kcomp-sec}

In the first two subsections below, we quickly review the construction of GIT quotients and then discuss the example which is used in Theorem~\ref{main_thm}. Readers familiar with the terminologies involved in GIT or not interested in the description of the spaces $X_{\mf{a}}$ in Theorem~\ref{main_thm} may skip to Section~\ref{RC}.

 \subsection{A quick tour of GIT}\label{GITtour}

Given a complex projective variety $Y$, an algebraic action of a complex reductive Lie group $G$ on $Y$, and a $G$-invariant ample line bundle $\cL$ with a lift of the $G$-action (called {\bf $G$-linearization}), the GIT quotient $Y\sslash_{\cL} G$ is the projective variety
$$
X=\tn{proj} \Big(R(Y)^{G}=\bigoplus_{k\geq 0} H^0(Y,\cL^{\otimes k})^{G} \Big).
$$
We say $\cL_1$ and $\cL_2$ are equivalent if $\cL_1^{\otimes k_1} = \cL_2^{\otimes k_2}$ as $G$-linearized line bundles. The GIT quotient  $X=Y\sslash_{\cL} G$ only depends on the equivalence class of $\cL$. In what follows, whenever the choice of $\cL$ is clear or fixed in the discussion, we will simply write  $Y\sslash G$ instead of $Y\sslash_{\cL} G$.\\

Geometrically, the \textbf{semistable} locus $Y^{ss}\subset Y$ is defined to be the complement of the base locus of the sub-linear system $|\cL^{\otimes m}|^{G}$ for some sufficiently large $m$ and $X=Y^{ss}\sslash G$ is a good quotient of that in the sense of \cite[Dfn.~2.3.6]{Hos}. Here, $|\cL^{\otimes m}|^{G}$ means the sub-linear system of G-invariant sections. \\

Restricted to $Y^{ss}$, we have $\cL^{\otimes m}=\varphi^*\cO_{\P^N}(1)$ where 
\bEq{varphimap_e}
\varphi\colon Y^{ss}\lra X \subset \C\P^N=\P\Big(H^0(Y,\cL^{\otimes m})^{G} \Big)
\eEq
is the quotient map realizing $Y\sslash G$ as a projective variety. We will denote the \textbf{unstable} locus, i.e., the complement of $Y^{ss}$ in $Y$, by $Y^{us}$. \\

The \textbf{stable} locus $Y^s\subset Y$ consists of points $y\in Y$ with finite stabilizer $G_y$ such that $G\cdot y$ is a closed orbit in some affine chart $Y_f=Y\!-\!(f^{-1}(0))$ for some nonzero section $f\in H^0(Y,\cL^{\otimes k})^{G}$. Equivalently, $y$ is stable iff it is semistable, the orbit $G\cdot y$ in $Y^{ss}$ is closed, and $G_y$ is finite. Restricted to $Y^{s}$, $\varphi$ is an orbifold $G$-bundle and $\varphi(Y^s)$ is a geometric quotient $Y^s/G$. The latter is a possibly empty Zariski open subset of $X$. More generally, a semistable point $y$ is said to be \textbf{polystable} if its orbit is closed in $Y^{ss}$. Let $Y^{ps}\subset  Y^{ss}$ denote the subset of polystable points. Note that every stable point is polystable. We say two semistable points are $S$-equivalent if their orbit closures meet in $Y^{ss}$.
The orbit closure of every semistable point $x$ contains a unique polystable orbit. Moreover, if $y$ is semistable but not stable, then this unique polystable orbit is also not stable. Finally, if $y,y'\in Y^{ss}$, then $\varphi(y) = \varphi(y')$ if and only if $y$ and $y'$ are $S$-equivalent. Consequently, there is a bijection of sets $X \cong Y^{ps}/G$.\\

If $Z$ is a $G$-invariant subvariety of $Y$, then $\cL|_{Z}$ is ample and $G$-linearized, and $Z^{ss}=Y^{ss}\cap Z$; therefore, 
$$
Z\sslash_{\cL|_{Z}} G \subset Y\sslash_{\cL} G
$$ 
is also a subvariety. \\

Finding semistable points, as described above, is hard. The Hilbert-Mumford criterion described below provides an effective way of detecting unstable points.  Consider a $1$-parameter subgroup (or 1-PS) 
$$
\vr\colon \C^*\lra G.
$$
For every $y\in Y$, the map
$$
\vr_y\colon \C^*\lra Y, \quad t \lra \vr(t)\cdot y
$$
extends (because $Y$ is compact) to a rational curve $\ov\vr_y\colon \P^1\lra Y$ (Note that $\vr_y$ is the constant map iff $\tn{im}(\vr)\subset G_y$). Let 
$$
y_0=\lim_{t\lra 0}\vr_y(t)\in Y \quad \tn{and}\quad y_\infty=\lim_{t\lra \infty}\vr_y(t) \in Y.
$$
Then $y_0$ and $y_\infty$ are fixed by the $\C^*$-action of $\vr$ and we get linear actions of $\C^*$ on $\cL_{y_0}$ and $\cL_{y_\infty}$ with weights $w_0$ and $w_\infty$.
Define $w(y,\vr)=w_0\in \Z$ (thus, $w(y,\vr^{-1})=w_\infty)$. {\it By the Hilbert-Mumford criterion, $y$ is semistable iff $w(y,\vr)\geq 0$ and stable iff $w(y,\vr)>0$ for all $1$-PS $\vr$}. Since $w(gy, g\vr g^{-1})=w(y,\vr)$, it is enough to check one $\vr$ in each conjugacy class.\\

In the case of $G=\SL(2,\C)$, there is only one conjugacy class of non-divisible $1$-parameter subgroups
\bEq{CTorus_e}
\vr\colon \C^*\lra \SL(2,\C),\quad \vr(t)=\begin{pmatrix}
t &0  \\
0 &t^{-1}
\end{pmatrix}.
\eEq
Therefore, the function $w(y)\defeq \tn{min}(w(y,\vr),w(y,\vr^{-1}))$ gives a stratification of Y, 
$$
Y=\coprod_{d\in \Z} Y_{d}, \quad Y_{d}=w^{-1}(d),
$$
such that 
$$
Y^s=\coprod_{d\in \Z_+} Y_{d}, \quad Y^{ss}=\coprod_{d\in \Z_{\geq 0}} Y_{d}, \quad \tn{and} \quad Y^{us}=\coprod_{d\in \Z_-} Y_{d}.
$$
Define $S_{d}\subset Y_{d}$ to be the subset of points $y$ such that the limit $z=\lim_{t\lra 0}(\vr(t)\cdot y)$ belongs to a component $Z_d$ in the intersection of the fixed locus of the $\C^*$-action of $\vr$ on $Y$ and $w^{-1}(d)$; c.f. Hesselink's theorem \cite[Thm.~5.36]{Hos}. Note that $Z_d$ depends on the particular choice of $\vr$ in its conjugacy class where as the weight function $w$ does not. 
We have $Y_{d}=\SL(2,\C)\cdot S_{d}$. Therefore, in order to find $Y^{us}$, 
\bIt
\item first, we find the unstable fixed points of the $\C^*$-action of $\vr$ in (\ref{CTorus_e}), 
\item then, we find the points converging to those fixed points, 
\item and, finally, we consider the $G$-orbits of the last step.
\eIt

\subsection{GIT quotients of products of $\P^1$}\label{P1quotient_sec}
The group $G=\SL(2,\C)$ acts on $\P^1$ by Mobius tansformations
\bEq{Mact_e}
\begin{pmatrix}
a & b\\
c& d
\end{pmatrix}
\cdot [x_1\colon x_2] = [ax_1+b x_2\colon cx_1+dx_2].
\eEq
For the tautological line bundle $\cO_{\P^1}(-1)\lra \P^1$, we consider the $\SL(2,\C)$-linearization induced by the canonical embedding 
$$
\cO_{\P^1}(-1)\subset \P^1\times \C^2.
$$
In fact, this is the only possible $\SL(2,\C)$-linearization of $\cO_{\P^1}(-1)$ because the center of $\SL(2,\C)$ has dimension $0$.
For any other $\cO_{\P^1}(a)$, we consider the linearization induced by the isomorphisms $\cO_{\P^1}(1)\cong \cO_{\P^1}(-1)^*=\cO_{\P^1}(-1)^{\otimes -1}$ and  $\cO_{\P^1}(a)=\cO_{\P^1}(-1)^{\otimes -a}$.\\

For  $m\geq 3$ and $\mf{a}=(a_1,\ldots,a_m)\in \Z^m$, let $\cL_\mf{a}\defeq \cO(a_1,\ldots,a_m)\lra (\P^1)^m$ denote the line bundle 
\bEq{BoxTimes}
\cO(a_1)\boxtimes\cdots \boxtimes \cO(a_m)
\defeq q_1^*\cO_{\P^1}(a_1)\otimes\cdots \otimes q_m^*\cO_{\P^1}(a_m)
\eEq
where $q_i\colon (\P^1)^m\lra \P^1$ is the projection to the $i$-th factor.
In the rest of this section we consider the diagonal action of $\SL(2,\C)$ on $Y=(\P^1)^m$
equipped with the very ample line bundle $\cL_{\mf{a}}$, with $a_i>0$ for all $1\leq i\leq m$. This example has been studied in numerous contexts in the past; see \cite{MFK,Th}. We will relate it to $\SL(2,\C)$ character varieties in the next section. \\

Corresponding to every point $(p_1,\ldots,p_m)\in (\P^1)^m$ we get a partition 
\bEq{partition_e}
I=I(p)\colon \Big(\{1,2,\ldots,m\}= I_1 \sqcup \cdots \sqcup I_k\Big)
\eEq
such that $i,j \in I_\ell$ if and only if $p_i=p_j$. Let $(\P^1)^m_I$ denote the set of points of type $I=\{I_1,\ldots,I_k\}$. 
For two partitions 
$$
I=\{I_1,\ldots,I_k\} \quad \tn{and}\quad J=\{J_1,\ldots,J_\ell\}, 
$$
we write $J\leq I$ if each $J_j$ is the union of a sub-collection of $I$. This defines a partial order on the set $\cP(m)$ of all such partitions with the unique maximal element 
$$
I_{\tn{max}}=\{\{1\},\ldots,\{m\}\}.
$$ 
The following lemma easily follows from the three-bullet process described after (\ref{CTorus_e}). For each $\mf{a}\in \Z_+^m$ and $I\!\in\! \cP(m)$ as in (\ref{partition_e}), define $\mf{a}_I \in \Z_+^{k}$ to be the tuple obtained by replacing each sub-tuple $(a_j)_{j\in I_\ell}$ with $\sum_{j\in I_\ell} a_j$ (arranged in some random order). For instance, if 
$$
m=4,\quad \mf{a}=(1,1,1,1),\quad \tn{and}~I=\{\{1,2\},\{3\}, \{4\}\},
$$ 
then $\mf{a}_{I}=(2,1,1)$ or some permutation of that.

\bLm{ssCondition_lmm}
With notation as in (\ref{partition_e}), the $\cL_{\mf{a}}$-unstable locus $Y^{us}$ of $Y=(\P^1)^m$ is the union of components $(\P^1)^m_{I}$ satisfying
\bEq{usI}
\sum_{i\in I_\ell} a_i > \sum_{i\notin I_\ell} a_i 
\eEq
for some $1\leq \ell \leq k$. In other words, $(p_1,\ldots,p_m)\in Y^{ss}$ iff
\bEq{ssI}
\sum_{i\in I_\ell} a_i \leq \sum_{i\notin I_\ell} a_i \quad \forall~1\leq \ell \leq  k,
\eEq
and $(p_1,\ldots,p_m)\in Y^{s}$ iff
\bEq{sI}
\sum_{i\in I_\ell} a_i < \sum_{i\notin I_\ell} a_i \quad \forall~1\leq \ell \leq k.
\eEq
\eLm

For a given $\mf{a}$, we say $I\!\in\! \cP(m)$ is unstable if (\ref{usI}) holds, is semistable if (\ref{ssI}) holds, is stable if (\ref{sI}) holds, and is strictly semi-stable if there is $\ell$ such that 
$$
\sum_{i\in I_\ell} a_i = \sum_{i\notin I_\ell} a_i.
$$
We denote the corresponding subsets of $\cP(m)$ by $\cP^{us}_\mf{a}(m)$, $\cP_\mf{a}^{ss}(m)$, $\cP_\mf{a}^s(m)$, and $\cP_\mf{a}^{ss-s}(m)$, respectively.\\

For $Y=(\P^1)^m$, $G=\SL(2,\C)$, and $\mf{a}\in \Z_+^m$ define
$$
X_{\mf{a}}\defeq Y\sslash_{\cL_{\mf{a}}} G = (\P^1)^{m}\sslash_{\cL_a} \SL(2,\C).
$$ 

\begin{example}
If one of $a_i$ is larger than the sum of all other $a_j$ then $Y^{ss}=\eset$ and $X_{\mf{a}}$ will be defined to be the empty space.\qed
\end{example}
\begin{example}\label{sss_ex}
If $a_{i_0}$ is equal the sum of all other $a_j$, for some $i_0$, then $Y^{s}=\eset$ and 
$$
Y^{ss}=\{(p_1,\ldots,p_m)\in (\P^1)^m\colon p_i\neq p_{i_0}~ \forall i\neq i_0\}.
$$ 
Therefore, we can fix $p_{i_0}$ to be $\infty$ to show that
$$
X_{\mf{a}}\cong \C^{m-1}\sslash \tn{Aff}_{\C}
$$
where $ \tn{Aff}_{\C}$ is the group of affine transformations $z\lra az+b$ acting component-wise on $\C^{m-1}$. Since there are no bounded non-constant polynomials on $\C^{m-1}$, the quotient space $X_{\mf{a}}$ is just a point. In this example, the only polystable orbit is the orbit  corresponding to the $2$-set partition 
$$
I=\{i_0\} \sqcup \{1,2,\ldots, i_0-1,i_0+1,\ldots,m\}.
$$
\qed
\end{example}

In the light of the two (bad) examples above and consistently with \cite[Page~11]{Ka}, we call $\mf{a}\in \Z_+^m$ \textbf{non-degenerate} if every $a_i$ is smaller than the sum of all other $a_j$; i.e., if $\cP^{s}_\mf{a}(m)\neq \eset$. If $\mf{a}$ is non-degenerate, then 
\bEq{Stable-a_ex}
(\P^1)^{m}_*\defeq \{(p_1,\ldots, p_m)\in (\P^1)^m\colon p_i\neq p_j \quad \forall~i\neq j\}\subset Y^s
\eEq
is a dense open subset; therefore, the un-compactified genus zero $m$-marked Deligne-Mumford space 
$$
\cM_{0,m}\cong X_{\mf{a}}^\circ\defeq (\P^1)^{m}_*/\SL(2,\C)
$$ 
is a dense open subset of $X_{\mf{a}}$. In other words, the set of all GIT quotients $X_{\mf{a}}$, when running through all non-degenerate $\mf{a}$, are (possibly different) compactifications of the open dense geometric quotient $\cM_{0,m}$. More generally, for every partition $I=\{I_1,\ldots,I_k\}\in \cP^{s}_\mf{a}(m)$, the image
$$
X_{\mf{a}_I}^\circ \cong (\P^1)^{m}_I\sslash \SL(2,\C)\subset X_{\mf{a}}
$$
of $(\P^1)^{m}_I$ in $X_{\mf{a}}$ is isomorphic to $\cM_{0,k}$. The closure  $X_{\mf{a}_I}$ of $X_{\mf{a}_I}^\circ$ has boundary 
$$
\partial X_{\mf{a}_I}= X_{\mf{a}_I}-X^{\circ}_{\mf{a}_I}= \bigcup_{J<I} X_{\mf{a}_J}=\bigcup_{J<I} X_{\mf{a}_J}^\circ.
$$
If $I\!\in\! \cP^{ss-s}_\mf{a}(m)$, or $I\!\in\! \cP^{s}_\mf{a}(m)$ and $k=3$, then $X_{\mf{a}_I}^\circ=X_{\mf{a}_I}$ is just a point. By Example~\ref{sss_ex}, in the first case, this is the point corresponding to a unique polystable orbit in $(\P^1)^m_I$.
Putting together, we conclude that
$$
X_{\mf{a}}-(Y^s/\SL(2,\C))
$$
is a finite set of points corresponding to polystable orbits, and each polystable orbit corresponds to a partition of $\{1,\ldots,m\}$ into two sets,
$$
\{1,\ldots,m\}=I_1 \sqcup I_2
$$
such that 
\bEq{ED_e}
\sum_{i\in I_1} a_i=\sum_{i\in I_2} a_i.
\eEq

In \cite{Th}, Thaddeus explains the nature of birational transformations or contractions corresponding to changing $\mf{a}$. \\

In the example of diagonal M\"obius action on $Y=(\P^1)^m$ explained above, one can try other techniques for constructing a nice quotient space such as ``Hilbert" and ``Chow" quotients. The latter does not need fixing a polarization and will be denoted by $Y\sslash_{\tn{Ch}} G$ below. Briefly, $Y\sslash_{\tn{Ch}} G$ is the closure in Chow variety of a geometric quotient $U/G$ for some nice open set $U\subset Y$ (usually smaller than $Y^s$). A point in $Y\sslash_{\tn{Ch}} G$ is a $G$-invariant cycle $Z=\sum c_i Z_i$ in $Y$ where each component is the closure of a single orbit. 
 In \cite{Ka}, Kapranov studies  Chow quotients of Grassmannians vs. products of $\P^1$ and proves the following results; see Theorems 0.4.3,~2.2.4,~2.4.7, and Proposition~4.3.17 in \cite{Ka}. See also \cite{GJM} for some more recent results.\\

Given $\mf{a}\in \Z_+^m$, the complex $(m-1)$-torus 
$$
(\C^*)^m_{\bullet}=\{(t_1,\ldots,t_m)\in (\C^*)^m\colon \prod_{i=1}^m t_i=1\}
$$
acts by $(x_1,\ldots,x_m)\lra (t_1^{a_1}x_1,\ldots,t_m^{a_m}x_m)$ on $\C^m$ and thus on the Grassmannian $\tn{Gr}(2,m)$. The line bundle $\cO_{\tn{Gr}(2,m)}(1)$ can be linearized compatibly with every such action. 

\begin{theorem}(\cite[Theorems~(4.1.8), (2,4,7), (0.4.3), Proposition~(4.3.7)]{Ka})\label{Kapranov} 
With notation as above, we have:
\bEn
\item the Chow quotients $(\P^1)^m\sslash_{\tn{Ch}} \SL(2,C)$ and $\tn{Gr}(2,m)\sslash_{\tn{Ch}} (\C^*)^m_{\bullet}$  are isomorphic to the smooth stable (Deligne-Mumford) compactification $\ov\cM_{0,m}$;
\item for every $\mf{a}\in \Z_+^m$, the GIT quotients $X_{\mf{a}}$ and $\tn{Gr}(2,m)\sslash_{\cO(1)_{\mf{a}}} (\C^*)^m_{\bullet}$ are isomorphic;
\item\label{ChtoGIT} for every GIT quotient $Y\sslash_{\cL} G$, if $Y^{s}$ is non-empty, there is a regular birational morphism $f\colon Y\sslash_{\tn{Ch}} G\lra Y\sslash_{\cL} G$;
\item consequently, if $\mf{a}\in \Z_+^m$ is non-degenerate, then there is a regular birational morphism 
$$
f_{\mf{a}}\colon \ov\cM_{0,m}\lra X_{\mf{a}}
$$
that restricts to the identity map on $\cM_{0,m}\subset \ov\cM_{0,m},X_{\mf{a}}$;
\item furthermore, $\ov\cM_{0,m}$ is the inverse limit of $\{X_{\mf{a}}\}_{\mf{a}\in \Z_+^m}$.
\eEn
\end{theorem}  

The map $f$ in \ref{ChtoGIT} has the following form. Given a $G$-invariant cycle\footnote{More precisely, in \cite[(0.4.8)]{Ka}, $Z$ is taken to be a limit position of closures of single orbits.} $Z=\sum c_i Z_i$ in $Y\sslash_{\tn{Ch}} G$, there is at least one $\cL$-semistable orbit in the support of $Z$, and every two such orbits are $S$-equivalent. The map $f$ which takes $Z$ to the point in $Y\sslash_{\cL} G$ represented by any of these equivalent orbits is a morphism of algebraic varieties. 
In $\tn{Gr}(2,m)\sslash_{\tn{Ch}} (\C^*)^m_{\bullet}$, every $Z_i$ is a toric variety which is the closure of some $(n-1)$-dimensional orbit and $c_i=1$; see \cite[Prop~1.2.11, 1.2.15]{Ka}. Similarly, in $(\P^1)^m\sslash_{\tn{Ch}} \SL(2,C)$ every $Z_i$ is the closure of some $3$-dimensional orbit and $c_i=1$.\\

In this paper, we will need the case where $m$ is even and $\mf{a}$ only involves $m/2$ independent integers in the following sense. In what follows, we assume $m$ is even.
\bDf{Syma_dfn}
Suppose $m$ is even, we say $\mf{a}\in \Z_+^m$ is \textbf{symmetric} if
\bEq{evena_e}
a_{2i-1}=a_{2i}=b_i \quad \forall~i=1,\ldots,m/2.
\eEq
\eDf

\begin{remark}
It is clear that every symmetric $\mf{a}$ is non-degenerate, therefore $Y^s\neq \eset$ and $X_{\mf{a}}$ is complex projective variety of the expected dimension.
\end{remark}

Let 
\bEq{De_curve}
\De_{\P^1}\subset \P^1\times \P^1
\eEq
denote the diagonal subspace, 
$$
\De_i \defeq \lrc{(p_1,\ldots, p_{m})\in (\P^1)^{m}\colon p_{2i-1}=p_{2i} }= (\P^1)^{\{1,\ldots, 2i-2\}}\times \De_{\P^1} \times (\P^1)^{\{2i+1,\ldots, m\}}\subset (\P^1)^m,
$$
and 
\bEq{De-divisor}
\De=\De_1\cup \cdots \cup \De_{m/2}.
\eEq
The $\SL(2,\C)$-invariant subspace $H^0(\cO(1,1))^{\SL(2,\C)}$ of the space of sections of the line bundle $\cO(1,1)$ on 
$$
\P^1[x_1\colon x_2]\times \P^1[y_1\colon y_2]
$$
is generated by the section $\mf{s}=x_1y_2-x_2y_1$. 
For $1\leq i\leq m/2$, let 
\bEq{mfsi}
\mfs_i\in H^0\big(q_{2i-1}^*\cO(1)\otimes q_{2i}^*\cO(1)\big)^{\SL(2,\C)}
\eEq
denote the corresponding pullback section of the pullback line bundle 
\bEq{Li}
\cL_i\defeq q_{2i-1}^*\cO(1)\otimes q_{2i}^*\cO(1)\lra (\P^1)^m.
\eEq
With notation as in (\ref{evena_e}), we have $\cL_{\mf{a}}=\cL_1^{\otimes b_1} \otimes \cdots \otimes \cL_{m/2}^{\otimes b_{m/2}}$ and  $\De_i= \mfs_i^{-1}(0)$. Let 
$$
D_{\mf{a},i}= \Big(\De_i \sslash_{\cL_{\mf{a}_i}} \SL(2,\C)\Big)\subset X_{\mf{a}},
$$
where $\cL_{\mf{a}_i}\defeq \cL_{\mf{a}}|_{\De_i}$ 
is the line bundle $\cL_{\mf{a}_i}\lra \De_i\cong (\P^1)^{m-1}$ given by 
$$
\mf{a}_i=(b_1,b_1,\ldots, b_{i-1},b_{i-1},2b_{i}, b_{i+1},b_{i+1},\ldots,b_{m/2},b_{m/2})\in \Z_+^{m-1}.
$$ 
Every individual divisor\footnote{For some $i$, $D_{\mf{a},i}$ could be empty. This happens when $b_i$ is larger than the sum of other $b_j$.} $D_{\mf{a},i}$ might not be a $\Q$-Cartier but the weighted $\sum_{i} b_i D_{\mf{a},i}$ is $\Q$-Cartier  and ample, see Remark~\ref{Descent-L} below.\\

Let 
\bEq{si_inv}
\si\colon (\P^1)^2 \lra (\P^1)^2
\eEq
denote the involution exchanging two copies of $\P^1$; i.e. $\si(p_1,p_2)=(p_2,p_1)$. For each $i=1,\ldots,m/2$, let 
$$
\si_i\colon (\P^1)^m \lra (\P^1)^m,
$$ 
denote the involution acting on $(p_{2i-1},p_{2i})$ by $\si$ and on the rest of the points by the identity map.
If $\mf{a}$ is symmetric, each $\si_i$ induces a similarly denoted involution 
$$
\si_i \colon X_{\mf{a}}\lra X_{\mf{a}}
$$
with $D_{\mf{a},i}\subset \tn{Fix}(\si_i)$. Furthermore, these involutions commute, and by Theorem~\ref{main_thm}, they appear as monodromy maps of fibers of $\ov\cR_{n}$ around certain hypersurfaces in $\C^n$.

\bLm{toric_lmm}
Suppose $\mf{a}$ is symmetric and, with notation as in (\ref{evena_e}),  $b_{i_0}\!>\! \sum_{j\neq i_0} b_j$ for some (unique) $1\leq i_0 \leq m/2$. Then $X_{\mf{a}}$ is a toric variety (of the expected dimension) and $D_{\mf{a},i_0}=\eset$. 
\eLm

\bPf
We may assume $i_0\!=\!1$.
By Lemma~\ref{ssCondition_lmm} and (\ref{Stable-a_ex}), $Y^s$ is non-empty and every point $(p_1,\ldots,p_m)\in Y^{ss}$ satisfies $p_1\!\neq\! p_2$. In other words, $Y^{ss}\!\subset\! Y-\De_1$, which implies the last statement of the lemma. Therefore, in each orbit, we may take a representative with $p_1=0$ and $p_2=\infty$. The subgroup of $\SL(2,\C)$ fixing $0$ and $\infty$ is the 1-parameter $\C^*$-subgroup  in (\ref{CTorus_e}). We conclude that
\bEq{Red-to-C*}
Y\sslash_{\cL_{\mf{a}}} \SL(2,\C) = (\P^1)^{m-2}\sslash_{\cL_{\mf{a}'}} \C^*
\eEq
where $\mf{a'}=(a_3,\ldots,a_m)$. The righthand side of (\ref{Red-to-C*}) is a toric variety with the moment polytope
$$
\lrc{(x_3,\ldots,x_m)\in \R^{m-2}\colon  -a_i \leq x_i \leq a_i, ~\sum_{i=3}^m x_i =0}.
$$
\ePf
\begin{example}
For $m=6$, the symmetric tuple $\mf{a}=(3,3,1,1,1,1)$ satisfies the assumption of Lemma~\ref{toric_lmm} and 
$$
X_{\mf{a}}=\tn{proj} \Big(\frac{\C[x_0,x_1,\ldots,x_5]}{(x_0x_1=x_2x_3=x_4x_5)}\Big)
$$
with the torus action 
$$
(t_0,t_1,t_2)\cdot (x_0,x_1,\ldots,x_5)\lra (t_0x_0,t_0^{-1}x_1,t_1x_2,t_1^{-1}x_3,t_3x_4,t_3^{-1}x_5), \qquad \forall(t_0,t_1,t_2)\in (\C^*)^3.
$$
This is a singular toric variety with $8$ codimension-$1$ boundary faces; see Figure~\ref{MomenPolytop}. \qed
\begin{figure}
\begin{pspicture}(-8,-2)(3,1.2) 
\psset{viewpoint=20 10 10 ,Decran=30}
\psSolid[object=new,
  action=draw,
  sommets= 
  1 1 -1 
  1 -1 1 
  -1 1 1 
   -1 1 -1 
  -1 -1  1 
 1 -1 -1 
 ,
  faces={%
    [1 0 2]
    [4 3 5]
    [0 1 5]
    [1 2 4]
    [2 0 3]
 }]%
\end{pspicture}
\caption{The polytope obtained as the convex hull of $\pm (-1,1,1), \pm (1,-1,1), \pm (1,1,-1)$ in $\R^3$.}\label{MomenPolytop}
\end{figure}
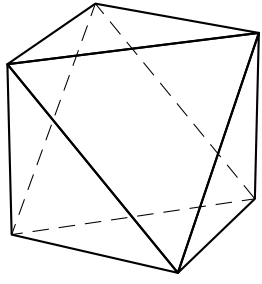
\end{example}

\begin{remark}\label{Descent-L}
As the referee pointed out to us, an individual divisor $D_{\mf{a},i}$ is often not $\Q$-Cartier. In \cite[Proposition~4.2]{KKV}, the authors find necessary and sufficient conditions for when a $G$-linearized line bundle (and thus the corresponding Cartier divisors) descends to a (GIT) quotient.
More precisely, \cite[Definition 4.1]{KKV} applies to the quotient map $Y^{ss}\lra X$ in (\ref{varphimap_e}). In our example, $Y=(\P^1)^m$ and there are two situations to consider. If no $b_{i_0}$ is larger than the sum of the other $b_i$s (as it happens in Lemma~\ref{toric_lmm}), then the unstable locus $Y^{us}$ has complex codimension higher than one. Therefore, 
$$
\tn{Pic}_{G}(Y^{ss})=\tn{Pic}_G\big((\P^1)^m\big),
$$
where $\tn{Pic}_{G}$ is the group of isomorphism classes of $G$-linearized line bundles, and $G=\SL(2,\C)$. Here, the connected
group $\SL(2,\C)$ acts trivially on the discrete group $\tn{Pic}\big((\P^1)^m\big)$, so any line
bundle is left invariant by the $G$-action and the lift of action is unique because the center of $G$ has dimension $0$. Consequently, we have
$$
\tn{Pic}_G((\P^1)^m)= \tn{Pic}((\P^1)^m)=\Z^m.
$$
Therefore, each $G$-linearized line bundle is the restriction of $\cL_{\mf{c}}$ in (\ref{BoxTimes}) to $Y^{ss}$ for some $\mf{c}\in \Z^m$. By \cite[Proposition~4.2]{KKV}, $\cL_{\mf{c}}$ descends to $X_{\mf{a}}$ if and only if the condition (PB) below holds:\\

(PB) For every $x\in Y^{ss}$ whose $G$-orbit is closed, the isotropy group $G_x$ acts trivially on the fiber $\cL_{\mf{c}}|_{x}$ of $\cL_{\mf{c}}$.\\

We have that closed orbits correspond to stable and poly-stable points. In fact, if there are at least three distinct point in $(\P^1)^{m}_I$, then the isotropy group $G_x$ is trivial and there is nothing to check. The only potentially problematic points are those $x$ in the pre-image of the finite set 
$$
X_{\mf{a}}-(Y^s/\SL(2,\C))
$$
of polystable orbits, each of which correponds to a partition of $\{1,\ldots,m\}$ into two sets,
$$
\{1,\ldots,m\}=I_1 \sqcup I_2
$$
such that (\ref{ED_e}) holds.
For such $x$, the isotropy group is $\C^*$ and the weight of the action of $\C^*$ on $\cL_{\mf{c}}|_{x}$ is 
$$
\pm \Big( \sum_{i\in I_1} c_i - \sum_{i\in I_2} c_i \Big).
$$
In conclusion, (PB) is satisfied at $x$ if and only if (\ref{ED_e}) implies 
$$
\sum_{i\in I_1} c_i - \sum_{i\in I_2} c_i =0.
$$
The divisor $D_{\mf{a},i}$ corresponds to $\mf{c}=(0,0,\ldots,1_{2i-1},1_{2i},\ldots, 0,...,0)$. Therefore, $D_{\mf{a},i}$ is $\Q$-Cartier at 
$\ov{x}\in X_{\mf{a}}-(Y^s/\SL(2,\C))$ if and only if the partition 
$$
\{1,\ldots,m\}=I_1 \sqcup I_2
$$
corresponding to $\ov{x}$ does not not satisfy $2i-1,2i\in I_1$ or $2i-1,2i\in I_2$.  Also, for instance, by (\ref{ED_e}) the divisor $\sum_i b_i D_{\mf{a},i}$ is Q-Cartier.\\

Similarly, in the exceptional case of Lemma~\ref{toric_lmm}, the complement of $Y^{ss}$ has complex codimension one; it consists of $\De_{i_0}$ and some complex codimension two subset. Then, one can follow a modification of the argument above or use the fact that the quotient is toric to study the Cartier property of the divisors $D_{\mf{a},i}$ at polystable orbits.\qed
\end{remark}

\bLm{autoss}
If $m$ is even and $\mf{a}$ is symmetric, with $Y=(\P^1)^m$ as before, we have $Y-\De \subset Y^{ss}$ and the affine variety
$$
U_{\mf{a}}\defeq X_{\mf{a}}-D_{\mf{a}}\cong \Big(Y-\De\Big)\sslash_{\cL_\mf{a}} \SL(2,\C)
$$
is a Zariski dense open subset of $X_{\mf{a}}$ independent of the choice of $\mf{a}$. In other words, if $\mf{a}$ and $\mf{a}'$ are symmetric, the birational transformation $X_{\mf{a}}\! \dashrightarrow \!X_{\mf{a'}}$ in \cite[Proposition~1.1]{Th}  restricts to an isomorphism of the affine subvarieties $U_{\mf{a}}$ and $U_{\mf{a}'}$.   
\eLm

\bPf 
For $(p_1,\ldots,p_{m})\in (\P^1)^{m}$, let 
$$
\{1,2,\ldots,m\}= I_1 \sqcup \cdots \sqcup I_k
$$
be the partition in (\ref{partition_e}). If $p\notin \De$, for every $1\leq \ell \leq k$ and $1\leq i \leq m/2$, only one of $2i\!-\!1$ or $2i$ belongs to $I_\ell$. Therefore, 
$$
\sum_{j\in I_\ell} a_j \leq \sum_{j\notin I_\ell} a_j \quad \forall~1\leq \ell \leq  k.
$$
It follows from Lemma~\ref{ssCondition_lmm} that $p$ is semistable. The second statement follows from the fact that $D_{\mf{a}}$ supports an ample $\Q$-Cartier divisor. More precisely, we have 
$$
\mfs_{\mf{a}}\defeq \prod_{i=1}^{m/2} \mfs_i^{b_i} \in H^0(Y,\cL_{\mf{a}})^{\SL(2,\C)} 
$$
and 
$$
U_{\mf{a}}=\tn{Spec}(R_{\mf{a}}), \quad R_{\mf{a}}=\bigoplus_{k\geq 0} \frac{H^0(Y,\cL_{\mf{a}}^{\otimes k})^{\SL(2,\C)}}{\mfs_\mf{a}^{k}}.
$$
For symmetric $\mf{a},\mf{a'}\in \Z_+^m$, the algebra isomorphism
$$
\frac{H^0(Y,\cL_{\mf{a}}^{\otimes k})^{\SL(2,\C)}}{\mfs_\mf{a}^{k}}\lra 
\frac{H^0(Y,\cL_{\mf{a'}}^{\otimes k})^{\SL(2,\C)}}{\mfs_\mf{a'}^{k}}, \quad f\lra \frac{\mfs_\mf{a'}^{k}}{\mfs_\mf{a}^{k}} f\qquad \forall~k\geq 0,
$$
gives the desired isomorphism of the affine subvarieties $U_{\mf{a}}$ and $U_{\mf{a'}}$.
\ePf

The morphisms $f_{\mf{a}}\colon \ov\cM_{0,m} \lra X_{\mf{a}}$ in Theorem~\ref{Kapranov} may/will blowup certain points/loci in the Zariski open set  
$$
U_{\mf{a}}=\big(Y-\De\big)\sslash_{\cL_{\mf{a}}} \SL(2,\C)\subset X_{\mf{a}}.
$$ 
Therefore, $U_{\mf{a}}$ does not embed as an open set in $\ov\cM_{0,m}$; however, by Lemma~\ref{autoss}, the inverse limit  of $\{X_{\mf{a}}\}$ with $\mf{a}$ only ranging over the set of symmetric tuples includes the affine variety $U_{\mf{a}}$ as a dense open set and is more likely to be smooth than any of $X_{\mf{a}}$ alone. This naturally raises the following question.

\begin{question}\label{InvLimXa}
If $m$ is even, what is the inverse limit $\ov\cM_{0,m}^{\tn{sym}}$ of the collection $\{X_{\mf{a}}\}$ with $\mf{a}$ only ranging over the set of symmetric tuples $\mf{a}\in \Z_+^m$?
\end{question}

\subsection{Relative compactification of $\cR_n$}\label{RC}

In \cite{K}, a natural compactification of $\tn{SL}(2,\C)$ to a projective quadratic threefold is used to construct a compactification of any fiber of $\cR_n \lra \C^n$. We show that this can be done (for various choices of an ample line bundle and) uniformly across the entire family. Subsequently, Proposition~\ref{Main_Prop} enables us to study a subfamily of $\cR_n$ as a family of divisors in the GIT quotients $X_{\mf{a}}$ constructed above (or their inverse limit). Proposition~\ref{Main_Prop} also simplifies the discussion of stability condition in \cite[Sec.~4]{K}.\\

\textbf{Notation.} In the following, we will use the matrix notation 
$
\begin{pmatrix}
*
\end{pmatrix}
$
for ordinary matrices and 
$
\begin{bmatrix}
*
\end{bmatrix}
$
for matrices defined up to scalar multiplication in a projective space (i.e., for points in some $\P\tn{GL}(N)$).

The smooth projective quadratic threefold 
$$
\ov{\SL(2,\C)}\defeq \lrc{ \begin{bmatrix}
a & b & \\
c & d & \\
&& e
\end{bmatrix}\colon ad-bc =e^2 }\subset \P^4[a:b:c:d:e]\subset \P(M_{3}(\C))
$$
includes $\SL(2,\C)\cong (e\neq 0)$ as a Zariski open subset. 
The pencil of quadratic surfaces
$$
|Q_t|_{t\in \P^1}=\lrc{ \begin{bmatrix}
a & b & \\
c & d & \\
&& e
\end{bmatrix}\colon ad-bc =e^2,~ a+d =t e }
$$
on $\ov{\SL(2,\C)}$ has the base locus 
$$
 C=\lrc{ \begin{bmatrix}
a & b & \\
c & d & \\
&& 0
\end{bmatrix}\colon ad-bc =0,~a+d =0 }\subset \ov{\SL(2,\C)} 
$$
and the fiber over $t=\infty$ is the smooth quadratic surface 
$$
Q_\infty=\lrc{ \begin{bmatrix}
a & b & \\
c & d & \\
&& 0
\end{bmatrix}\colon ad-bc =0 }= \ov{\SL(2,\C)} \cap (e=0).
$$
We primarily care about $t\in \C\subset \P^1$ but the discussion below extends to the entire parameter space $\P^1$ and is useful for studying the degenerations of $\cR_{n,\mf{t}}$.  \\

For $t \neq \pm 2$, $Q_t$ is a smooth quadratic surface; therefore, it is isomorphic to $\P^1\times \P^1$; see~(\ref{IsoP1P1}).  For $t\!=\!\pm 2$, $Q_t$ has an isolated A1-singularity (i.e. $x^2=yz$)  at 
\bEq{sing-point}
\begin{bmatrix}
 I_2 & \\
&\pm 1
\end{bmatrix} 
\eEq
and it is a cone over $C$; see the discussion after the proof of Proposition~\ref{Main_Prop}. 
The conjugation action of $\SL(2,\C)$ on itself extends to an action on $\ov{\SL(2,\C)}$ defined by 
\bEq{InnerAction}
\begin{bmatrix}
A & \\
& e
\end{bmatrix}\lra  \begin{bmatrix}
PAP^{-1}& \\
& e
\end{bmatrix}.
\eEq
Note the $P=\pm I_2$ acts as identity; therefore, the $\SL(2,\C)$-action in (\ref{InnerAction}) can be exchanged by the action of 
$$
\PGL(2,\C)\cong \SL(2,\C)/\pm I_2.
$$

For $\mf{b}=(b_1,\ldots,b_{n-1})\in \Z^{n-1}$, let $\cL_\mf{b}\defeq \cO(b_1,\ldots,b_{n-1})\lra \ov\SL(2,\C)^{n-1}$ denote the line bundle 
$$
\cO(b_1)\boxtimes \cdots \boxtimes \cO(b_{n-1})
$$  
where $\cO(b)$ is the restriction to $\ov\SL(2,\C)$ of the line bundle $\cO(b)$ on $\P^4$ and $\boxtimes$ is defined similarly to (\ref{BoxTimes}).\\

For $g=0$ and $n>0$, the fundamental group $\pi_1(\Si_{0,n})$ is freely generated by the peripheral curves around the first $n-1$ punctures. Therefore, 
$$
\cR_n=\SL(2,\C)^{n-1}\sslash\SL(2,\C)
$$ 
and a point in $\cR_{n}$ is the equivalence class $[A_1, \ldots , A_{n-1}]$ of an $(n-1)$-tuple of matrices\footnote{ in the semi-stable locus of the action of $\SL(2,\C).$} 
such that $A_n=(A_1 \ldots A_{n-1})^{-1}$ is the matrix associated to $a_n$.
Since $\tn{tr}(A)=\tn{tr}(A^{-1})$ for all $A\in \SL(2,\C)$, the projection map $\Pi_{n}\colon  \cR_{n} \lra~\C^n$ is then
\bEq{pin0_e}
[A_1, \ldots , A_{n-1}]\lra  \mft=(t_1,\ldots,t_n)\\
=\Big(\tn{tr}(A_1),\ldots,\tn{tr}(A_{n-1}),\tn{tr}(A_n)=\tn{tr}\big(A_1 \ldots A_{n-1}\big)\Big).
\eEq

By Proposition~\ref{Main_Prop} below and Lemma~\ref{autoss} (or simply the fact that $\cR_{n}=\SL(2,\C)^{n-1}\sslash \SL(2,\C)$ is an affine GIT quotient), any GIT quotient 
\bEq{Whole}
\ov\SL(2,\C)^{n-1}\sslash_{\cL_\mf{b}} \SL(2,\C)
\eEq
is a compactification\footnote{This is the compactification in \cite{BLR}.} $\ov\cR_n$ of $\cR_n$; however, the fibration map~(\ref{pin0_e}) does not extend to $\ov\cR_n$ (i.e., (\ref{Whole}) is not a relative compactification of $\cR_n$). Nevertheless, 
for each fixed $\mft=(t_1,\ldots,t_n) \in \C^{n}$, restricted to 
$$
Q_{\mf{t}}\defeq Q_{t_1}\times \cdots \times Q_{t_{n-1}}\subset \ov\SL(2,\C)^{n-1},
$$
the expression $F_{\mft}$ defined by 
\bEq{section-equation_e}
F_{\mft}\lrp{
\begin{bmatrix}
A_1 & \\
& e_1
\end{bmatrix}
,\ldots, \begin{bmatrix}
A_{n-1} & \\
& e_{n-1}
\end{bmatrix}}\defeq 
 \tn{tr}(A_1\ldots A_{n-1})- t_n e_1\ldots e_{n-1}
\eEq
is an $\SL(2,\C)$-invariant section of $\cL_{(1,\ldots,1)}|_{Q_{\mft}}$, and 
\bEq{ovcRK}
 \ov\cR_{n,\mft}\defeq  \wt{\cR}_{n,\mft} \sslash_{\cL_{\mf{b}}|_{Q_{\mft}}} \SL(2,\C),\quad \wt{\cR}_{n,\mft}\defeq (F_{\mft}^{-1}(0))\subset Q_{\mf{t}},
\eEq
defines a compactification $\ov\cR_{n,\mft}$ of $\cR_{n,\mft}$ (depending on $\mf{b}$) that sits as a divisor in $Q_{\mf{t}} \sslash \SL(2,\C)$.
For $\mf{b}=(1,\ldots,1)$, this is the compactification appearing \cite{K}. For $n=5$, it is further blown up to make its boundary divisor normal crossing. The boundary divisor $\cD_{n,\mft}$ of $\ov\cR_{n,\mft}$ is the intersection of $\ov\cR_{n,\mft}$ with the image of the $\SL(2,\C)$-invariant divisor 
\bEq{partialQ}
\partial Q_{\mft}\defeq \big(C\times Q_{t_2}\times \cdots \times Q_{t_{n-1}}\big) \cup \cdots \cup \big(Q_{t_1}\times \cdots \times Q_{t_{n-2}}\times C\big) 
\eEq
in $Q_{\mf{t}} \sslash \SL(2,\C)$.
\\

To obtain a relative compactification of $\cR_n\lra \C^n$, we need to resolve the base locus of the pencil $|Q_t|_{t\in \P^1}$. Let 
\bEq{pi_blowup}
\pi\colon M=\tn{Bl}_{C} \ov\SL(2,\C)\lra \ov\SL(2,\C),
\eEq
denote the blowup of $\ov\SL(2,\C)$ along the base curve $C$ with the exceptional divisor $E=C\times \P^1$. 
Obviously, $M$ admits a projection 
\bEq{tr-on-M}
\tn{tr}\colon M\lra \P^1
\eEq
for which $M_t=\tn{tr}^{-1}(t)\cong Q_t$ is a compactification of $\SL(2,\C)_{\tn{tr}=t}$ with the boundary divisor $M_t \cap E= C\times \{t\}\cong C$. Let $M_{\C}\defeq M-M_\infty=\tn{tr}^{-1}(\C)$ and note that $M-(E\cup M_\infty)\cong \SL(2,\C)$. 

\bLm{action-on-M}
The conjugation action in (\ref{InnerAction}) lifts to an action on $M$ that preserves each fiber $M_t$ (i.e., the projection map $M\lra \P^1$ is $\SL(2,\C)$-invariant).
\eLm 
\bPf
The action (\ref{InnerAction}) preserves $C\cong \P^1$ and acts by M\"obius transformations on it. By the first definition in \cite[Page 163]{HS}, we have $$
M=\tn{Proj} \Big(\bigoplus_{d\geq 0} \mc{I}^d\Big),
$$
where $\cI$ is the ideal sheaf of $C\subset \ov\SL(2,\C)$.
Since $C$ is $\SL(2,\C)$-invariant, $\SL(2,\C)$ acts on $\cI$ and thus on $M$.
Furthermore, the projection map $M\lra \ov\SL(2,\C)$ is $\SL(2,\C)$-equivariant.

Since the action (\ref{InnerAction}) preserves each $Q_t$, the extended action on $M$ preserves the proper-transform $M_t$ of $Q_t$, for all $t\in \P^1$. Note that the induced  action of $\SL(2,\C)$ on the normal bundle 
$$
\cN_{\ov\SL(2,\C)}C \cong \cO_{C}(2)\oplus \cO_{C}(2)
$$
is diagonal and that the restriction of the extended action on $M$ to the exceptional divisor
$$
E=\P(\cN_{\ov\SL(2,\C)}C)\cong \P(\cO_{C}(2)\oplus \cO_{C}(2))\cong \P(\cO_{C}\oplus \cO_{C})=C\times \P^1
$$
is trivial on the second factor.
\ePf

We consider the product $M^{n-1}\times \P^1$ with the projection map
$$
\tn{Tr}\defeq (\tn{tr}_1,\ldots,\tn{tr}_{n-1},\tn{tr}_n)\colon M^{n-1}\times \P^1\lra (\P^1)^{n},
$$
where $\tn{tr}_i$, with $1\leq i \leq n-1$, is the projection map (\ref{tr-on-M}) on the $i$-th copy of $M$ and $\tn{tr}_n$ is simply projection to the last component. The action of $\SL(2,\C)$ on $M^{n-1}\times \P^1$ by the diagonal action on $M^{n-1}$ and trivial action on $\P^1$ makes the projection map $\tn{Tr}$ $\SL(2,\C)$-invariant. The line bundle 
\bEq{L'b}
\cL'_{\mf{b}}=\Big(\pi^*\cO(b_1)\boxtimes \cdots \boxtimes \pi^*\cO(b_{n-1})\Big) \otimes \bigotimes_{i=1}^n \tn{tr}_i^*\cO(1)\lra M^{n-1}\times \P^1
\eEq
is ample and admits a natural linearization acting trivially on $\bigotimes_{i=1}^n \tn{tr}_i^*\cO(1)$. Here, $\pi^*\cO(b_i)$ is the pullback of $\cO(b_i)$ on $\ov\SL(2,\C)$ to $i$-th copy of $M$ via $\pi$ in (\ref{pi_blowup}). \\

With $F_{\mf{t}}$ as in (\ref{section-equation_e}), the subset 
$$
M^{n-1}\times \P^1
\supset \wt\cR^{\tn{rel}}_n = \lrc{ (x_1,\ldots,x_{n-1},t_n) \colon F_{\mft}(\pi(x_1),\ldots,\pi(x_{n-1}))=0, \quad \mft=\tn{Tr}(x_1,\ldots,x_{n-1},t_n) }  
$$
is an $\SL(2,\C)$-invariant divisor in $M^{n-1}\times \P^1$.  In fact, the following lemma follows from the explicit equations (\ref{Apqt}) and (\ref{limitiota_e}) in the next section and some straightforward comparison of line bundles.

\begin{lemma}\label{Fmft_lmm}
The expression $F_{\mft}(\pi(x_1),\ldots,\pi(x_{n-1}))$ is an $\SL(2,\C)$-invariant section of the line bundle $\cL'_{(1,\ldots,1)}\lra M^{n-1}\times \P^1$.
\end{lemma}

By definition, the intersection of $\wt\cR_n^{\tn{rel}}$ with $\tn{Tr}^{-1}(\mft)$ is isomorphic to $\wt{\cR}_{n,\mft}$ in (\ref{ovcRK}). We conclude that
\bIt
\item the compact GIT quotient 
$$
\ov{\mf{M}}_n\defeq (M^{n-1}\times \P^1)\sslash_{\cL'_{\mf{b}}} \SL(2,\C)
$$
includes $\wt\cR^{\tn{rel}}_n\sslash \SL(2,\C)$
as a divisor, 
\item the projection map $\tn{Tr}$ descends to a projection map 
$$
\tn{Tr}\colon \ov{\mf{M}}_n\lra (\P^1)^n,
$$ 
\item and, the restriction 
$$
\Pi_n=\tn{Tr}|_{\ov\cR^{\tn{rel}}_n} \colon \ov\cR^{\tn{rel}}_n\lra \C^n, \qquad \ov\cR^{\tn{rel}}_n\defeq \big(\wt\cR_n^{\tn{rel}}\sslash \SL(2,\C)\big)\cap \tn{Tr}^{-1}(\C^n),
$$
is a quasi projective relative compactification of $\cR_n$ (depending on $\mf{b}$) sitting as a divisor in 
$$
\mf{M}_n\defeq \tn{Tr}^{-1}(\C^n)=(M_\C^{n-1}\times \C)\sslash_{\cL'_{\mf{b}}} \SL(2,\C).
$$
\eIt

This construction defines the relative compactification $\ov\cR^{\tn{rel}}_n\subset \mf{M}_n$ stated in Theorem~\ref{main_thm} and Remark~\ref{mfM_rmk} depending on the parameters $\mf{b}\in \Z_+^{n-1}$. By lifting to the inverse limit of 
$$
\lrc{(M_\C^{n-1}\times \C)\sslash_{\cL'_{\mf{b}}}}_{\mf{b}\in \Z_+^{n-1}}
$$
we can get a more refined relative compactification of $\cR_n$.
In the next section, at the cost of restricting to $\C^n_\circ=(\C-\{\pm 2\})^{n-1}\times \C$, we find a uniform and explicit parametrization of $\mf{M}_n$ and complete the proof of Theorem~\ref{main_thm}.

\subsection{Uniform parametrization of conjugacy classes}\label{uniform_s}

In what follows, we provide an explicit parametrization of $Q_{\mf{t}}$ to uniformly compactify (almost every)  $\cR_{n,\mft}$  as linearly equivalent divisors in any of the GIT quotients $X_{\mf{a}}$ defined in Section~\ref{P1quotient_sec}. These calculations also describe the monodromy of the family $\ov\cR^{\tn{rel}}_n$ around the hypersurfaces $t_i=\pm 2$ in the base space $\C^n$.  \\

For $s= \sqrt{t^2-4}$, we have
$$
\aligned
Q_t=&\lrc{ \begin{bmatrix}
a & b & \\
c & d & \\
&& e
\end{bmatrix}\colon  \Big(\frac{a+d}{2}\Big)^2-\Big(\frac{a-d}{2}\Big)^2- bc = e^2,~ a+d=t e }\\
=&\lrc{ [a-d,b,c,e] \in \P^3 \colon  \Big(\frac{(a-d)+ s e}{2}\Big) \Big(\frac{(a-d)- s e}{2}\Big)  +bc=0 }.
\endaligned
$$
By definition, $t$ determines $s$ up to a sign factor. Therefore, to be precise, we take $(t,s)$ to be a point in the curve $\cA\cong \C^*$ defined by $(t^2-s^2=4)\subset \C^2$.
The projection map 
\bEq{pcA_e}
\cA\lra \C, \quad (t,s)\lra t,
\eEq
is a two-fold branched covering with the symmetry group 
\bEq{siA_e}
\si\colon \cA\lra \cA, \quad \si(t,s)=(t,-s),
\eEq
and the branching loci $(t,s)=(\pm 2, 0)$. The completion $\ov\cA\cong \P^1$ of $\cA$ in $\P^2$ is a smooth conic and includes two extra points at the infinity denoted by $\pm \infty =[0,0,\pm 1]\in\ov\cA\subset \P^2$. The twofold branched covering map (\ref{pcA_e}) extends to $\ov\cA\lra \P^1$
and is un-ramified at $t=\infty$.\\

For $t\neq \pm 2$ (or equivalently $s \neq 0$), the map
\bEq{IsoP1P1}
\aligned
&\iota_{(s,t)}\colon \P^1\times \P^1\lra Q_t \subset \ov\SL(2,\C),\\ 
&\iota_{(s,t)}(p,q)=\begin{bmatrix}
A(p,q,(t,s)) & \\
& e(p,q)
\end{bmatrix}\qquad \forall (p,q)=\big([x_1\colon x_2], [y_1\colon y_2]\big)\in \P^1\times \P^1 
\endaligned
\eEq
such that 
\bEq{Apqt}
A(p,q,(t,s))=
 \begin{pmatrix}
 \beta_2 x_2y_1-\beta_1 x_1 y_2 & sx_1y_1 \\
-sx_2y_2 & \beta_1 x_2y_1-\beta_2 x_1y_2 
\end{pmatrix},\quad 
e(p,q)=(x_2y_1-x_1y_2),
\eEq
$$
\beta_1=\frac{1}{2}(t+s) \quad \tn{and}\quad \beta_2 =\frac{1}{2}(t -s),
$$
identifies $\P^1\times \P^1$ and $Q_t$. 
Since
$$
(e=0) \Rightarrow  (x_2y_1-x_1y_2=0)  \Rightarrow  \frac{x_1}{x_2}= \frac{y_1}{y_2},
$$
the pre-image of $C\subset Q_t$ in $\P^1\times \P^1$ is the diagonal  $(1,1)$-curve $\De_{\P^1}$ in (\ref{De_curve}). 

\begin{remark}\label{Qinfty}
The case $t=\pm \infty$ corresponds to the points $\pm \infty=[0,0,\pm 1]\in\ov\cA\subset \P^2$  and 
$$
\begin{bmatrix}
A(p,q, \pm \infty) &\\
& e
\end{bmatrix}=
\begin{bmatrix}
- x_1 y_2 & x_1y_1 & \\
-x_2y_2 &  x_2y_1 &\\
&&0
\end{bmatrix}\in Q_\infty.
$$
extends the parametrization to $t=\pm \infty$.
\end{remark}

\begin{remark}\label{change-sign_rmk}
Changing the sign of $s$ in (\ref{IsoP1P1}) corresponds to the involution
$$
[x_1\colon x_2]\times [y_1\colon y_2]\stackrel{\si}{\lra} [y_1\colon y_2]\times [x_1\colon x_2]
$$
on $\P^1\times \P^1$ with the fixed point set $\De_{\P^1}$; see (\ref{si_inv}). In other words, the map 
\bEq{iotacA_e}
\iota \colon \P^1\times \P^1 \times \big(\ov\cA-\{(\pm 2,0)\}\big) \lra \ov\SL(2,\C), \quad (p,q,(t,s))\lra \iota_{(t,s)}(p,q),
\eEq
is $\si$-invariant, where the action of $\si$ on the $\cA$-factor is by (\ref{siA_e}).
\end{remark}
\begin{remark}
Note that 
$$
\beta_1/\beta_2=\frac{t +s}{t -s} = h \in  \C^*-\{ 1\}
$$
such that $h$ is well-defined up to inversion (because $s$ is well-defined up to sign) and 
\bEq{t-notation}
t=  (h^{\frac{1}{2}}+h^{-\frac{1}{2}}), \quad s=  (h^{\frac{1}{2}}-h^{-\frac{1}{2}}), \quad \beta_1= h^{\frac{1}{2}}, \quad \beta_2=h^{-\frac{1}{2}}.
\eEq
Thus, we can write $A$ as a function of $h$ by
$$
A(p,q,h)=\begin{pmatrix}
h^{-\frac{1}{2}} x_2y_1-h^{\frac{1}{2}} x_1 y_2 & (h^{\frac{1}{2}}-h^{-\frac{1}{2}}) x_1y_1 \\
-(h^{\frac{1}{2}}-h^{-\frac{1}{2}}) x_2y_2 & h^{\frac{1}{2}} x_2y_1-h^{-\frac{1}{2}} x_1y_2 
\end{pmatrix}.
$$
The matrix associated to $0\times \infty = [0\colon 1]\times [1\colon 0]$ is 
$$
\begin{bmatrix}
h^{-\frac{1}{2}} & 0 & \\
0 & h^{\frac{1}{2}} & \\
&& 1
\end{bmatrix}
$$
which corresponds to the $\C^*$-action $z=\frac{x_1}{x_2} \lra h^{-1}z$. \end{remark}

The following proposition is the key observation of this section.

\bPr{Main_Prop}
For every $(p,q,(t,s))$ and every $\rho \in \SL(2,\C)$, we have 
$$
A(\rho(p),\rho(q),(t,s))= \rho^{-1} A(p,q,(t,s))\rho
$$
In other words, pullback to $\P^1\times\P^1$ of the conjugation action of $\SL(2,\C)$ on $Q_t$ by the parametrization map $\iota_{(s,t)}$ in (\ref{IsoP1P1}) \underline{is independent of} $(s,t)\in \cA$, and coincides with the diagonal M\"obius action on $\P^1\times\P^1$.
\ePr
\bPf
With notation as in Remark~\ref{change-sign_rmk}, for $(p,q)\in (\P^1)^2-\De_{\P^1}$, 
$$
\frac{1}{e(p,q)}A(p,q,h)\in \SL(2,\C)
$$ 
is the unique M\"obius transformation that fixes $p$ and $q$ and acts by the $\C^*$-action $z \lra h^{-1}z$ on $\P^1-\{p,q\}\cong \P^1-\{0,\infty\}\cong \C^*$. 

\ePf

\begin{remark}\label{Qpm2_rmk}
There are a priori two ways to desingularize/parametrize $Q_{\pm2}$ and extend the discussion above.
The blowup of $Q_{\pm 2}$ at its only singular point yields the ruled surface $\P(\cO_{C}\oplus \cO_{C}(2))$ that include a $(-2)$-curve in place of the singular point (\ref{sing-point}). Alternatively, $Q_{\pm 2}$ can be obtained as a quotient of $\P^1\times \P^1$  by $\Z_2\times\Z_2$ in the following way. First, $Q_{\pm 2}$ is the quotient with respect to the $\Z_2$-action of
\bEq{P^2inv_e}
[u,v,w]\lra [u,v,-w]
\eEq
on $\P^2$. The twofold branched covering map $\P^2\lra Q_{\pm 2}$ realizing the latter can be written as 
$$
[u, v, w] \lra 
\begin{bmatrix}
 w^2+uv & u^2 & \\
-v^2 &  w^2-uv& \\
&& \pm w^2
\end{bmatrix} 
$$
On the other hand $\P^2$ is the $\Z_2$-quotient of $\P^1\times\P^1$ by the involution $\si$ in (\ref{si_inv}). The two-fold covering map $\P^1\times \P^1\lra \P^2$ realizing the latter can be written as 
$$
[x_1\colon x_2]\times [y_1\colon y_2] \lra [u=x_1y_1\colon v=x_2y_2\colon w=x_1y_2+x_2y_1].
$$
 Under this map, the $\Z_2$-action (\ref{P^2inv_e}) lifts to the $\Z_ 2$ action 
 $$
  [x_1\colon x_2]\times [y_1\colon y_2]\stackrel{\de}{\lra} [-x_1\colon x_2]\times [-y_1\colon y_2].
 $$
 Therefore, the $\Z_2\times\Z_2$ quotient of $\P^1\times\P^1$ generated by $\si$ and $\de$
yields $Q_{\pm 2}$. The degree four covering map $\P^1\times \P^1\lra Q_{\pm 2}$ realizing the latter is the degree $4$ morphism
$$
[x_1\colon x_2]\times [y_1\colon y_2]\lra \begin{bmatrix}
 (x_1y_2+x_2y_1)^2+x_1x_2y_1y_2 & (x_1y_1)^2 & \\
-(x_2y_2)^2 &  (x_1y_2+x_2y_1)^2-x_1x_2y_1y_2& \\
&& \pm (x_1y_2+x_2y_1)^2
\end{bmatrix}
$$
which does not look like an extension of the degree $2$ maps defining $\iota$. \end{remark}

Taking the naive limit of $\iota_{(t,s)}$ as $s\lra 0$ yields the constant map 
$$
\aligned
&\P^1\times \P^1 \lra Q_{\pm 2}, \\ 
&[x_1\colon x_2]\times [y_1\colon y_2]\lra \begin{bmatrix}
\pm(x_2y_1- x_1 y_2) &0 & \\
0 & \pm (x_2y_1- x_1y_2) &\\
&& x_2y_1-x_1y_2
\end{bmatrix}= \begin{bmatrix}
 I_2 & \\
&\pm 1
\end{bmatrix}.
\endaligned
$$
Therefore, the naive extension of the map $\iota$ in (\ref{iotacA_e}) to the fibers at $\pm2$ does not provide a parametrization of $Q_{\pm 2}$. 
The discussion below shows that one needs to pass to a blowup to obtain a non-trivial limit.\\

By Remark~\ref{change-sign_rmk}, the map $\iota \colon \P^1\times \P^1 \times \big(\ov\cA-\{(\pm 2,0)\}\big)\lra \ov\SL(2,\C)$
in (\ref{iotacA_e}) descends to a similarly denoted map 
$$
\iota \colon  \Big(\P^1\times \P^1 \times \big(\ov\cA-\{(\pm 2,0)\}\big)\Big)\big/\si\lra \ov\SL(2,\C).
$$
Furthermore, the latter lifts to a similarly denoted $\SL(2,\C)$-equivariant inclusion map into $M$ such that the following diagram commutes:
\bEq{iotaM_e}
\xymatrix{
 \Big(\P^1\times \P^1 \times \big(\ov\cA-\{(\pm 2,0)\}\big)\Big)\big/\si \ar[rrr]^{~~~~~~~~~~~\iota} \ar[d]&&& M \ar[d]^{\tn{tr}}\\
\big(\ov\cA-\{(\pm 2,0)\}\big)/\si=\P^1-\{\pm 2\}\ar [rrr] &&& \P^1\;.
}
\eEq
Here, by Proposition~\ref{Main_Prop}, the $\SL(2,\C)$-action on the domain of $\iota$ is induced by the diagonal M\"obius action on $\P^1\times\P^1$ and the trivial action on $\cA$. The action on the target is that of Lemma~\ref{action-on-M}.
The map $\iota$ in (\ref{iotaM_e}) defines a birational morphism $M'\dasharrow M$ between the compact quotient space
$$
M'\defeq \big(\P^1\times \P^1 \times \cA\big) /\si
$$
and $M$. The threefold $M'$ fibers over $\ov\cA/\si \cong \P^1$, has $\P^2$ double-fibers\footnote{Recall from Remark~\ref{Qpm2_rmk} that $(\P^1\times \P^1)/\si\cong \P^2$.} over $\pm 2$, and has A1-singularities along the image curves $\De_\pm$ of the $\si$-fixed loci $\De_{\P^1}\times \{\pm 2\}\subset \P^1\times \P^1 \times \cA$. The rational curves $\De_{\pm}$ are $\SL(2,\C)$-invariant conics inside the $\P^2$ double-fibers.\\

Blowing up $M'$ along $\De_\pm$ gives us a smooth threefold $M''\lra \P^1$ with reducible fibers $\P^2\cup \P(\cO\oplus \cO(2))$ over $\pm 2$. 
The $\SL(2,\C)$ action lifts to the blowup. Blowing down the $\P^2$ components gives us $M$. Therefore, we get an $\SL(2,\C)$-equivariant diagram 
$$
\xymatrix{
&M'' \ar[dl]\ar[dr]& \\
M'\ar@{-->}[rr]^{\iota}&&M
}
$$
describing the birational morphism $\iota$ such that $M''\lra M$ is blowup at the singular points of $Q_{\pm 2}\subset M$ and $M''\lra M'$ is the blowup along the singular locus $\De_{\pm}\subset M'$.\\

More explicitly, the contraction map $M''\lra M$ is obtained from $\iota$ in (\ref{iotaM_e}) by rescaling the domain  $\P^1\times \P^1$ of $\iota_{(t,s)}$ in the normal direction to $\De_{\P^1}$ as $s\lra 0$, that results in blowing up the domain $M'$ to $M''$ in the following sense.\\

In terms of affine coordinates $x=\frac{x_1}{x_2}$ and $y=\frac{y_1}{y_2}$, the map $\iota_{(t,s)}$ takes the form
$$
\iota_{(s,t)}(x,y)=
\begin{bmatrix}
t\,\frac{y-x}{2}- s\,\frac{x+y}{2} & s xy& \\
-s & t\,\frac{y-x}{2}+ s\,\frac{x+y}{2} & \\
&& y-x
\end{bmatrix}.
$$
Blowing up $\P^1\times \P^1\times \cA$ along $\De_{\P^1}\times (\pm 2, 0)$ corresponds to replacing the local coordinates 
$$
(u=\frac{x+y}{2},v=\frac{y-x}{2}, s)
$$
with 
$$
(u,v, s)\times [v'\colon s']\in \C^3\times \P^1, \quad \tn{with}~~vs'=sv'.
$$ 
On the chart $s'\neq 0$, we can put $s'=1$, $v=sv'$ and we get 
\bEq{limitiota_e}
\iota_{(s,t)}(u,v')=
\begin{bmatrix}
ts v'- s u & s (u^2-s^2 v'^2)& \\
-s & tsv' + s u & \\
&& 2sv'
\end{bmatrix}=
\begin{bmatrix}
t v'- u & (u^2-s^2 v'^2)& \\
-1 & tv' +  u & \\
&& 2v'
\end{bmatrix}.
\eEq
The limit of the latter as $s\lra 0$ is 
$$
\begin{bmatrix}
2v'- u & u^2& \\
-1 & 2v' +  u & \\
&& \pm 2v'
\end{bmatrix}\in Q_{\pm 2}\;.
$$
On the chart $v'\neq 0$, we can put $v'=1$, $s'v=s$ and we get 
 $$
\iota_{(s'v,t)}(u,v)=
\begin{bmatrix}
t v- s'v u & s'v (u^2-v^2)& \\
-s'v & tv + s'v u & \\
&& 2v
\end{bmatrix}=
\begin{bmatrix}
t- s'u & s'(u^2-v^2)& \\
-s' & t +  s'u & \\
&& 2
\end{bmatrix}\;.
$$
The limit of the latter as $s'\lra 0$ is 
$$
\begin{bmatrix}
I_2 \\
& \pm 1
\end{bmatrix}\in Q_{\pm 2}\;.
$$
We conclude that only after lifting from $M'$ to $M''$, the morphism $\iota$ in (\ref{iotaM_e}) defined on 
$$
M''|_{\P^1-\{\pm 2\}}\cong \big(\P^1\times \P^1 \times \big(\ov\cA-\{(\pm 2,0)\}\big)\big)/\si
$$
extends to
$$
\iota \colon M'' \lra M
$$
such that $\iota$ restricted to the $\P^2$ components is the constant map $\begin{bmatrix}
I_2 \\
& \pm 1
\end{bmatrix}$ and restricted to each $\P(\cO \oplus \cO(2))$ is the blowdown map onto $Q_{\pm 2}$. This explains Remark~\ref{mfM_rmk}. \\

Lemma~\ref{Fmft_lmm} follows from (\ref{Apqt}), (\ref{limitiota_e}), and comparison of line bundles on $M$, $M'$, and $M''$.

\subsection{Proof of Theorem~\ref{main_thm}}
In this section, we put the results proved in the previous sections together to complete the proof of Theorem~\ref{main_thm}.\\

Fix $\mf{b}=(b_1,\ldots,b_{n-1})\in \Z_+^{n-1}$ and let $\mf{a}=(b_1,b_1,\ldots,b_{n-1},b_{n-1})\in \Z_+^{2(n-1)}$ be the associated symmetric tuple as in Definition~\ref{Syma_dfn}.\\

In Section~\ref{RC}, we constructed a quasi projective relative compactification
$$
\Pi_n \colon \ov\cR^{\tn{rel}}_n\lra\C^n,
$$
that sits as a divisor in 
$$
\mf{M}_n=(M_\C^{n-1}\times \C)\sslash_{\cL'_{\mf{b}}}  \SL(2,\C).
$$
By Proposition~\ref{Main_Prop}, and with notation as in (\ref{Bcirc}), for 
$$
\mu=((t_1,s_1),\ldots,(t_{n-1},s_{n-1}), t_n)\in \cB^\circ,\quad \tn{with}~\mft=(t_1,\ldots,t_n)\in \C^n_\circ,
$$  
under the identification
\bEq{imu}
\iota_{\mu}\defeq \iota_{(t_1,s_1)}\times \cdots \times \iota_{(t_{n-1},s_{n-1})}\colon (\P^1)^{2n-2}\lra Q_{\mft}=Q_{t_1}\times \cdots \times Q_{t_{n-1}},
\eEq
the conjugation action of $\SL(2,\C)$ on $Q_{\mft}$ pulls back to the $\mu$-independent diagonal action of $\SL(2,\C)$ on $(\P^1)^{2(n-1)}$ by M\"obius transformations. 
Furthermore, pullback of the defining equation 
$$
 \tn{tr}(A_1\ldots A_{n-1})-t_n e_1\ldots e_n=0
$$
of $\wt{\cR}_{n,\mft}$ in (\ref{section-equation_e}) is of the form 
$$
\big(F_{\mu}(p_1,\ldots,p_{2n-2})=0\big)\subset(\P^1)^{2(n-1)},
$$
where 
$$
F_\mu \in H^0\big((\P^1)^{2(n-1)},\cO(1,\ldots,1)\big)^{\SL(2,\C)}
$$ 
is an $\SL(2,\C)$-invariant section of the line bundle $\cO(1,\ldots,1)\lra (\P^1)^{2(n-1)}$ depending algebraically (linearly) on 
$\mu$. By (\ref{change-sign_rmk}), the function 
$$
\mu \lra F_{\mu}
$$
is $\Gamma=\ll \si_1,\ldots,\si_{n-1}\rr$-equivariant.
The divisor $D_{\mf{a}}\subset X_{\mf{a}}=(\P^1)^{2(n-1)}\sslash_{\cL_{\mf{a}}} \SL(2,\C)$ is also the $\SL(2,\C)$-quotient of the zero set of an $\SL(2,\C)$-invariant section of $\cO(1,\ldots,1)$; i.e., the section $\prod_{i=1}^{n-1} \mf{s}_i$ in the terminology of (\ref{mfsi}). 
Therefore,  the map 
\bEq{Psimu_e}
\Psi\colon \cB^\circ \lra |D_{\mf{a}}|, \qquad \mu \lra F_{\mu}^{-1}(0)\sslash_{\cL_{\mf{a}}} \SL(2,\C)
\eEq
defines an $\Gamma$-equivariant morphism 
$
\Psi\colon \cB^\circ \lra |D_{\mf{a}}|
$
into the linear system of the divisor $D_{\mf{a}}\subset X_{\mf{a}}$ such that the divisor $\Psi\big(\mu\big)\subset X_{\mf{a}}$ is the compactification $\ov\cR_{n,\mft}$ in (\ref{ovcRK}). Since $\iota_{\mu}^{-1}( \partial Q_{\mf{t}})=\De$ (see (\ref{partialQ}) and the sentence before Remark~\ref{Qinfty}),
we have
$$
\cD_{n,\mf{t}}=\ov\cR_{n,\mft}-\cR_{n,\mft}=\Psi\big(\mu\big) \cap D_{\mf{a}}.
$$
Putting all the maps $\iota_{\mu}$ together, we get an $\SL(2,\C)$-equivariant and $\Gamma$-invariant morphism 
$$
\iota\colon (\P^1)^{2n-2}\times \cB^\circ \lra M_\C^{n-1} \times \C
$$
such that
\bEn
\item $\iota^* \big(\pi^*\cO(b_1)\boxtimes \cdots \boxtimes \pi^*\cO(b_{n-1})\big)=\cL_{\mf{a}}\boxtimes \cO_{\cB^\circ}$ 
\item\label{iso} $\big((\P^1)^{2n-2}\times \cB^\circ\big)/\Gamma \lra (M_\C^{n-1} \times \C)|_{\C^n_\circ} $  is an isomorphism, 
\eEn

%
In conclusion, the isomorphism in \ref{iso} descends to an isomorphism
$$
(X_{\mf{a}}\times \cB^\circ)/\Gamma \lra (\mf{M}_n)|_{ \C^n_\circ}$$
such that the pre-image of $\ov\cR^{\tn{rel}}_n$ is a divisor linearly equivalent to $(D_{\mf{a}}\times \cB^\circ)/\Gamma$, where the linear equivalence is given by the $\SL(2,\C)$ and $\Gamma$-invariant rational function 
$$
(\mu,p_1,\ldots,p_{2n-2})\lra \frac{F_{\mu}(p_1,\ldots,p_{2n-2})}{s_1\mfs_1(p_1,p_2)\cdots s_n\mfs_{n-1}(p_{2n-3},p_{2n-2})}
$$
on $(\P^1)^{2(n-1)}\times \cB^\circ$.\qed

\section{$\SL(2,\R)$ and $\tn{SU}(2)$ character varieties}\label{real}

Define the \textbf{real locus} $\cR^{\R}_{g,n}(\SL(2,\C))$ of $\cR_{g,n}(\SL(2,\C))$ to be the subset of (equivalence classes) of representations 
$$
\rho\colon \pi_1(\Si_{g,n})\lra \SL(2,\C)\quad \tn{s.t.}\quad \tn{tr}(\rho(\gamma))\!\in\! \R \qquad \forall~\gamma\in \pi_1(\Si_{g,n}).
$$ 
When $g=0$, we will simply denote $\cR^\R_{0,n}(\SL(2,\C))$ by $\cR^\R_{n}.$
The restriction of the projection map in (\ref{pin0_e}) to the latter is 
$$
\Pi_n\colon \cR^{\R}_{n}\lra \R^n.
$$
It is known (c.f. \cite{MS}) that
\bIt
\item $\cR^{\R}_{g,n}(\SL(2,\C))\cong \cR_{g,n}(\SL(2,\R))\cup \cR_{g,n}(\SU(2))$
\item each connected component of $\cR^{\R}_{g,n}(\SL(2,\C))$ is a symplectic variety with singularities (cf. \cite{G});
\item the $\SL(2,\R)$-component has finitely many components one of which is the Teichm\"uller space of $\Si_{g,n}$;
\item for certain values of $\mft\in (-2,2)^{n}$, the genus $0$ relative $\SL(2,\R)$-character variety
$$
\cR_{n,\mft}(\SL(2,\R))\defeq \Pi_n^{-1}(\mft)\cap \cR_{0,n}(\SL(2,\R))
$$ 
contains compact components isomorphic to $\C\P^{n-3}$ (see \cite{DeTh} and Remark~\ref{theta-to-t});
\item The $\cR_{g,n}(\SU(2))$ character variety is a closed symplectic manifold with singularities.
\eIt

In this section, we will only focus on the elliptic ($|\tn{tr}|<2$) elements in $\SL(2,\R)$ and the corresponding character variety 
\bEq{elliptic-ChVar}
\cR_{n}(\SL(2,\R))^{\tn{elliptic}}=\Pi_n^{-1}((-2,2)^n)\cap \cR_{0,n}(\SL(2,\R)).
\eEq

In what follows, first, we introduce the anti-holomorphic involutions $\tau$ and $\eta$ in Theorem~\ref{real-main_thm}. Then, we use the parametrization map (\ref{Apqt}) to identify the fixed loci of $\tau$ and $\eta$ on $\ov\cR_{n,\mft}\!\subset\! X_{\mf{a}}$ with $\ov\cR_{n,\mft}(\SL(2,\R))^{\tn{elliptic}}$ and $\ov\cR_{n,\mft}(\SU(2))$, respectively.

\subsection{Anti-holomorphic involutions}
There are two conjugacy classes of anti-holomorphic involutions on $\P^1$ characterized by their fixed loci.
The involution 
$$
\tau\colon \P^1[x_1\colon x_2]\lra\P^1[x_1\colon x_2], \quad [x_1\colon x_2]\lra [\ov{x}_1\colon \ov{x}_2]
$$
is the natural extension of the standard complex conjugation on $\C$ with $\tn{Fix}(\tau)\cong S^1$. The involution
$$
\eta\colon \P^1[x_1\colon x_2]\lra\P^1[x_1\colon x_2], \quad [x_1\colon x_2]\lra [\ov{x}_2\colon -\ov{x}_1],
$$
however, is the antipodal map on $S^2$ and has no fixed-point. 
These involutions naturally appear in real Gromov-Witten theory, Floer homology, and Mirror Symmetry. \\

For $c=$ $\tau$ and $\eta$, let 
$$
\wt{c}\colon \P^1\times \P^1\lra \P^1\times \P^1, \quad (p,q)\lra (c(q),c(p))
$$
denote the corresponding twisted anti-holomorphic involutions on $\P^1\times \P^1$. While
$$
\tn{Fix}(\wt\tau)\cong \P^1\quad \tn{and}\quad \tn{Fix}(\wt\eta)\cong \P^1,
$$
they are different by the fact that 
\bEq{tildeFixloci}
\tn{Fix}(\wt\tau)\cap \De_{\P^1}\cong S^1\quad \tn{and}\quad \tn{Fix}(\wt\eta)\cap \De_{\P^1}=\eset.
\eEq

\begin{lemma}
For $c=\tau$ or $\eta$, the anti-holomorphic involution
$$
\wt{c}_{n-1}\defeq \underbrace{\wt{c}\times \cdots \times \wt{c}}_{n-1~\tn{times}}\colon (\P^1)^{2(n-1)}\lra (\P^1)^{2(n-1)}
$$ 
preserves the $\SL(2,\C)$-orbits of the diagonal M\"obius action on $(\P^1)^{2(n-1)}$ as well as each divisor $\De_i$. Therefore, $\tau$ and $\eta$ induce similarly denoted anti-holomorphic involutions 
$$
\tau,\eta \colon X_{\mf{a}}=(\P^1)^{2(n-1)}\sslash_{\cL_{\mf{a}}}\SL(2,\C)\lra X_{\mf{a}}
$$
that preserve each component of $D_{\mf{a}}$. Moreover, $\tn{Fix}(\eta)\subset X_{\mf{a}}-D_{\mf{a}}$.
\begin{proof}
For every
$$
\rho=\begin{pmatrix}
a & b \\
c & d
\end{pmatrix}\in \SL(2,\C),
$$
we have 
$$
\rho \circ \wt\tau_{n-1}= \wt\tau_{n-1} \circ \ov\rho\qquad \tn{and}\qquad 
\rho \circ \wt\eta_{n-1}= \wt\eta_{n-1} \circ \ov\rho',
$$
when acting on $(\P^1)^{2n-2}$,
where 
$$
\ov\rho=\begin{pmatrix}
\ov{a} & \ov{b} \\
\ov{c} & \ov{d}
\end{pmatrix}\in \SL(2,\C) \quad \tn{and}\quad 
\ov\rho'=\begin{pmatrix}
\ov{d} & -\ov{c} \\
-\ov{b} & \ov{a}
\end{pmatrix}\in \SL(2,\C).
$$
This proves the first statement. It is clear that each divisor $\De_i$ is preserved by the two involutions. Finally, by the second identity in (\ref{tildeFixloci}), the fixed locus of the involution $\eta\colon X_{\mf{a}}\!\lra\! X_{\mf{a}}$
has no intersection with~$D_{\mf{a}}$.
\end{proof}
\end{lemma} 

\begin{remark}
The untwisted component-wise extension  
$$
(\P^1)^{2(n-1)}\lra (\P^1)^{2(n-1)}, \quad (p_1,\ldots,p_{2n-2})\lra (\tau(p_1),\ldots,\tau(p_{2n-2})),
$$
also induces an anti-holomorphic involution on each $X_{\mf{a}}$ that corresponds to (not so interesting) representations in the hyperbolic ($|\tn{tr}|>2$) range of $\SL(2,\R)$. 
\end{remark}

\subsection{Proof of Theorem~\ref{real-main_thm}}
For $t\in \R$, with $|t|<2$, we have that $s=\sqrt{t^2-4}=\mfi y$ is a pure imaginary number and we choose the choice satisfying $y\!>\!0$. Therefore, restricted to $(-2,2)^n \subset \C^n_\circ$, (\ref{imu}) and (\ref{Psimu_e}) simplify to
$$
\iota_{\mft} \colon (\P^1)^{2n-2}\lra Q_{\mft}, \qquad \forall~\mft\in(-2,2)^n  
$$
and 
$$
\Psi\colon(-2,2)^n \lra |D_{\mf{a}}|.
$$
To avoid writing square roots, in what follows, $s$ or $y$, as defined above, will appear in the formulas but they should be considered as functions of $t$.
\begin{lemma}
With notation as in (\ref{section-equation_e}), for $\mft\in(-2,2)^n$ the divisor 
$$
\big(F_{\mft}(p_1,\ldots,p_{2n-2})=0\big)\subset (\P^1)^{2n-2}
$$ 
is invariant under the actions of $\wt{\tau}_{n-1}$ and $\wt{\eta}_{n-1}$. Therefore, $\ov\cR_{n,\mft}$ is preserved by $\tau$ and $\eta$.
\end{lemma}

\bPf
For $t\in \R$, with $|t|<2$, we have $\beta_2=\ov{\beta_1}$ and 
$$
 \begin{bmatrix}
A(\wt\tau(p,q),t) & \\
&e(\wt\tau(p,q))
 \end{bmatrix}=
 \begin{bmatrix}
  \ov{\beta_1 y_2x_1}- \ov{\beta_2 y_1 x_2} & s \ov{x_1y_1} &\\
-s\ov{x_2y_2} & \ov{\beta_2 y_2x_1}- \ov{\beta_1 y_1x_2} &\\
 && -\ov{(x_2y_1-x_1y_2)}
 \end{bmatrix}.
$$
In other words, if 
$$
 \begin{bmatrix}
A(p,q,t) & \\
&e(p,q)
 \end{bmatrix}=
 \begin{bmatrix}
a & b &\\
c& d &\\
 && e 
 \end{bmatrix},
$$
then 
$$
 \begin{bmatrix}
A(\wt\tau(p,q),t) & \\
&e(\wt\tau(p,q))
 \end{bmatrix}=
 \begin{bmatrix}
\ov{a} & \ov{b} &\\
\ov{c}& \ov{d} &\\
 && \ov{e}
 \end{bmatrix}.
$$
On the other hand, if $\tn{tr}\big(A_1\ldots A_{n-1}\big)=t_n e_1 \ldots e_n$
then $
\tn{tr}\big(\ov{A}_1\ldots \ov{A}_{n-1}\big)=t_n \ov{e}_1 \ldots \ov{e}_n.
$
Therefore, $F_{\mft}(\wt\tau_{n-1}(p_1,\ldots,p_{2n-2}))=\ov{F_{\mft}(p_1,\ldots,p_{2n-2})}$ which proves the first claim.\\

The other case is similar but less trivial. We have
$$
 \begin{bmatrix}
A(\wt\eta(p,q),t) & \\
&e(\wt\eta(p,q))
 \end{bmatrix}=
 \begin{bmatrix}
  \ov{\beta_2 y_2x_1}- \ov{\beta_1 y_1 x_2} & \ov{-s x_2y_2} &\\
\ov{s x_1y_1} &  \ov{\beta_1 y_2 x_1}-\ov{\beta_2 y_1x_2} &\\
 && -\ov{(x_2y_1-x_1y_2)}
 \end{bmatrix}.
$$
In other words, if 
$$
 \begin{bmatrix}
A(p,q,t) & \\
&e(p,q)
 \end{bmatrix}=
 \begin{bmatrix}
a_{11} & a_{12} &\\
a_{21}& a_{22} &\\
 && e 
 \end{bmatrix},
$$
then 
$$
 \begin{bmatrix}
A(\wt\eta(p,q),t) & \\
&e(\wt\eta(p,q))
 \end{bmatrix}=
 \begin{bmatrix}
\ov{a}_{22} & \ov{a}_{21} &\\
\ov{a}_{12}& \ov{a}_{11} &\\
 && \ov{e}
 \end{bmatrix}.
$$
The matrices $A(p,q,t)$ and $A(\wt\eta(p,q),t) $ are related by the composition of complex conjugation and 
$$
\begin{pmatrix}
a_{11} & a_{12} \\
a_{21}& a_{22} 
 \end{pmatrix}\stackrel{f}{\lra} \begin{pmatrix}
a_{22} & a_{21} \\
a_{12}& a_{11} 
 \end{pmatrix}.
$$
Therefore, it remains to show that $\tn{tr}\big(f(A_1)\ldots f(A_{n-1})\big)=\tn{tr}\big(A_1\ldots A_{n-1}\big).$ 
If we write
$$
A_{i}=
\begin{pmatrix}
a_{i,11}& a_{i,12} \\
a_{i,21}& a_{i,22}\\
\end{pmatrix}, \qquad \forall~1\leq i\leq n-1,
$$
then
\bEq{expansion}
\tn{tr}\big(A_1\ldots A_{n-1}\big)=\sum_{j} \prod_{i=1}^{n-1} a_{i,j_{i}j_{i+1}} 
\eEq
where $j$ runs over sequences of the form $j=(j_1,\ldots,j_{n})\in \{1,2\}^n$ satisfying $j_1=j_n$.
For every such $j$, there is a dual $j$ denoted by $\ov{j}$ obtained by switching $1$s to $2$s and vice-versa in the sequence. The action $f$ on (\ref{expansion}) corresponds to switching  $j \lra \ov{j}$. The equation (\ref{expansion}) is clearly invariant under this operation.  Therefore, $\tn{tr}\big(f(A_1)\ldots f(A_{n-1})\big)=\tn{tr}\big(A_1\ldots A_{n-1}\big).$ 
\ePf

To finish the proof of Theorem~\ref{real-main_thm}, it remains to identify the fixed loci of $\tau$ and $\eta$ on $\ov\cR_{n,\mft}\!\subset\! X_{\mf{a}}$ with $\ov\cR_{n,\mft}(\SL(2,\R))^{\tn{elliptic}}$ and $\ov\cR_{n,\mft}(\SU(2))$, respectively.\\

For $t\in \R$, with $|t|<2$, and 
$$
(p,\tau(p))=[x_1\colon x_2]\times [\ov{x}_1\colon \ov{x}_2]\in \tn{Fix}(\wt\tau),
$$
the matrix assigned in (\ref{IsoP1P1}) to $(p,t)$ is the real matrix
\bEq{Atau_e}
\aligned
&\begin{bmatrix}
A_\tau(p,t) & \\
&e_\tau(p)
 \end{bmatrix}\defeq
 \begin{bmatrix}
A(p,\tau(p),t) & \\
&e(p,\tau(p))
 \end{bmatrix}\\
&=\begin{bmatrix}
 \beta_2 x_2\ov{x}_1-\beta_1 x_1 \ov{x}_2 & |x_1|^2 & \\
-|x_2|^2 & \beta_1 x_2\ov{x}_1-\beta_2 x_1\ov{x}_2 & \\
&& s^{-1}(x_2\ov{x}_1-x_1\ov{x}_2)
\end{bmatrix}\\
&=\begin{bmatrix}
-y\mf{Re}(x_2\ov{x}_1)+ t  \mf{Im}(x_2\ov{x}_1)& y|x_1|^2 & \\
-y |x_2|^2 & y \mf{Re}(x_2\ov{x}_1)+t \mf{Im}(x_2\ov{x}_1)& \\
&& 2 \mf{Im}(x_2\ov{x}_1)
\end{bmatrix}\in \ov\SL(2,\R) \cap Q_t.
\endaligned
\eEq
Away from the circle $(\mf{Im}(x_2\ov{x}_1)=0)=\tn{Fix}(\wt{\tau})\cap \De$, we can normalize (by dividing by $e_\tau$) to get 
$$
A'_\tau(p,t)\defeq \frac{1}{e_\tau(p)}A_\tau(p,t)=\frac{1}{2}\begin{pmatrix}
t-y \frac{\mf{Re}(x_2\ov{x}_1)}{\mf{Im}(x_2\ov{x}_1)}& \frac{y |x_1|^2}{\mf{Im}(x_2\ov{x}_1)} \\
\frac{-y |x_2|^2}{\mf{Im}(x_2\ov{x}_1)} &  t+y \frac{\mf{Re}(x_2\ov{x}_1)}{\mf{Im}(x_2\ov{x}_1)}\\
\end{pmatrix}\in \SL(2,\R)_{\tn{tr}=t}.
$$

\begin{remark}\label{theta-to-t}
In \cite{Ma}, relative elliptic $\SL(2,\R)$-character varieties are defined by fixing the rotation angles $(\theta_1,\ldots,\theta_n)$ of the $n$ matrices 
$$
A_1,\ldots,A_{n-1},A_n \defeq A_1\ldots A_{n-1}
$$
around their respective fixed points in the upper-half plane $\H$. In terms of the affine coordinate $z=\frac{x_1}{x_2}=a+ib$ in the upper-half plane $\H$ and the angle parameter $\theta=\theta(t)$ defined by 
\bEq{theta-to-t_e}
t=2\cos(\theta/2),\qquad 0< \theta  < 2\pi.
\eEq
We have 
$$
A'_\tau(z,\theta)=\frac{1}{b}\begin{pmatrix}
\cos(\theta/2) b- \sin(\theta/2) a& -\sin(\theta/2)(a^2+b^2) \\
\sin(\theta/2)& \cos(\theta/2) b+\sin(\theta/2)a\\
\end{pmatrix}
\in \SL(2,\R)_{\tn{tr}=2\cos(\theta/2)}.
$$
The fixed point of the action of $A'_\tau(z,\theta)$ on $\H$ is $z$ and the rotation angle of that around $z$ is $\theta$. Therefore, the real relative character varieties considered in \cite{Ma,DeTh} are the same as here. By \cite[Thm.~4]{DeTh}, if 
$$
2\pi(n-1)<\theta_1+\cdots+\theta_n  <2\pi n,
$$
then $\cR_{n,\mft}(\SL(2,\R))$ includes a closed component isomorphic to $\C\P^{n-3}$. 
\end{remark}

By~(\ref{Atau_e}), for $\mft\in (-2,2)^{2n-2}$, the function 
$\iota_{\mft} \colon (\P^1)^{2n-2}\lra Q_{\mft}$
maps $\tn{Fix}(\wt\tau_{n-1})$ to 
$$
Q_\mft \cap \ov\SL(2,\R)^{n-1}.
$$
By Remark~\ref{theta-to-t}, it is surjective. Therefore, $\tn{Fix}(\tau)\cap \ov\cR_{n,\mft} =\ov\cR_{n,\mft}(\SL(2,\R))^{\tn{elliptic}}$ is a compactification of $\cR_{n,\mft}(\SL(2,\R))^{\tn{elliptic}}$ obtained by attaching the unbounded components of $\cR_{n,\mft}(\SL(2,\R))^{\tn{elliptic}}$ along the limiting real hypersurfaces 
$$\ov\cR_{n,\mft}(\SL(2,\R))^{\tn{elliptic}}\cap D_{\mf{a}}.
$$

Similarly to (\ref{Atau_e}), for $t\in \R$ with $|t|<2$, and
$$
(p,\eta(p))=[x_1\colon x_2]\times [\ov{x}_2\colon -\ov{x}_1]\in \tn{Fix}(\wt\eta),
$$
the matrix assigned in (\ref{IsoP1P1}) to $(p,t)$ is the (normalized) matrix
\bEq{Aeta_e}
A_{\eta}(p,t)\defeq \frac{1}{e(p,\eta(p))}A(p,\eta(p),t)=\frac{1}{2}\begin{pmatrix}
 t + s\frac{|x_1|^2-|x_2|^2}{|x_1|^2+|x_2|^2}&\frac{2s x_1\ov{x_2}}{|x_2|^2+|x_1|^2}  \\
\frac{2s\ov{x_1}x_2}{|x_2|^2+|x_1|^2} &   t + s\frac{|x_2|^2-|x_1|^2}{|x_1|^2+|x_2|^2} 
\end{pmatrix}\in \SU(2)_{\tn{tr}=t}.
\eEq
As in Remark~\ref{theta-to-t}, in terms  $\theta=\theta(t)$ defined by (\ref{theta-to-t_e}), the matrix
$$
A_\eta(p,\theta)=  \frac{1}{|x_1|^2+|x_2|^2}  \begin{pmatrix}
 \tn{e}^{-\mfi \theta/2} |x_2|^2+ \tn{e}^{\mfi \theta/2} |x_1|^2 & 2\mfi \sin(\theta/2) x_1\ov{x_2} \\
-2\mfi \sin(\theta/2)\ov{x_1}x_2&  \tn{e}^{\mfi \theta/2}  |x_2|^2+ \tn{e}^{-\mfi \theta/2}  |x_1|^2
\end{pmatrix}\in \tn{SU}(2).
$$ 
is the unique matrix in $\SU(2)$ such that $(p,\eta(p))$ are the fixed points of the M\"obius action of $A_\eta(p,\theta)$ on $\P^1$ and $\theta$ is the rotation angle around the $(p,\eta(p))$-axis.\\

By~(\ref{Aeta_e}), for $\mft\in (-2,2)^{2n-2}$, the function 
$\iota_{\mft} \colon (\P^1)^{2n-2}\lra Q_{\mft}$
maps $\tn{Fix}(\wt\eta_{n-1})$ to 
$$
Q_\mft \cap \SU(2)^{n-1}
$$
By the last paragraph, it is surjective. Therefore, $\tn{Fix}(\eta)\cap \ov\cR_{n,\mft}=\tn{Fix}(\eta)\cap \cR_{n,\mft}=\cR_{n,\mft}(\SU(2))\subset\cR^\R_{n,\mft}.$ 

\begin{example}
We have
$$
A(\infty,\theta)=  \begin{bmatrix}
 \tn{e}^{\mfi \theta/2} & 0 \\
0&  \tn{e}^{-\mfi \theta/2}  
\end{bmatrix}\in \tn{SU}(2).
$$
\end{example}

Since $\tn{SU}(2)$ is compact, for each $\vartheta=(\theta_1,\ldots,\theta_n)\in (0,2\pi)^n$, the relative $\tn{SU}(2)$ character variety $\cR_{n,\vartheta}(\tn{SU}(2))$ is isomorphic to the compact geometric quotient
$$
\lrp{(F_\vartheta =0)\subset (\P^1)^{n-1}}/\tn{SU}(2) .
$$
Here, we have replaced the $t$ variables in the notation $\cR_{n,\mft}(\tn{SU}(2))$ with the equivalent $\theta$ variables. 
By fixing the first point to be $p_1=\infty$, we realize that
$$
\cR_{n,\vartheta}(\SU(2))=H_{\vartheta}/S^1, \qquad 
H_{\vartheta}\defeq \big(F_\vartheta(A(\infty,\theta_1),-) =0\big) \subset (\P^1)^{n-2}
$$
is the quotient of some $S^1$-invariant real hypersurface $H_{\vartheta}\subset (\P^1)^{n-2}$.
The hypersurface $H_{\vartheta}$ includes a fixed point of the diagonal $S^1$-action on $(\P^1)^{n-2}$ if and only if 
$$
\sum_{i=1}^{n} \ve_i \theta_i \equiv \pi \qquad \tn{mod}~2\pi
$$
for some choice of signs $(\ve_i)\in (\pm)^{n}$. Therefore, topology of $\cR_{n,\vartheta}(\tn{SU}(2))$ is fixed in the connected components of the complement of the hyperplanes above in $(0,2\pi)^n$.

\bibliographystyle{amsalpha}

\end{document}